\documentclass[twoside,leqno,10pt]{article}
\usepackage{graphics,ifthen}
\usepackage{latexsym}
\begin{document}
\input{latex.sty}
\input{references.sty}
\input epsf.sty
\def\ind{\stackrel{\mathrm{ind}}{\sim}}
\def\iid{\stackrel{\mathrm{iid}}{\sim}}
\def\Prodi{\mathop{{\lower9pt\hbox{\epsfxsize=15pt\epsfbox{pi.ps}}}}}
\def\prodi{\mathop{{\lower3pt\hbox{\epsfxsize=7pt\epsfbox{pi.ps}}}}}

\def\Definition{\stepcounter{definitionN}
    \Demo{Definition\hskip\smallindent\thedefinitionN}}
\def\EndDefinition{\EndDemo}
\def\Example#1{\Demo{Example [{\rm #1}]}}
\def\EndExample{\qed\EndDemo}
\def\Category#1{\centerline{\Heading #1}\rm}

\newcommand{\eps}{\epsilon}
\newcommand{\proof}{\noindent {\bf Proof:\ }}
\newcommand{\remarks}{\noindent {\bf Remarks:\ }}
\newcommand{\note}{\noindent {\bf Note:\ }}
\newcommand{\examp}{\noindent {\bf Example:\ }}

\newcommand{\ben}{\begin{enumerate}}
\newcommand{\een}{\end{enumerate}}
\newcommand{\bi}{\begin{itemize}}
\newcommand{\ei}{\end{itemize}}
\newcommand{\hp}{\hspace{.2in}}

\newtheorem{PCR}{Algorithm}
\newtheorem{lw}{Proposition 3.1, Lo and Weng (1989)}
\newtheorem{thm}{Theorem}[section]
\newtheorem{defin}{Definition}[section]
\newtheorem{prop}{Proposition}[section]
\newtheorem{lem}{Lemma}[section]
\newtheorem{cor}{Corollary}[section]
\newcommand{\rb}[1]{\raisebox{1.5ex}[0pt]{#1}}
\newcommand{\mc}{\multicolumn}
 
\newenvironment{lmm}[1]{\begin{lemma}\label{#1}}{\end{lemma}}
\newcommand{\Tinf}[2]{ T_{\inf}^{#1,#2} }
\newcommand{\Tsup}[2]{ T_{\sup}^{#1,#2} }
\newcommand {\besd}{ \mbox{BES$(\delta)$}}
\newcommand {\besdbr}{ \mbox{BES$(\delta)^{\br}$}}
\newcommand {\besfd}{ \mbox{BES$(4-\delta)$}}
\newcommand {\best}{ \mbox{BES$(3)$}}
\newcommand {\besta}{ \mbox{BES$(2 - 2 \al)$}}
\newcommand {\bestbr}{ \mbox{BES$(3)^{\br}$}}
\newcommand {\Pbly}{ P^\dagger }
\newcommand {\Pbr}{ P^{\rm br}}
\newcommand{\fst}{ F_* }
\newcommand{\mula}{\mu_\la}
\newcommand{\efst}{ E( \fst ^{\alpha/\gamma} ) }
\newcommand{\fell}{ {\cal F }_\ell }
\newcommand {\Ed}{ P^{(\delta)}}
\newcommand {\gi}{\eta}
\newcommand {\ph}{\phi}
\newcommand {\ala}{{\alpha,\lambda}}
\newcommand {\mual}{\mu_{\ala}}
\newcommand {\gal}{c_\alpha}
\newcommand {\tht}{\phi}
\newcommand {\ps}{\psi}
\newcommand {\deln}{\Delta_n}
\newcommand {\delj}{\Delta_j}
\newcommand {\delo}{\Delta_1}
\newcommand {\xxi}{\chi}
\newcommand {\phs}{\xxi}
\newcommand {\CF}{C_F}
\newcommand {\Vj}{V_j}
\newcommand {\Mj}{M_j}
\newcommand {\hM}{\widehat{M}}
\newcommand {\vj}{V_j}
\newcommand {\hvj}{\widehat{V}_j}
\newcommand {\bvj}{\bar{V}_j}
\newcommand {\bvo}{\bar{V}_1}
\newcommand {\bv}{\bar{V}}
\newcommand {\yj}{\bvj}
\newcommand {\yn}{\bar{V}_n}
\newcommand {\yo}{\bvo}
\newcommand {\fjt}{\widetilde{F}_j}
\newcommand {\La}{\Lambda}
\newcommand {\al}{\alpha}
\newcommand {\Lal}{\Lambda_\al}
\newcommand {\giv}{ \,|\,}
\newcommand {\epo}{ \varepsilon_0}
\newcommand {\epso}{ \varepsilon}
\newcommand {\epone}{ \varepsilon_1}
\newcommand {\eptwo}{ \varepsilon_2}
\newcommand {\epj}{ \varepsilon_j}
\newcommand {\epi}{ \varepsilon_i}
\newcommand {\epFi}{ \varepsilon_{F,i}}
\newcommand {\epsi}{ \tilde{\varepsilon}_k}
\newcommand {\tjs}{ T_{j}^*}
\newcommand {\tjsm}{ T_{j-1}^*}
\newcommand {\tks}{ T_{k}^*}
\newcommand {\tksm}{ T_{k-1}^*}
\newcommand{\bX}{{\bf X}}
\newcommand{\kl}{{\lambda}}
\newcommand{\Te}{{T_\kl}}
\newcommand{\GT}{{G_{T_\la}}}
\newcommand{\FF}{{\cal F}}
\newcommand{\var}{{\rm Var}}
\font\bb=msbm10
\def\bR{\hbox{\bb R}}
\newcommand {\BLY}{ \mbox{$B^\dagger$}}
\newcommand{\br}{\mbox{$\scriptstyle{\rm br}$}}
\newcommand{\ex}{\mbox{$\scriptstyle{\rm ex}$}}
\newcommand{\me}{\mbox{$\scriptstyle{\rm me}$}}
\newcommand{\birth}{\mbox{$\scriptstyle{\rm birth}$}}
\newcommand{\BB}{B^{\br}}
\newcommand{\Bex}{B^{\ex}}
\newcommand{\Bme}{B^{\me}}
\newcommand{\exj}{B^{\ex}_j}
\newcommand{\tb}{T^{(3)}}
\newcommand{\ts}{T^{*}}
\newcommand{\la}{\lambda}
\newcommand{\exo}{e_{(F > x)}}
\newcommand{\exone}{e^{(F > 1)}}
\newcommand{\ins}{\nu}
\newcommand{\nuF}{\nu_F^{\br}}
\newcommand{\nuV}{\nu_V^{\br}}
\newcommand{\nb}{N^{\birth}}
\newcommand{\nst}{N^{*}}
\newcommand{\nstm}{N^{*}(t-)}
\newcommand{\muinv}{\mu^{-1}}
\newcommand{\muu}{\mu}
\newcommand{\ggla}{\la^{-\gamma} \Gamma_{\alpha}}
\newcommand{\gamal}{\Gamma_{\alpha}}
\newcommand{\hal}{h_\alpha}
\newcommand{\KAL}{k_\alpha}
\newcommand{\kal}{K_\alpha}
\newcommand{\calph}{c_\alpha}
\newcommand{\Pal}{P_\alpha}
\newcommand{\Palt}{\widetilde{P}_\alpha}
\newcommand{\Eal}{E_\alpha}
\newcommand{\gamj}{\Gamma_j}
\newcommand{\gamn}{\Gamma_n}
\newcommand{\pdalth}{\mbox{PD}(\alpha,\theta) }
\newcommand{\pdalal}{\mbox{PD}(\alpha,\alpha) }
\newcommand{\pdalzero}{\mbox{PD}(\alpha,0) }
\newcommand{\hf}{ \mbox{${1 \over 2}$}}
\newcommand{\hl}{ - \mbox{${1 \over 2}$} \lambda^2}
\newcommand{\hnu}{\mbox{${1 \over 2}$} \nu^2}
\newcommand{\hfpi}{ \mbox{${\pi \over 2}$}  }
\newcommand{\BBt}{B_t^{\br}}
\newcommand{\mbrpj}{M_j^{\br +}}
\newcommand{\mbrpk}{M_k^{\br +}}
\newcommand{\mbrpo}{M_1^{\br +}}
\newcommand{\mbrnj}{M_j^{\br -}}
\newcommand{\mbrno}{M_1^{\br -}}
\newcommand{\mbrj}{M_j^{\br}}
\newcommand{\mjbr}{M_j^{\br}}
\newcommand{\vjbr}{V_j^{\br}}
\newcommand{\fbrj}{F_j^{\br}}
\newcommand{\tfbrj}{\widetilde{F}_j^{\br}}
\newcommand{\fbrjh}{\widehat{F}_j^{\br}}
\newcommand{\fbrjhm}{\widehat{F}_{j-1}^{\br}}
\newcommand{\fj}{F_j}
\newcommand{\vbrj}{V_j^{\br}}
\newcommand{\fbrjo}{F^{\br}_1}
\newcommand{\fbrjt}{F^{\br}_2}
\newcommand{\vbro}{V^{\br}_1}
\newcommand{\vbrt}{V^{\br}_2}
\newcommand{\vbrth}{V^{\br}_3}
\newcommand{\mbr}{M^{\br}}
\newcommand{\mbrp}{M^{{\br} +}}
\newcommand{\mbrn}{M^{{\br} -}}
\newcommand{\mbro}{M_1^{\br}}
\newcommand{\mbrt}{M_2^{\br}}
\newcommand{\fbro}{F_1^{\br}}
\newcommand{\fbrt}{F_2^{\br}}
\newcommand{\mbrpt}{M_2^{\br + }}
\newcommand{\mbrnt}{M_2^{\br - }}
\newcommand{\mbrnseq}{(M_j^{\br - })}
\newcommand{\mbrpseq}{(M_j^{\br + })}
\newcommand{\mbrseq}{(M_j^{\br})}
\newcommand{\tfbrseq}{(\widetilde{F}_j^{\br})}
\newcommand{\vbrseq}{(V_j^{\br})}
\newcommand{\fbrseq}{(F_j^{\br})}
\newcommand{\mbrseqst}{(M_j^{*})}
\newcommand{\mbrjst}{M_j^{*}}
\newcommand{\mex}{M ^{\ex}}
\def\qed{\mbox{\rule{0.5em}{0.5em}}}
\def\proof{\noindent{\bf Proof.\ \ }}
\def\note{\noindent{\bf Note.\ \ }}
\def\endpf{$\Box$}

\def\Beta{\text{Beta}}
\def\Dir{\text{Dirichlet}}
\def\DP{\text{DP}}
\def\P{{\bf p}}
\def\fhat{\widehat{f}}
\def\GA{\text{gamma}}
\def\ind{\stackrel{\mathrm{ind}}{\sim}}
\def\iid{\stackrel{\mathrm{iid}}{\sim}}
\def\K{{\bf K}}
\def\N{\text{N}}
\def\p{{\bf p}}
\def\U{{\bf U}}
\def\u{{\bf u}}
\def\w{{\bf w}}
\def\X{{\bf X}}
\newcommand{\reals}{{\rm I\!R}}
\newcommand{\PR}{{\rm I\!P}}
\def\Z{{\bf Z}}
\def\yy{{\mathcal Y}}
\def\rr{{\mathcal R}}
\def\BP{\text{beta}}
\def\ts{\tilde{t}}
\def\js{\tilde{J}}
\def\gs{\tilde{g}}
\def\fs{\tilde{f}}
\def\ys{\tilde{Y}}
\def\ps{\tilde{\mathcal {P}}}

\def\Report{James}
\def\Author{Poisson Process Calculus}
\pagestyle{myheadings}
\markboth{\Author}{\Report}
\thispagestyle{empty}

\bct\Heading
Poisson Process Partition Calculus with applications to\\
Exchangeable models and Bayesian Nonparametrics.
\lbk\lbk\smc
Lancelot F. James\footnote{
\indent\eightit AMS 2000 subject classifications. 
               \rm Primary 62G05; secondary 62F15.\\
\indent\eightit Keywords and phrases. 
                \rm  
          Bayesian Nonparametrics,           
          Brownian excursions,
          Chinese restaurant process,
          exchangeable partition probability functions,
          inhomogeneous Poisson process,
          L{\'e}vy-Cox models,    
          multiplicative intensity models,
          Neutral to the right processes, 
          generalised gamma processes,
          two-parameter Poisson-Dirichlet 
          }
\lbk\lbk

\BigSlant 
The Hong Kong University of Science and Technology\rm 
\lbk\lbk (\today) 
\ect 

\Quote
This article discusses the usage of  a partiton based Fubini calculus
for Poisson processes. 
The approach is an amplification of Bayesian
techniques developed in Lo and Weng for gamma/Dirichlet processes.
Applications to models are considered which all fall within an
inhomogeneous spatial extension of the size biased framework used in 
Perman, Pitman and Yor. Among some of the results;
an explicit partition based calculus is then developed for such models,
which also includes a series of important exponential change of measure
formula. These results are then applied to solve the mostly unknown
calculus for
spatial L{\'e}vy-Cox moving average models. The analysis then proceeds to
exploit a structural feature of a  scaling
operation which arises in
Brownian excursion theory. From this a series of new mixture
representations and posterior characterizations for large classes of 
random measures, including probability measures, are given.
These results are applied to yield new
results/identities related to the large class of
two-parameter Poisson-Dirichlet models. 
The results also yields easily perhaps the most
general and certainly quite informative characterizations of
extensions of the Markov-Krein correspondence exhibited by the
linear functionals of Dirichlet processes. 
This article then defines a natural extension of Doksum's Neutral to the
Right priors (NTR) to a spatial setting. NTR models are practically
synonymous with exponential functions of subordinators and arise in Bayesian
non-parametric survival models. It is shown that manipulation of the
exponential formulae makes what has been otherwise formidable analysis
transparent. Additional interesting results related to the Dirichlet
process and other measures are developed. Based on practical
considerations, computational procedures
which are extensions of the Chinese restaurant process are also developed.
\EndQuote

\baselineskip14pt


\rm
\tableofcontents
\section{Introduction}

This paper discusses the active usage of a Poisson process partition
based Fubini calculus to solve a variety of problems. That is, a
method which will be used to solve problems associated with large
classes of random
partitions ${\bf p}:=\{C_{1},\ldots, C_{n({\bf p})}\}$ of the integers
  $\{1,\ldots,n\}$. 
The method is
based on the formal statement of two results, which are known in
various levels of generality,
 concerning a Laplace functional change of measure and a partition
based Fubini representation. In terms of technique, this is an
amplification of the methods discussed in Lo and Weng (1989)[see also
Lo (1984)] for a
class of Bayesian nonparametric weighted gamma process mixture
models. The idea to choose
a Poisson process framework was based on suggestions from Jim Pitman.
The utility  of the approach is demonstrated by its application to a
suite of problems which are within
the general size-biased framework of Pitman, Perman and Yor (1992,
Section 4), with now a spatial inhomogeneous
component. Methodologically this article may be viewed as a treatment of
combinatorial stochastic processes from a Bayesian (infinite-dimensional
calculus) technical viewpoint. 

A key role
in the works of Lo (1984) and Lo and Weng (1989), is played by a
partition distribution on the integers $\{1, \ldots, n\}$ which is
a variant of Ewens sampling formula [see Ewens~(1972) and Antoniak~(1974)]
associated with the Poisson-Dirichlet partition distribution. One
particular feature is that posterior quantities written with respect to a
Blackwell and MacQueen (1973) urn scheme can be further simplified to
calculations which amount to sums over partitions ${\bf p}$ of 
$\{1,\ldots,n\}$. This makes ${\bf p}$ what I term a {\it separating
  class}. The methodology discussed there amounts to a
partition based Fubini calculus for Dirichlet[see Ferguson (1973, 1974),
Freedman (1963) and Fabius (1963)] and gamma processes. 
Pitman (1996), extends the description of the Blackwell-MacQueen sampling from the Dirichlet process to a large
class of species sampling random measures. Importantly, he develops
ideas surrounding the two-parameter Poisson-Dirichlet family of
distributions in a Bayesian context. This provides a nice
bridge to related work, where the two-parameter family appears,  in
for instance Pitman (1995a, b, 1997a, 1999), Pitman, Perman and Yor (1992), Pitman
and Yor (1992, 1997, 2001). Those works are non-Bayesian and center
around topics such as Brownian excursion
theory and Kingman's theory of partition structures as developed in
Aldous (1985) and Pitman (1995a). Returning to a Bayesian setting,
Ishwaran and James (2001a) recently develop the calculus for a
class of species sampling mixture models analogous to and extending the
model of Lo (1984), based on the work of Pitman (1995a, 1996). 

The interest here is to
extend these ideas to other random measures, not necessarily mixture
models, via more general partition structures. 
The question of how indeed to obtain
information via as yet possibly unknown partition structures related
to classes of random measures  suggests that one may
need considerable expertise in combinatorial calculations. Here, 
this issue is circumvented by usage of what is referred to as a {\it Poisson
Process Partition Calculus}.  Note that Poisson Palm calculus is
employed in Pitman, Perman and Yor (1992) and Pitman and Yor
(1992). [See also Fitzsimmons, Pitman and Yor (1992)]. The intersection with those results manifests itself in
section 5.

One may infer
from Lo (1984) and Lo and Weng (1989), that  Bayesian
infinite-dimensional calculus 
is a calculus based on the {\it disintegration} of joint product
structures on abstract spaces which exploits properties of partitions
${\bf p}$ or some other {\it separating class}. Examples of other {\it
  separating classes}, which will not be discussed, are the ${\bf s}$-path models found in Brunner
and Lo (1989)
and the classification based methods for generalised Dirichlet
stick-breaking models
discussed in Ishwaran and James (2001b).
Specifically these ideas are applied to classes of {\it boundedly finite} random measures, say
$\mu$ , on complete and separable spaces[see Daley and Vere-Jones
(1986)] which are linked to Poisson random measures. I will mention early that this is not synonymous with the
notion of stochastic integration which 
suggests simply a marginalization over the infinite-dimensional
component; although these types of ideas will play a role. Here
the primary interest is in the derivation and various
characterizations of the joint strucuture in terms of an
infinite-dimensional posterior law of $\mu$ and its marginal components. A key aspect of
disintegration of measures on Polish spaces is the availabilty of a well
defined Fubini's theorem. The notion of Bayesian models and
disintegrations is made quite clear in Le Cam (1986,
Chapter 12). The term Bayesian is meant primarily in terms of
technique.
The treatment of problems is more in line with a broader point
of view such as in Kingman (1975) rather than Ferguson (1973). Of course
I shall cover quite thoroughly models which arise in Bayesian
nonparametrics.
Readily accessible general disussions  
on disintegrations of measures may be found in Pollard
  (2001) and Kallenberg (1997).  See also Blackwell and Maitra (1984)
, Dellacherie and Meyer (1978), and Pachl~(1978). 
Additionally, Daley and Vere-Jones~(1988),
Matthes, Kerstan and Mecke~(1978) and Kallenberg~(1986) provide details
about Fubini's theorem for random measures cast 
within the language of Palm calculus. 

\subsection{Basic principles and motivation}
   For motivation  the typical but rich mixture model
setup is described. Let $K$ denote a non-negative integrable kernel on a complete and separable
space ${\mathcal X}\times{\mathcal Y}$ and let $\mu$ denote a boundedly
finite discrete random measure on ${\mathcal Y}$. A mixture model is defined as
follows
\Eq
f(X|\mu):=\int_{\mathcal Y}K(X|Y)\mu(dY).
\label{mixturemodel} 
\EndEq
for $X$ varying in ${\mathcal X}$. When $\mu$ is a Dirichlet process, a
two-parameter Poisson-Dirichlet process or a weighted gamma process then we are
in the setting described earlier . When $K$ is chosen to be a density such as a
normal kernel and $\mu$ is a probability measure then $f$ is a random
density. This presents one way to describe random measures  over 
spaces of densities and is analogous to the idea behind classical
density estimation where one convolves an empirical distribution function
with a kernel.  A similar construct holds for random hazard rates
which might be useful in models in a multiplicative intensity
setting. However, the fact that the  kernel may be specified rather
arbitrarily leads to the description of a large body of models which
appear in nonparametric statistics, spatial modelling and general
inverse problems. Letting $K(A_{ij},Y):=I\{Y\in A_{i,j}\}$ corresponds to
partially observed models such as interval censoring,right censoring,
double censoring described, from a fequentist viewpoint,
in Turnbull (1976) and Groeneboom and Wellner (1992). See Lindsay (1995)
and  Groeneboom (1996) for much more general models.
As discussed in Lo (1984)
and Lo and Weng(1989) one can induce specific shapes such as the class
of monotone densities or hazards via a uniform kernel or  completely
monotone models via mixtures of exponential kernels. In statistical
terms, when $\mu$ is a probability measure $P$, the general mixture model has an interpretation where $X|Y,P$ is
$K(X|Y)$ and $Y|P$ is missing information with distribution
$P$. In spatial
statistics, Wolpert and
Ickstadt (1998b) propose to specify $\mu$ as a L{\'e}vy random field
where ~\mref{mixturemodel} may represent the intensity of a Poisson
process. Brix (1999) proposes a class of generalized gamma shot-noise
processes which is a flexible class of such models. 
Such models are used rather than a raw L{\'e}vy-Cox process to introduce spatial
dependence. The term L{\'e}vy-moving average process has been used in
analogy to Gaussian moving averages where an example is 
 a model of the form $\int K(x-y)\mu(dy)$ or more specifically 
\Eq
\int_{0}^{t}{\mbox e}^{(s-t)}\mu(ds)
\EndEq
which is reminiscent of a stationary Ornstein-Uhlenbeck process. Barndorff-Nielsen
and Shepard (2001) propose the usage of non-Gaussian
Ornstein-Uhlenbeck processes which also can be viewed as a mixture
representation where $\mu$ is a L{\'e}vy process. See also Le Cam (1961)
for an early mathematical discussion of random measures and shot-noise
models.
While such models are
indeed rich in terms of flexibility and diversity in terms of
applications, many questions remain open about their properties.
  Suppose that now one has the joint product measure of $\{X_1,\ldots,X_n,\mu\}$
\Eq
{\mathcal P}{(d\mu)}\prod_{i=1}^{n}\int_{\mathcal Y}K(X_i|Y_i)\mu(dY_i)
\EndEq
which arises in a variety of contexts. The quantity above represents
one possible disintegration of the joint strucuture
$\{X_1,\ldots,X_n,\mu\}$. Most often direct evaluation of
$\{X_1,\ldots,X_{n},\mu\}$ is not simple or practically
implementable.  Stripping away the kernel $K$ one is left with the
joint product structure,  $\{Y_{1},\ldots, Y_{n};\mu\}$
which due to versatility of an available Fubini's theorem becomes the
main object of interest. That is 
knowledge of this structure reduces the problem of the mixture model
above to a special case of a cadre of possibilities. Hence the goal is to
find the following disintegration
\Eq
{\mathcal P}(d\mu)\prod_{i=1}^{n}\mu(dY_i):={\mathcal P}(d\mu|{\bf Y})M_{\mu}(dY_1,\ldots, dY_n)
\EndEq 
where $M_{\mu}$ is a possibly sigma-finite joint moment measure of
${\bf Y}$ and
${\mathcal P}(d\mu|{\bf Y})$ can be thought of as the posterior
distribution of ${\mu|{\bf Y}}$. However, the result is best understood
and its utility is revealed by an equivalent statement via Fubini's
theorem, 
\Eq
\int g({\bf Y},{\mu})\prod_{i=1}^{n}\mu(dY_i){\mathcal P}(d\mu)
:=\int g({\bf Y},{\mu}){\mathcal
P}(d\mu|{\bf Y})M_{\mu}(dY_1,\ldots, dY_n)
\label{Fubini}
\EndEq
for $g$ an integrable or positive function. As is certainly known it
is sufficient to check for $g$ specified to be an indicator of appropriate cylinder
sets or other characterizing function.  The structure $M_{\mu}(dY_1,\ldots, dY_n)$ is an
urn-type structure that can be generated sequentially via conditional
moment measures. If $\mu$ is a probability measure the  notion of
conditional moment measures is synonymous with the notion of 
Bayesian prediction rules.  However the exchangeable urn structure
becomes a bit un-wieldly and one seeks a further disintegration of ${\bf
  Y}$. A natural one is based on the often quite informative decomposition  ${\bf Y}:={({\bf Y}^{*},{\bf p})}$ where ${\bf Y}^{*}=\{Y_{1},\ldots, Y_{n({\bf
    p})}\}$ are the unique random values given a random partition ${\bf
p}$ of the integers $\{1,\ldots, n\}$. In other words the joint (sigma-finite)
measure admits a disintegration in terms of a conditional measure on its
unique values  and a measure on ${\bf p}$. Neither of which need be a
proper probability measure.

The main point boils down to the
following basic principles; suppose that a random
measure~$\mu^{*}$ is some
function of $\mu$ , i.e. $\mu^{*}=g(\mu)$. Then via~\mref{Fubini} its 
posterior law, marginal law and partition strucuture can all be
derived from those corresponding aspects of $\mu$. 
This is of course provided that one has
explicit information about $\mu$. The utility of such a procedure for
$\mu$ is then amplified by its richness. That is, a measure of how
many interesting processes $\mu^{*}$ can it capture.
The interesting aspect of
this is that the measure of richness of $\mu$ must correspond to the 
simplicity of its posterior laws, marginal moment and partition
strucutures, while still being informative.
The
structures must indeed act in a way like canonical basis functions. In
other words complex structures can be derived from simple ones.
The Poisson random measure, $N$,  emerges as a natural candidate given its
prevalence in various theories of random measures and its basic
connections (via the Poisson random variable and Bell's number) to
random partitions of the integers. Albeit there is a duality
to an approach using combinatorial arguments, an exploitation of the Poisson random
measure analogue of ~\mref{Fubini}, with a further partition
disintegration, allows one
to proceed in a pure framework of disintegration of measures to
directly derive many aspects of large classes of measures $\mu$.

The notation $\sum_{\bf
  p}$ will be used to denote the sum over all partitions of the
integers $\{1,\ldots,n\}$. As is well known[see Rota (1964)], this sum is
equivalent to Bell's number. For papers which discuss the natural
relationships of the Poisson process/random variable to
partitions, see for instance Constantine and Savits (1994), Pitman
(1997b), Constantine (1999) and Di Nardo and Senato (2001). The papers
by Constantine and Savits (1994) and Constantine (1999), and references
therein, are certainly
related to this one. Constantine and Savits (1994) discuss methods to evaluate
identity/moment formulae for compound Poisson processes via Faa Di
Bruno's formula. Certainly one can infer from Theorem 2.1 of Constantine
and Savits (1994) that many of the formulae here, expressed in terms
of $\sum_{\bf p}$, can be re-expressed in tems of infinite-sum notation
related to Dobinski's formula or more obviously cycle notation. 

\subsection{Notation and preliminaries}

Again let ${\p}=\{C_{1},\ldots, C_{n(\p )}\}$ denote a partition of size $
n(\p)$ of the integers $\{1,\ldots ,n\}$, let $e_{j,n}$ denote the
cardinality of each cell $C_{j}$ for  $j=1,\ldots
,n(\mathbf{p})$. This partition structure is related to a description of
general analogues of a Chinese restaurant scheme to generate partitions
described in terms of a sequential seating of customers[see
Aldous~(1985), Pitman~(1996) and Kerov~(1998)]. The results will be closely connected
to such a structure generated from the {\it exchangeable partition
probability function} (EPPF) (partition distribution). See Pitman (1995a,b, 1996)
for a thorough description of the EPPF concept. Additionally, for
$r>1$, let $\p_r=\{C_{1,r},\ldots,C_{n(\p_r),r}\}$ denote a partition of $\{1,2\ldots,r\}$, where $C_{i,r}$ denotes the
current configuration of table $i$ after $r$ customers have been
seated and $e_{i,r}$ denotes the number of customers seated at
$C_{i,r}$ .  The partition $\p_{r+1}$ then denotes the (updated) one
step larger partition on $\{1,2,\ldots,r+1\}$. 

Now specific notation is given for models which shall be looked at in
some detail. That is, for the two-parameter Poisson-Dirichlet
models and closely connected generalised gamma family of random
measures. In addition notation is given for a general spatial variation of the
Beta process of Hjort (1990). 
First, we  briefly describe the two-parameter Poisson-Dirichlet 
class of models. See Pitman and Yor (1997) and Pitman (1996) for more
details. Let $(Z_i)$ denote a collection of iid random variables
whose distribution is a diffuse probability measure $H$ and
independently of $(Z_i)$, let $(P_i)$ denote a collection of ranked
probabilities which sum to one and have a two-parameter
Poisson-Dirichlet distribution denoted as $PD(\alpha, \theta)$
with parameter values $0\leq\alpha<1$ and $\theta>-\alpha$. The
corresponding random probability measure has a representation, 
\Eq
P_{\alpha,\theta}(\cdot):=\sum_{i=1}^{\infty}P_{i}\delta_{Z_{i}}(\cdot).
\EndEq
The law of $P_{\alpha,\theta}$, denoted ${\mathcal
  PD}_{\alpha,\theta}(dP|H)$ is uniquely associated with its 
prediction rule and {\it exchangeable partition probability function}
(EPPF) given as, 
\Eq
\PP\{Y_{n+1}\in\cdot\;|Y_1,\ldots,Y_n\} =
{\frac{\theta +n(\p)\alpha}{\theta +n}}H(\cdot) +
\sum_{j=1}^{n(\p)}{\frac{(e_{j,n}-\alpha)}{\theta +n}}\,\d_{Y_j^*}(\cdot),
\label{prediction.rule}
\EndEq  
and,
\Eq
{PD}(\p|\alpha, \theta) = 
\frac{\(\prod_{j=1}^{n({\bf p})-1}(\theta+\alpha\,j)\)
\(\prod_{j=1}^{n(\p)}\Gamma(e_j-\alpha)\)}
{\prod_{j=1}^{n-1}(\theta+j)}.
\label{twoparameterPY.EPPF}
\EndEq
The extreme cases are the normalised stable law process $P_{\alpha,0}$ and the
Dirichlet process, $P_{0,\theta}$ with shape parameter $\theta H$.  The
laws of the $(P_i)$ and 
Ewens sampling EPPF formula for the Dirichlet process are denoted as
$PD(\theta)$ and $PD({\bf p}|\theta)$. Similarly for the stable process
write $PD(\alpha)$ and $PD({\bf p}|\alpha)$. Now for $b\ge 0$ the rich family of generalized
gamma random measures [see Brix (1999)] is generated by the L{\'evy} measure
\Eq 
\rho_{\alpha, b}{(ds)}={\frac{1}{\Gamma(1-\alpha)}}s^{-\alpha-1}{\mbox e}^{-bs}ds,
\label{sub}
\EndEq
which includes the stable law subordinator, $b=0$, gamma processs
subordinator $\alpha=0$, and the inverse-Gaussian law, $\alpha=.5$, $b>0$,
among
others. The notation $\rho_{\alpha}:=\rho_{\alpha,0}$ will be reserved
for the stable law and the choice $\theta\rho_{0,1}$ will be used to
generate the Dirichlet process family of models.
The general subordinator has increments which have a distribution
belonging to an exponential family of distributions with a power
variance function, introduced by Tweedie (1984) and further discussed
in Hougaard~(1986), Bar-Lev and Enis (1996),  and Jorgensen (1997). See
K{\"{u}}chler and Sorensen (1997). Classes
of compound Poisson process models based on this distribution are
discussed in Aalen (1992), Lee and Whitmore (1993) and Hougaard, Lee
and Whitmore (1997). 
Lastly, a spatial version of Hjort's~(1990) (two-parameter) Beta
process corresponds to the L{\'e}vy process generated
by the inhomogeneous L{\'e}vy measure, 
\Eq
u^{-1}(1-u)^{c(s)-1}duA_{0}(ds,dx)
\label{Hjort}
\EndEq
for $(u,s,x)\in(0,1]\times(0,\infty)\times{\mathcal X}$. The quantity
$c(s)$ is a decreasing function on $(0,\infty)$ and $A_{0}$ is a hazard
measure.
The symbols ${\mathcal G}(a,b), {\beta}(a,b)$ will be used to denote
gamma and beta random variables respectively. ${\mathcal G}(dx|a,b)$
will denote a gamma density. 

Some other 
references connected to the Poisson-Dirichlet family, not mentioned
later, include McCloskey~(1965), Engen~(1978), Carlton~(1999),
Donnelly and Tavar{\'e} (1987),
Gyllenberg and Koski (2001). See the article by  Ewens and Tavar{\'e}
(1997)  for a discussion of the wide applicability of the
two-parameter model. See also Pitman (1995b), which will be referenced later.

\section{Poisson Process Partition Calculus}
Let $N$ denote an inhomogeneous Poisson process (measure) on a complete and
separable space ${\mathcal X}$ with (diffuse) mean measure $\nu(\cdot)$. 
That is, the Laplace functional of $N$ is of the form 
\Array
{\mathcal {L}}_{N}(f|\nu)&=&\int_{\mathcal{M}}\exp\left\{-\int_{
\mathcal{X}}f(x)N(dx)\right\}{\mathcal{P}}(dN|\nu)\\
&=&
\exp \(-\int_{\mathcal{X}}(1-{\mbox e}^{-f(x)})\nu (dx)\)
\EndArray
for non-negative functions $f\in BM(\mathcal{X})$ on $(\mathcal{X},\mathcal{B
}(\mathcal{X}))$ where $BM(\mathcal{X})$ denotes the collection of
measurable functions of bounded support on $\mathcal{X}$. See Daley and
Vere-Jones (1988) for a description of these concepts. For brevity we
use the shorthand notation of the type, 
$$
{\mbox e}^{-N(f)}=\exp\left\{-\int_{
\mathcal{X}}f(x)N(dx)\right\}.
$$
The exposition
of this paper centers around the utilization of disintegration results
related to the joint measure
\Eq
{\mathcal{P}}(dN|\nu){\prod_{i=1}^{n}N(dX_i)},
\label{jointPoisson}
\EndEq
where~\mref{jointPoisson} 
represents a disintegration of the joint product measure of 
$\{X_1,\ldots,X_n,N\}$. Moreover, the  collection {\bf X}=$\{X_1,\ldots,X_n\}$
can be considered as conditionally independent given $N$. However
importantly once integration is done over N the collection {\bf X} will
usually consist of tied values. It follows that one can always
represent ${\bf X}={({\bf X}^{*},{\bf p})}$ where ${\bf
  X}^{*}=\{X^{*}_{1},\ldots,X^{*}_{n({\bf p})}\}$ denotes the unique values
and {\bf p} dictates which variables are equal according to the
relationship $X_{i}=X^{*}_{j}$ if and only if $i\in C_j$. The main
purpose of this section is to describe two results  concerning the
Poisson process which are fashioned as {\it tools} to be tailor-made
to solve a variety of problems in an expeditious manner.

\subsection{Basic tools}
First an (exponential) change of measure or disintegration formulae
  based on Laplace functionals is given below. Such an operation is
  commonly called {\it exponential tilting}.  

\begin{lem}
For non-negative functions $f\in BM(\mathcal{X})$ on $(\mathcal{X}$, $
\mathcal{B}(\mathcal{X)})$ and $g$ on $(\mathcal{M},\mathcal{B}(\mathcal{M}))
$
$$
\int_{\mathcal{M}}g(N){\mbox e}^{-N(f)}\mathcal{P}(dN|\nu) \\
={\mathcal {L}}_{N}(f|\nu)
\int_{\mathcal{M}}g(N){\mathcal{P}}(dN|{\mbox e}^{-f}\nu),
$$ 
where ${\mathcal{P}}(dN|{\mbox e}^{-f}\nu)$ is the law of a Poisson Process with intensity
$$
{\mbox e}^{-f(x)}\nu(dx).
$$ 
In other words the following absolute continuity result holds,  
\Eq
{\mbox e}^{-N(f)}\mathcal{P}(dN|\nu)={\mathcal {L}}_{N}(f|\nu)
{\mathcal{P}}(dN|{\mbox e}^{-f}\nu).
\label{poissontilt}
\EndEq

\end{lem}
\Proof
By the unicity of of Laplace functionals for random measures on
$\mathcal{X}$ it suffices to check this result for the case
$g(N)={\mbox e}^{-N(h)}$.
Thus it follows that,
$$
\int_{\mathcal{M}}{\mbox e}^{-N(f+h)}\mathcal{P}(dN|\nu) \\
={\mathcal{L}}_{N}(f|\nu)
\int_{\mathcal{M}}{\mbox e}^{-N(h)}
{\mathcal{P}}_{f}(dN)
$$ 
where for the time being ${\mathcal P}_{f}$ denotes some law on $N$.
Simple algebra shows that 
$$
\int_{\mathcal{M}}{\mbox e}^{-N(h)}
{\mathcal{P}}_{f}(dN)
=\frac{{\mathcal{L}}_{N}(f+h|\nu)}{{\mathcal{L}}_{N}(f|\nu)}
$$
and hence ${\mathcal P}_{f}(dN):={\mathcal{P}}(dN|{\mbox e}^{-f}\nu)$
which concludes the result.
\EndProof

\Remark
Lemma 2.1 is a simple functional extension, mod the Gaussian
and drift component, of an analogous result for
L{\'e}vy processes on ${\mathcal R}$ or more generally ${\mathcal R}^{d}$
which may be found in K{\"{u}}chler and Sorenson (1997)
Proposition  2.1.3.  The utility of Lemma 2.1 will be demonstrated throughout.
Poisson processes with laws described by ${\mathcal{P}}(dN|{\mbox e}^{-f}\nu)$
can be found in Pitman and Yor~(1992)[See Section 5 of this manuscript].
\EndRemark

\Remark
Naturally Lemma 2.1 extends to the
following somewhat more vague generalisation; Suppose that ${\mu }$ is a
random measure with Laplace functional ${\mathcal L}_{\mu}$, then given
the setup in Lemma 2.1 with $\mu$ in place of $N$, 
 ${\mbox e}^{-\mu(f)}\mathcal{P}(d
\mu)={\mathcal {L}}_{\mu}(f){\mathcal P}_{f}(d\mu)$ where ${\mathcal
  P}_{f}$ is characterized by its Laplace functional 
$$
\int_{\mathcal{M}}{\mbox e}^{-\mu(h)}
{\mathcal{P}}_{f}(d\mu)
=\frac{{\mathcal{L}}_{\mu}(f+h)}{{\mathcal{L}}_{\mu}(f)}.
$$
One can replace the argument with characteristic functionals.
\EndRemark

Results which identify the disintegration of~\mref{jointPoisson}
in terms of the posterior
distribution of the Poisson process and the marginal joint measure
\Eq
M(dX_1,\ldots,dX_n)=\int_{\mathcal{M}}
  \[{\prod_{i=1}^{n}N(dX_i)}\]{\mathcal{P}}(dN|\nu)
\label{Pmoment}
\EndEq
are well known in the literature via Palm calculus for Poisson
processes. The quantity $M$ in~\mref{Pmoment} is also known as the joint moment
measure. These existing results are customized in Lemma 2.2 below
where emphasis is placed on the partition structure. 

\begin{lem}  
Let $g$ be a non-negative or integrable function on
${{\mathcal{X}}^{n}\times{\mathcal M}}$, then
for each $n\geq 1$,
\Eq
\int_{{\mathcal M}\times{\mathcal{X}}^{n}} g(\mathbf{x},N)
{\prod_{i=1}^{n}N(dx_i)}{\mathcal{P}}(dN|\nu)
=\sum_{{\bf
    p}}\int_{{\mathcal{X}}^{n({\bf p})}}\[\int_{{\mathcal M}}g({\bf X}^{*},{\bf p},N+{\sum_{j=1}^{n(\p)}\delta_{X^{*}_j}){\mathcal{P}}(dN|\nu)}\]
\prod_{j=1}^{n(\bf p)}\nu(dX^{*}_{j}),
\label{FubiniP}
\EndEq
with 
\Eq
\int_{{\mathcal M}}g({\bf X}^{*},{\bf
  p},N+{\sum_{j=1}^{n(\p)}\delta_{X^{*}_{j}}}) {\mathcal{P}}(dN|\nu)
=\int_{{\mathcal M}}g({\bf X}^{*},{\bf p},N){\mathcal{P}}(dN|\nu,{\bf X}).
\EndEq
The moment measure $M$ is also expressible via the conditional moment
measures as, 
\Eq
M(d{\bf X})=\nu(dX_1)\prod_{i=2}^{n}
\[\nu(dX_{i})+
\sum_{j=1}^{n({\bf p}_{i-1})}\delta_{X^{*}_{j}}(dX_{i})\]
\EndEq
\end{lem}

\Remark
It is of course true that Lemma 2.2 is not 
entirely novel. However, the partition representation that is used is certainly
not readily seen in the literature. Moreover, it has been tailor made
to assume its present purpose as a general tool. 
One way to deduce the partition
representation is to examine carefully Daley and Vere-Jones
(1988),[equation (5.517), Lemma 5.2.VI, and the discussion on page
192]. A simple minded but informative approach is to simply refer back
to the case of Poisson random variables. 
For clarity and also to showcase what is believed to be 
interesting side results involving partly Lemma 2.1 we prove this result in its
entirety in the next section using alternate means.
\EndRemark

\subsection{Supporting results}
Note that Lemma 2.1 implies the following result for each bounded set $B$,
\Eq
\int_{\mathcal{M}}N(B)\exp \left\{-\int_{\mathcal{X}}f(x)N(dx)\right\}\mathcal{P}(dN|\nu) \\
={\mathcal {L}}_{N}(f|\nu)
\int_{B}{\mbox e}^{-f(X)}\nu(dx).
\label{cobra}
\EndEq

This is reminiscent of the expression which appears in Lemma
10.6 in Kallenberg (1986) and perhaps more clearly in Proposition
12.1.V in Daley and Vere-Jones (1998). That
is, the expression~\mref{cobra} identifies the conditional Laplace
functional of the Poisson process given one observation[see Daley and Vere-Jones (1988), p. 458]. A point to
note is that in contrast to Daley and Vere-Jones (1988) Proposition
12.I.V, this result is not obtained  by taking derivatives. This
suggests that Lemma 2.1 can be used repeatedly to obtain the conditional
Laplace functional given n observations. Thus providing an alternative
to an argument using repeatedly say Lemma 12.1.V. 
The general dual of Lemma 12.1.V. can be deduced from Remark 2
as follows; Suppose that $\mu$ is an  random measure(as in Remark 2) with finite 1st
moment measure, say $m_{\mu}$, then 
\Eq
\int_{\mathcal{M}}\mu(B){\mbox e}^{-\mu(f)}{\mathcal{P}}(d\mu) \\
={\mathcal {L}}_{\mu}(f)
\int_{B}\int_{\mathcal M}\mu(dx)P_{f}(d\mu):={\mathcal {L}}_{\mu}(f)
\int_{B}r(f|x)m_{\mu}(dx)
\label{derivative2}
\EndEq
for some function $r$ determined by the second expression. That is, the
evaluation of $\int_{B}\int_{\mathcal M}\mu(dx)P_{f}(d\mu)$.
Hence the conditional Laplace functional
of $\mu|x$ is 
$$
{\mathcal L}_{\mu}(f|x):={\mathcal L}_{\mu}(f)r(f|x).
$$
This general form can be applied repeatedly to (conditional) random
measures $\mu|x_1,\ldots, x_{i}$ etc, where the requirement is
the existence of a finite conditional measure
$m_{\mu}(\cdot|x_1,\ldots,x_{i})$. All such results can be deduced
from an argument similar to what is used in Proposition 2.1 below.

\Remark
A result for general $\mu$ is quite applicable. As an  example
consider finite random mixtures of infinitely divisible random variables. 
That is,
\Eq
\sum_{k=1}^{m}W_{k,m}\delta_{Z_{k}}
\EndEq
where $W_{k,m}$ are iid infinitely divisible random variables and $Z_{k}$
are iid random variables. Such models 
can be used as approximations to many of the models discussed
here. The emphasis on the change of measure interpretation should
also prove useful. In fact, for such an infinitely divisible class the
results in Section 3 apply with small modification.
\EndRemark

\Remark
James (2001a, b) using the analogue of Lemma 2.1 for
weighted gamma and generalised weighted gamma process obtained their
posterior characterizations in this manner without any specific
mention of Poisson proceses. The idea for this approach is based on an
extension of the arguments in Lo and Weng (1989). In Section 3 it is
shown that Lemma 2.1 actually implies these analogues. 
\EndRemark

\begin{prop}
Lemma 2.1 implies that the conditional Laplace functional of $N|X_1,\ldots,X_n$
based on the model is, 
\Eq
{\mathcal {L}}_{N}(f|\nu)\[\prod_{j=1}^{n({\bf p})}{\mbox e}^{-f(X^{*}_j)}\]
\label{conditionalLaplace}
\EndEq 
\end{prop}
\Proof
The result proceeds by induction. Let $s,f\in
BM({\mathcal X})$ and choose
$g(N)=\int_{\mathcal {X}}s(v)N(dv)$. The case for $n=1$ follows from~\mref{cobra}.
For general $n=r$ it follows from Lemma 2.1 that the conditional Laplace
functional of $N$ given $({\bf X}_{r}, X_{r+1})$ is determined by the 
expression 
$$
{\mathcal {L}}_{N}(f|\nu)\[\prod_{j=1}^{n(\p_{r})}{\mbox e}^{-f(X^{*}_j)}\]
\[{\int_{\mathcal {X}}s(x_{r+1}){\mbox e}^{-f(X_{r+1})}
\nu(dX_{r+1})}+\sum_{j=1}^{n({\bf p}_{r})}s({X^{*}_{j}})\]
$$
Now define a function $t(X_{r+1})$ to be  ${\mbox e}^{-f(X_{r+1})}$ if $X_{r+1}$ is not equal to
$\{X^{*}_1,\ldots,X^{*}_{n({\bf p}_{r})}\}$ and is set to be one
otherwise. Then, 
$$
{\int_{\mathcal {X}}s(x_{r+1})t(X_{r+1})
\[\nu(dX_{r+1})+
\sum_{j=1}^{n({\bf p}_{r})}\delta_{X^{*}_{j}}(dX_{r+1})\]}
=
{\int_{\mathcal {X}}s(x_{r+1}){\mbox e}^{-f(X_{r+1})}
\nu(dX_{r+1})}+\sum_{j=1}^{n({\bf p}_{r})}s({X^{*}_{j}}).
$$
Hence,
the conditional Laplace functional is,
\Eq
{\mathcal {L}}_{N}(f|\nu)\[\prod_{j=1}^{n(\p_r)}{\mbox e}^{-f(X^{*}_j)}\]
t(X_{r+1})
={\mathcal {L}}_{N}(f|\nu)\[\prod_{j=1}^{n(\p_{r+1})}{\mbox e}^{-f(X^{*}_j)}\]
\label{condLaplace}
\EndEq
as desired.
\EndProof
\Remark
The proof of Proposition 2.2 below follows closely an unpublished
proof by Albert Y. Lo for
the case of gamma processes. That is, it is an alternate proof for Lemma
2 in Lo (1984) which yields the appropriate partition representation
for integrals with respect to a Blackwell-MacQueen urn distribution
derived from a Dirichlet process. 
The style of proof exploits properties of partitions
similar to those stated in Pitman (1995a, Proposition 10). In particular
see~\mref{pitmanpart} below. See also the
proof of Lemma 5 in Hansen and Pitman (2001) for general species
sampling models. Details in the proof of Proposition 2.2 translate
into generalizations of a weighted Chinese restaurant algorithm given
in the next section. Proposition 2.1 and 2.2 combine to yield Lemma 2.2.
\EndRemark

\begin{prop}
For $i=1,\ldots,n$, let $g_{i}$ be non-negative functions in
$BM(\mathcal{X})$ then,  
\Eq
\int_{\mathcal {M}}\[\prod_{i=1}^{n}\int_{{\mathcal
    {X}}}g_{i}(x_i)N(dx_i)\]{\mathcal{P}}(dN|\nu)
=\sum_{{\bf
    p}}\prod_{j=1}^{n(\bf p)}
\int_{\mathcal{X}}\[\prod_{i\in
  C_j}g_{i}(x^{*}_{j})\]\nu(dx^{*}_{j}).
\label{momentFP}
\EndEq
Equivalently,
$ 
M(d{\bf X})=\prod_{j=1}^{n(\bf p)}\nu(dX^{*}_{j}).
$
\end{prop}
\Proof
The proof of~\mref{momentFP} proceeds by induction. Case $n=1$ is obvious, 
Now suppose it is true for $n=r$. 
Let $\p_{r+1}$ denote a partition of $\{1,\ldots, r+1\}$ , and define
for each $r>0$, 
\Eq
\phi_{g}(\p_r)=
\prod_{j=1}^{n(\p_r)}
\int_{{\mathcal{X}}}\[\prod_{i\in
  C_{j,r}}g_{i}(x^{*}_{j})\]\nu(dx^{*}_{j})
\label{definefip}
\EndEq
It follows that $\phi_{g}(\p_{r+1})$
is
$$ 
\phi_{g}(\p_r)\int_{{\mathcal Y}}g_{r+1}(v)\nu(dv)
$$
if $n(\p_{r+1})=n(\p_{r})+1$, otherwise  if the index {r+1} is in an
existing cell/table $C_{i,r}$ then it is equivalent to
$$
\phi_{g}(\p_r)\int_{{\mathcal
    Y}}g_{r+1}(v)\pi_{g}(dv|C_{i,r})
$$
where 
$$
\pi_{g}(dv|C_{i,r})=\frac{\[\prod_{l \in C_{i,r}}g_l(v)\]\nu(dv)}
{\int_{{\mathcal Y}}\[\prod_{l \in C_{i,r}}g_l(v)\]\nu(dv)}
$$
for $i=1,\ldots,n(\p_r)$.
Note that this implies that,
\Eq
\sum_{{\bf
    p_{r+1}}}\phi_{g}(\p_{r+1})=
\sum_{{\bf p}_r}\phi_{g}(\p_r)\[\int_{{\mathcal X}}g_{r+1}(v)\nu(dv)
+\sum_{i=1}^{n({\bf p}_{r})}\int_{{\mathcal
    X}}g_{r+1}(v)\pi_{g}(dv|C_{i,r})\].
\label{pitmanpart}
\EndEq
Now by (simple algebra) and the induction hypothesis on $r$ it follows
that,
$$
\sum_{{\bf
    p_{r+1}}}\phi_{g}(\p_{r+1})=\int_{{\mathcal{X}}^n}
\[\int_{\mathcal {X}}g_{r+1}(v)\nu(dv)+
\sum_{j=1}^{n({\bf p}_{r})}g_{r+1}(X^{*}_{j})\]
\[\prod_{i=1}^{r}g_{i}(X_i)\]M(d{\bf
  X}_{r}).
$$
Now utilizing the fact that,
$
M(d{\bf X}_{r+1})=
\[\nu(dX_{r+1})+
\sum_{j=1}^{n({\bf p}_{r})}\delta_{X^{*}_{j}}(dX_{r+1})\]
M(d{\bf X}_{r}),
$
concludes the proof.
\EndProof

\Remark
Of course~\mref{momentFP} in its most basic form leads to well-known
results for moments and cumulants of a Poisson random variable. For
instance, setting $g_i$ to be indicators of a bounded set $A$ yields,  
$$
E(N(A)^n)=\sum_{{\bf
    p}}{\nu(A)}^{n(\bf p)}.
$$
Where $N(A)$ is a Poisson random variable with mean measure
$\nu(A)$.
\EndRemark

\subsection{Chinese restaurant like approximation methods}
In this section a new algorithm for approximating complex
integrals and in fact posterior distributions is described. 
This algorithm works by sequentially sampling from a 
partition distribution and  
structurally behaves similar to the Chinese restaurant process seating
algorithm discussed in Aldous~(1985), Pitman~(1996) and Kerov~(1998). In 
particular, the proposed scheme is influenced by the weighted Chinese
restaurant(WCR) procedure developed in Lo, Brunner and Chan (1996)
[see also Brunner, Chan, James and Lo (2001)] for
mixtures of Dirichlet process and weighted gamma process posterior
models. Ishwaran and James (2001a) subsequently generalise the WCR to
include the class of species sampling mixture models. 
The essence of all these algorithms will be revisted in this section. 
  
Notice that the left hand side of~\mref{momentFP}, 
\Eq
\sum_{{\bf
    p}}\prod_{j=1}^{n(\bf p)}
\int_{\mathcal{X}}\[\prod_{i\in
  C_j}g_{i}(x^{*}_{j})\]\nu(dx^{*}_{j}),
\label{basemodel}
\EndEq
is a function of partitions of the form,
$
\sum_{{\p}}t({\bf p})$. This result is analogous to Lo (1984) where he
points out that complex multiple integrals with respect to a Blackwell-MacQueen urn are equivalent to considerably more manageable sums over
partitions. 
However, it is known that the compexity of the number of partitions
behaves like Bell's number as $n$ increases and hence one needs some method to
approximate such quantities. In order to further illustrate a
connection to a Chinese restaurant we introduce the following expressions;
\Eq
\sum_{{\bf
    p}} \theta^{n({\bf p})}\prod_{j=1}^{n({\bf p})}\Gamma(e_j-\alpha)\prod_{j=1}^{n({\bf p})}
\int_{\mathcal{Y}}\[\prod_{i\in
  C_j}K(X_i|Y^{*}_{j})\]\eta(dY^{*}_{j}),
\label{generalgamma}
\EndEq
and
\Eq
\sum_{{\bf
    p}}{PD}({\bf p}|\alpha,\theta)\prod_{j=1}^{n({\bf p})}
\int_{\mathcal{Y}}\[\prod_{i\in
  C_j}K(X_i|Y^{*}_{j})\]H(dY^{*}_{j}).
\label{generalDirichlet}
\EndEq
The expressions above may arise respectively from a generalized gamma
mixture model and a two-parameter Poisson-Dirichlet process mixture
model. The first expresssion~\mref{generalgamma} appears in James
(2001b) and reduces to an expression for the gamma process in Lo and
Weng (1989) and James (2001a). The latter expression appears in Ishwaran and
James (2001a) and extends the analogous result for the the Dirichlet
process in Lo (1984). That is, the latter corresponds to the marginal
likelihood of $X|Y$ when $Y_1,\ldots,Y_n|P$ are iid $P$, the
law of $P$ is ${\mathcal P}_{\alpha, \theta}(dP|H)$.
Allowing for more flexibility in the interpretation of
$g_i$ and $\nu$ these terms can be written as special cases of
~\mref{basemodel}. Consequently, an understanding of the mechanism behind
the WCR algorithm as outlined in Brunner, Chan and Lo (1996) translates
into a general algorithm which is now described.
From~\mref{pitmanpart} set
\Eq
l (r)= \[\int_{{\mathcal X}}g_{r+1}(v)\nu(dx)
+\sum_{i=1}^{n({\bf p}_{r})}\int_{{\mathcal
    X}}g_{r+1}(x)\pi_{g}(dx|C_{i,r})\].
\EndEq
The procedure relies on a method to generate  partitions ${\bf p}$
based on the following  rule, described in terms of customers entering
a restaurant,  

\begin{PCR}
\Enumerate
\item[Step]1: Seat the first customer to a table with probability 
${l(0)}/{l(0)}=1$.
\item[Step]($r+1$): Given $\p_{r}$, customer $r+1$ sits at
table $C_{j,r}$ with probability
$$
\PR(\p_{r+1}|\p _{r})={l (r)}^{-1}
\int_{\mathcal X}
g_{r+1}(x)
\pi_{g}(dx|C_{j,r}),
$$
where $\p_{r+1}=\p_r\cup\{r+1\in C_{i,r}\}$ for $i=1,\ldots,n(\p _r)$.
Otherwise, customer $r+1$ sits at a new table with probability
$$
\PR (\p _{r+1}|\p _{r}) ={l (r)}^{-1}
\int_{\mathcal X} g_{r+1}(x)\nu(dx).
$$
\EndEnumerate

The completion of Step $n$ produces a
${\bf p}=\{C_1,\ldots,C_{n({\bf p})}\}={\bf p}_n$,
where now ${\bf p}$ is drawn from a density $q({\bf p}|{\bf
  g})$ whose form is described in the Lemma 2.3. 
\end{PCR}

\begin{lem} 
The $n$-step seating algorithm results in a partition ${\bf p}$ drawn from a
density/distribution  given 
by $q({\bf p}|{\bf g})$ that satisfies,
$$
{\mathcal {I}}({\bf p}|{\bf g})q({\bf p}|{\bf g})=
{\prod_{j=1}^{n(\p)}}\int_{\mathcal{X}}\[\prod_{i\in
    C_j}g_{i}(X^{*}_j)\]\nu(dX^{*}_j),
$$
where ${\mathcal{I}}(\p|{\bf g})={\prod_{r=1}^{n}l(r-1)}$.
\end{lem}

\Proof
As in the proof of Proposition 2.2, define 
$\phi_{g}(\p_{r})$ from~\mref{definefip}. Now note carefully that, 
$$
\phi_{g}(\p_{r})\pi_{g}(dX_{j}^{*}|C_{j,r})=
\nu(dX_{j}^{*})\prod_{i\in C_{j,r}}g_{i}(X_{j}^{*})\times 
\prod_{l\neq j}\int\prod_{i\in C_{l,r}}g_{i}(u)\nu(du).
$$
Hence it follows that if $\p_{r+1}=\p_r\cup\{r+1\in C_{j,r}\}$, then 
\Eq
\int g_{r+1}(u)\pi_{g}(du|C_{j,r})=
\frac{\phi_{g}(\p_{r+1})}{\phi_{g}(\p_{r})}
\quad{\mbox and}\quad 
\PR(\p_{r+1}|{\bf p} _r)=
\frac{\phi_{g}(\p_{r+1})}{l (r)\phi_{g}(\p_{r})}.
\label{kaplan}
\EndEq
A similar argument for $\p_{r+1}$ forming a new table shows
that~\mref{kaplan} holds in general. 
Now notice that since $\p _{r+1}$ contains all the information in
$\p_r$, the product rule of probability gives
$$
q(\p_{n}|{\bf g})=\PR({\bf p}_{1})\prod_{r=1}^{n-1}
\PR(\p_{r+1}|{\bf p}_r)=\frac{
\phi_{g}(\p_{n})}{{\mathcal {I}}({\bf p}|{\bf g})},
$$
where $\PR (\p_{1})=\phi_{g}(\p_{1})/l (0)=1$. Now setting $\p=\p_{n}$ yields the desired result.
\EndProof
Now to approximate terms such as ~\mref{basemodel}, draw an iid sample,
say ${\bf p}^{(1)},\ldots,{\bf p}^{(B)}$, of size $B$ from $q({\bf
  p}|{\bf g})$ and use 
$B^{-1}\sum_{b=1}^{B}{\mathcal {I}}({\bf p}^{(b)}|{\bf g})$. The
Chinese restaurant algorithm to generate a draw from $PD({\bf
  p}|\alpha,\theta)$ is recovered by setting $l(r)=\theta+r$,
 $\int_{{\mathcal X}}g_{r+1}(x)\nu(dx):=\theta +n({\bf p}_{r})\alpha$,
and $\int_{{\mathcal
    X}}g_{r+1}(x)\pi_{g}(dx|C_{i,r}):=e_{i,r}-\alpha$. In other words
under this specification $q({\bf p}|{\bf g})=PD({\bf p}|\alpha,\theta)$.
\Remark
The algorithm above is an example of a sequential importance sampling
procedure. The efficient Dirichlet process 
algorithm discussed in MachEachern, Clyde and Liu (1999) can be seen
as a special case of the WCR when using a
binomial kernel. The Chinese restaurant process structure however is
not emphasized in that work. There are also analogous MCMC methods which can
now readily be deduced from the descriptions given in Brunner, Chan, and
Lo (1996) or Ishwaran and James (2001a). See Ishwaran and James
(2001b) and Ishwaran, James and Sun (2001) for applications and ideas
for other algorithms.  
\EndRemark
   
\section{Size-biased generalizations of completely random measures}
In this section it is shown how specific applications of Lemma 2.1 and 2.2
yield explicit disintegration results for a class of random measures which
includes completely random measures. One feature of the analysis
reveals that cumulants assume in many respects the role played by the EPPF in
posterior calculus for random probability measures. The present
construction is influenced by section 4 of Pitman, Perman and
Yor (1992). The results will be
applied throughout.

Let $N$ denote a Poisson process on an arbitrary Polish space ${\mathcal X}={\mathcal S}\times{\mathcal
  Y}$ with intensity $\nu(ds,dy)=\rho(ds|y)\eta(dy)$ for $\rho$ a L{\'
  e}vy measure on the Polish space ${\mathcal S}$  
depending on $y$ in a fairly arbitrary way and $\eta$ a
sigma-finite (non-atomic) measure on ${\mathcal Y}$. Denote the law of $N$ as
${\mathcal P}(dN|\rho,\eta)$. As in Pitman, Perman and Yor (1992, section 4) let $h$ denote an
arbitrary  strictly positive function on ${\mathcal S}$. Furthermore
it is assumed that $h,\rho,\eta$ are selected such that for each bounded
  set $B$ in ${\mathcal Y}$,   
\Eq
\int_{B}\int_{\mathcal S}\min{(h(s),1)}\rho(ds|y)\eta(dy)<\infty.
\label{boundfinite}
\EndEq
Now define a random measure $\mu$ on ${\mathcal Y}$ such that it may
be represented in a distributional sense as, 
\Eq
\mu(dy)=\int_{\mathcal S}h(s)N(ds,dy). 
\label{repcomplete}
\EndEq
The law of $\mu$ is denoted as ${\mathcal P}(d\mu|\rho, \eta)$. When
  $\rho$ does not depend on $y$, then write
  $\rho(ds|y):=\rho(ds)$. Similar to Tsilevich, Vershik and Yor~(2001)
  the term homogeneous will sometimes be applied to this special case
  of $\rho$ and  $\mu$.  
In the case that ${\mathcal S}=(0,\infty)$ and $h(s):=s$ then
following Kingman (1967, 1993) $\mu$ is a completely random measure without a deterministic
component. 
\Remark
Related to completely random measures, Ferguson and Klass (1972) discuss constructions for the class of  L{\'e}vy
processes on $(0,\infty)$ without a Gaussian component but allowing
for fixed points of discontinuity.
See also Wolpert and Ickstadt
(1998a,b), Brix (1999) for recent applications of completely random
measures to spatial statistics. 
The condition~\mref{boundfinite} gaurantees that $\mu$
in~\mref{repcomplete} is a {\it boundedly finite} measure in the
language of Daley and Vere-Jones (1988, Definition 6.1.I.). 
Theorem 6.3.VIII of that work discusses the representation of
completely random measures. Kallenberg~(1997, chapter 10) describes conditions
under which Poisson functionals are finite.  Certain aspects of the 
presentation below are of course implicit in Kallenberg (1986).
\EndRemark 

The Laplace
functional for $\mu$ can be represented as follows
\Eq
{\mathcal {L}}_{\mu}(g)=
\exp \(-\int_{\mathcal{Y}}\int_{\mathcal S}(1-{\mbox
  e}^{-g(y)h(s)})\rho(ds|y)\eta(dy)\).
\label{Laplacec}
\EndEq

As in Pitman, Perman and Yor (1992) the notation 
$T:=\mu({\mathcal Y})$ will be used to denote an almost surely finite
total mass.
 
\subsection{Disintegrations and posterior distributions}

Define, the following moments with respect to the measure $\rho(s|y)$
as for each fixed $n$ and $y$, 
$$
\kappa_{n}(\rho|y)=\int_{\mathcal S}{h(s)}^{n}\rho(ds|y)\quad{\mbox
  and }\quad
\kappa_{n}(\rho)=\int_{\mathcal S}{h(s)}^{n}\rho(ds).
$$
Note that, 
$
m_{\mu}(dv|\rho, \eta)=\[\int_{\mathcal S}h(s)\rho(ds|v)\]\eta(dv)
$
denotes the first moment measure of $\mu$.  

Now similar to the case of Poisson processes let 
$
{\mathcal{P}}(d\mu|\rho, \eta){\prod_{i=1}^{n}\mu(dY_i)}
\label{distc}
$
represent a disintegration of the joint product measure of 
$\{Y_1,\ldots,Y_n,\mu\}$; where
further $\{Y_1,\ldots ,Y_n\}$ can be viewed as conditionally independent
given $\mu$.  The techniques in Section 2 will now be applied to identify the  disintegration which describes the posterior
distribution of $\mu$ given $\{Y_1,\ldots ,Y_n\}$, and the marginal
joint measure of  $\{Y_1,\ldots ,Y_n\}$ and its corresponding
disintegration $({\bf Y}^{*},{\bf p})$. A description of the posterior
law of $\mu$ will now be given. The result will then be justified in
Theorem 3.1.

Now for each $n$, let ${\mathcal{P}}(d\mu|\rho, \eta, {\bf Y})$ denote
the conditional law of $\mu$ corresponding to the random measure, 
\Eq
\mu(\cdot)+\sum_{j=1}^{n(\p)}h(J_{j,n})\delta_{Y^{*}_j}(\cdot) 
\label{post}
\EndEq
where $J_{j,n}$ are independent random variables each with (conditional)
distribution depending on $Y^{*}_j$,
\Eq
\PR(J_{j,n}\in ds|\rho, Y^{*}_j)=\frac{{h(s)}^{e_{j,n}}\rho(ds|Y^{*}_j)}
{\int_{\mathcal S}{h(u)}^{e_{j,n}}\rho(du|Y^{*}_j)}
= \frac{{h(s)}^{e_{j,n}}\rho(ds|Y^{*}_j)}{\kappa_{e_{j,n}}(\rho|Y^{*}_{j})}
\label{post}
\EndEq
and chosen independently  of $\mu$ which is ${\mathcal{P}}(d\mu|\rho,
\eta )$.

The corresponding(conditional) moment measure
and Laplace functional for~\mref{post} are given by
\Eq
m_{\mu}(dv|\rho, \eta, {\bf Y})=m_{\mu}(dv|\rho, \eta)+
\sum_{j=1}^{n({\bf p})}E[h(J_{j,n})|Y^{*}_j]\delta_{Y^{*}_j}(dv)
\EndEq
and
\Eq
{\mathcal {L}}_{\mu}(g|\rho, \eta, {\bf Y})={\mathcal
  {L}}_{\mu}(g|\rho, \eta)
\prod_{j=1}^{n(\p)}\int_{0}^{\infty}{\mbox
  e}^{-g(Y^{*}_j)h(s)}\PR(J_{j,n}\in ds|\rho,Y^{*}_j)
\label{condLaplacem}
\EndEq
The joint marginal measure of ${\bf Y}$ is expressible as 
$$
M_{\mu}(d\mathbf{Y}|\rho, \eta)=
m_{\mu}(dY_1|\rho,\eta)\prod_{i=2}^{n}m_{\mu}(dY_i|\rho,\eta, Y_1,\ldots,Y_{i-1}).
$$

\begin{thm}
Suppose that $\mu$ is a random measure defined by~\mref{Laplacec} and
assume that $\kappa_{n}(\rho|y)<\infty$ for each fixed $y$. Let $g$ be a non-negative or integrable function on
${{\mathcal{Y}}^{n}\times{\mathcal M}}$, then
for each $n\geq 1$,
\Eq
\int_{{\mathcal M}\times{\mathcal{Y}}^{n}} g(\mathbf{Y},\mu )
{\prod_{i=1}^{n}\mu(dY_i)}{\mathcal{P}}(d\mu|\rho,\eta)
=\int_{{\mathcal M}\times{\mathcal{Y}}^{n}} g(\mathbf{Y},\mu )
{\mathcal{P}}(d\mu|\rho,\eta,{\bf Y})
M_{\mu}(d\mathbf{Y}|\rho,\eta).
\label{basicdisint}
\EndEq
The expressions in~\mref{basicdisint} are equivalent to, 
\Eq
\sum_{{\bf
    p}}\int_{{\mathcal{Y}}^{n({\bf p})}}\[\int_{{\mathcal M}}g({\bf
  Y}^{*},{\bf p},\mu){\mathcal{P}}(d\mu|\rho, \eta ,{\bf Y})\]
\prod_{j=1}^{n(\bf p)}
\kappa_{e_{j,n}}(\rho|Y^{*}_{j})\eta(dY^{*}_{j})
\label{FubiniM}
\EndEq
and 
$M_{\mu}(d\mathbf{Y}|\rho, \eta):=\prod_{j=1}^{n(\bf p)}
\kappa_{e_{j,n}}(\rho|Y^{*}_{j})\eta(dY^{*}_{j})$
\end{thm}
\Proof 
First by definition, $M_{\mu}$ is completely determined by 
$$
\int_{{\mathcal {Y}}^{n}}
\[\prod_{i=1}^{n}{\gs}_i(Y_i)\]
M_{\mu}(d{\bf Y}|\rho,\eta)=\int_{{\mathcal M}\times {\mathcal {Y}}^{n}}
\[\prod_{i=1}^{n}{\gs}_{i}(Y_i)\]{\prod_{i=1}^{n}\mu(dY_i)}{\mathcal{P}}
(d\mu|\rho,\eta)
$$
for ${\gs}_i$ in $BM({\mathcal {Y}})$.
But this is equivalent to 
$$
\int_{{\mathcal M}\times {\mathcal {Y}}^{n}}
\[\prod_{i=1}^{n}{\gs}_i(Y_i)\]
\[{\prod_{i=1}^{n}\int_{0}^{\infty}h(s_{i})N(ds_i,dY_i)}\]{\mathcal{P}}(dN|\rho,\eta).
$$
A direct application of Lemma 2.2, or Proposition 2.2,  yields all the desired forms of
$M_{\mu}$. This is seen immediately by setting $g_i(x_i)=h(s_{i}){\gs}_{i}(y_i)$ and
replacing $\nu(dx)$ with $\rho(ds|y)\eta(dy)$.  Now it simply remains to show
that the conditional Laplace functional of $\mu$ is~\mref{condLaplacem}. Since the form of the marginal measure $M_{\mu}$ is
established, Lemma 2.2 now shows that the conditional Laplace
functional is obtained by using the fact that $\mu(g):=\int_{\mathcal
  Y}\int_{\mathcal S}h(s)g(y)N(ds,dy)$ and 
replacing $f(x^{*}_j)$ in Proposition 2.1  with
$g(y^{*}_j)h(s_{j})$. Hence, the conditional Laplace functional of
$\mu$ is equivalent to, 
\Eq
\int_{{\mathcal S}^{n({\bf p})}}\[\int_{\mathcal M}{\mbox e}^{-\mu(g)}
{\mathcal P}(dN|\rho, \eta, {\bf s},{\bf Y}^{*})\]
\prod_{j=1}^{n(\bf p)}
\PR({\js}_{j,n}\in
ds|\rho,\eta, Y^{*}_j).
\label{condLaplace}
\EndEq
\EndProof

\subsection{Cumulants and moment representations}
The joint marginal measure $M_{\mu}$
disintegrates into 
$$
M_{\mu}(d{\bf Y}|\rho,\eta)
=\[\prod_{j=1}^{n({\bf p})}\kappa_{e_{j,n}}(\rho|Y^{*}_{j})\]
\prod_{j=1}^{n({\bf p})}\eta(dY^{*}_{j}),
$$
which shows that the possibly sigma-finite measure of ${\bf Y}$,
disintegrates into an appropriate joint measure of $({\bf Y^{*}},{\bf
  p})$.
Such structures will play a fundamental role throughout.
As a simple application the joint structure  can be
used to obtain expressions for the moments of the corresponding random
variable $\mu(B)$, for $B$ finite, as follows, 
\Eq
E[{\mu(B)}^{n}|\rho,\eta]:=\sum_{\bf p}\prod_{j=1}^{n({\bf p})}
{\int_{B}}\kappa_{e_{j,n}}(\rho|Y^{*}_{j})\eta(dY^{*}_{j})
\label{generalcum}
\EndEq
which corresponds to the classical relationship between moments and
cumulants. Additionally for integrable linear functionals $\mu(f_i)$,
one might be interested in calculating the joint moments, 
$$
E\[\prod_{l=1}^{q}{(\mu(f_l))}^{n_l}|\rho,\eta\]=\int_\mm
\[\prod_{l=1}^{q}\prod_{i=1}^{n_l}\int_\yy
f_l(y_{i,l})\mu(dy_{i,l})\]
{\mathcal P}(d\mu|\rho,\eta)
$$
for integers $n_i$ such that without loss of generality
$n=\sum_{l=1}^{q}n_l$.
An application of Theorem 3.1 easily yields,
\Eq
E\[\prod_{l=1}^{q}{(\mu(f_{l}))}^{n_l}\]=
\sum_{\bf p}\prod_{j=1}^{n({\bf p})}
\int_{{\mathcal Y}}\kappa_{e_{j,n}}(\rho|u)
\[\prod_{l=1}^{q}{f_{l}}^{e^{l}_{j,n}}(u)\]\eta(du)
\label{momcalc1}
\EndEq
where $e^{l}_{j,n}$, satisfying $e_{j,n}:=\sum_{l=1}^{q}e^{l}_{j,n}$, denotes the number of indices associated with $f_l$ in
$C_j$.  
Suppose that $E[T^{n}|\rho,\eta]<\infty$, then an
important case of $~\mref{generalcum}$ is,
\Eq
E[T^{n}|\rho,\eta]
:=\sum_{\bf p}\prod_{i=1}^{n({\bf p})}\kappa_{e_{i,n}}(\Omega)
=\sum_{{\bf p}}\prod_{j=1}^{n({\bf
      p})}\int_{\mathcal Y}\kappa_{e_{j,n}}
(\rho|Y^{*}_{j})\eta(dY^{*}_{j}).
\label{rosetta}
\EndEq
The consequences of this representation will play a major role in
Section 5.

\Example1
As an important example consider the generalised gamma Levy measure
$\rho_{\alpha,b}$ for $b>0$. In this case the $J_{j,n}$ are ${\mathcal
  G}{(e_{j,n}-\alpha, b)}$,
\Eq
\prod_{i=1}^{n({\bf p})}\kappa_{e_{i,n}}(\rho_{\alpha,
  b}):=\prod_{j=1}^{n({\bf
    p})}\Gamma(e_{j,n}-\alpha)b^{-(e_{j,n}-\alpha)}:=b^{-(n-n({\bf p})\alpha)}
\prod_{j=1}^{n({\bf p})}\Gamma(e_{j,n}-\alpha),
\label{cumgen}
\EndEq
and hence 
\Eq
E[{\mu(B)}^{n}|\rho_{\alpha,b},\eta]:=
b^{-n}\sum_{\bf p}b^{n({\bf p})\alpha}{\[\eta(B)\]}^{n({\bf p})}
\prod_{j=1}^{n({\bf p})}\Gamma(e_{j,n}-\alpha)
\EndEq
When $\alpha:=0$ and $b=1$ this expression combined
with $~\mref{generalcum}$ corresponds to the gamma process with shape
measure $\eta$. In this case, where $\eta({\mathcal Y})=\theta$ is finite, it follows that 
normalising the expression~\mref{cumgen} by $E[T^{n}|\rho_{0,1},\eta]$
yields the EPPF, ${PD}({\bf p}|\theta)$ of the Dirichlet process ${\mathcal {PD}}_{0,\theta}(dP|H)$. Otherwise in the
un-normalised case one obtains expressions
for the generalised gamma random measure in Brix (1999).
In the weighted version of this model, discussed in
James (2001b), set
$b:=b(Y^{*}_{j})$, which reduces to expressions for the weighted
gamma process when $\alpha=0$ in Lo and Weng(1989). Note that although
dividing by $E[T^{n}|\rho_{\alpha,b},\eta]$ yields a proper distribution for ${\bf p}$ it
is not an EPPF except for the case of the Dirichlet process. 

\EndExample
\Example2
For the Beta process, with parameters $c$, $A_{o}$,the $J_{j,n}|Y^{*}_{j}$ are
${\mathcal B}(e_{j,n}, c(Y^{*}_{j}))$. Hence,  
\Eq
\kappa_{e_{j,n}}(\rho|Y^{*}_{j})
:={\frac{\Gamma(e_{j,n})\Gamma(c(Y^{*}_{j}))}{\Gamma(e_{j,n}+c(Y^{*}_{j}))}},
{\quad}
{\mbox and}{\quad}
E[J_{j,n}|Y^{*}_{j}]:=\frac{\kappa_{1+e_{j,n}}(\rho|Y^{*}_{j})}{\kappa_{e_{j,n}}(\rho|Y^{*}_{j})}:=\frac{e_{j,n}}{e_{j,n}+c(Y^{*}_{j})}
\EndEq
\EndExample

These types of integrable operations identify joint distributional
structures for $({\bf Y^{*}},{\bf p})$ which have product form. This
fact is  summarized in the next result. The result is important
for applications of mixture models where $\bf Y$ again are missing values and
not observables. The result will also play a significant role in Section 5.

\begin{cor}
Let $\prod_{i=1}^{n}\int_{\mathcal Y}{\gs}_i(Y_i)\mu(dY_i)$ be an integrable function
of ${\mathcal P}(d\mu|\rho, \eta)$, then there exists a conditional
  distribution of ${\bf Y}|{\bf p}$ such that the unique values
  $\{Y^{*}_1,\ldots, Y^{*}_{n(\bf p)}\}$ are independent with distribution 
$$
\PR(dY^{*}_{j}|\rho,\eta,{\gs})\propto\[\prod_{i\in
  C_j}{\gs}_{i}(Y^{*}_{j})\]\kappa_{e_{j,n}}(\rho|Y^{*}_{j})\eta(dY^{*}_{j}).
$$
and ${\bf p}$ has a distribution proportional to $\prod_{j=1}^{n({\bf
    p})}{\int_{\mathcal Y}\[\prod_{i\in  C_j}{\gs}_{i}(Y^{*}_{j})\]\kappa_{e_{j,n}}(\rho|Y^{*}_{j})\eta(dY^{*}_{j})}$.
If the integrability condition still holds when the ${(\gs _i)}$ are
equal to one, then
\Eq
\PR(dY^{*}_{j}|\rho,\eta)\propto\kappa_{e_{j,n}}(\rho|Y^{*}_{j})\eta(dY^{*}_{j})
\label{jointpart},
\EndEq
the distribution of $\bf p$, is proportional to $\prod_{j=1}^{n({\bf p})}\kappa_{e_{j,n}}(\Omega)$.
In the homogeneous case the integrability condition holds only if
$\eta$ is a finite measure. When $\eta$ is a finite measure 
it follows also that the unique values $Y^{*}_{1},\ldots,Y^{*}_{n(\bf p)}$ are
iid $\eta(\cdot)/\eta({\mathcal Y}):=H(\cdot)$.
\end{cor}

Hereafter, based on~\mref{jointpart} and Theorem 3.1,  denote a joint
law (conditional on ${\bf p}$) of  $\{(J_{j,n},Y^{*}_{j})\}$ as 
\Eq
\PR(d{\bf J},d{\bf Y}^{*}|\rho,\eta):=\prod_{j=1}^{n({\bf p})}
\PR(dJ_{j,n}|\rho,Y^{*}_j)\PR(dY^{*}_{j}|\rho,\eta).
\label{jointlaw}
\EndEq
It follows that there exists joint laws of $N,{\bf J},{\bf
  Y}^{*}|{\bf p}$ denoted as
\Eq
{\mathcal P}(dN,d{\bf J},d{\bf Y}^{*}|\rho,\eta):=
{\mathcal P}(dN|\rho,\eta,{\bf J},{\bf Y}^{*})\PR(d{\bf J},d{\bf
  Y}^{*}|\rho,\eta):={\mathcal P}(dN,d{\bf J}|\rho,\eta,{\bf Y}^{*})
\PR(d{\bf Y}^{*}|\rho,\eta)
\label{mixpoisson}
\EndEq
and also a joint law of $\mu,{\bf Y}$ given ${\bf p}$ denoted as 
\Eq
{\mathcal P}(d\mu,d{\bf Y}^{*}|\rho,\eta):={\mathcal P}(d\mu|\rho,\eta,
{\bf Y})\prod_{j=1}^{n({\bf p})}\PR(dY^{*}_{j}|\rho,\eta).
\label{mixmu}
\EndEq

\subsection{Updating and moment formulae}
This section presents a series of important exponential based updating
(change of measure) formulae
which will be used throughout.

\begin{prop}(Updating and moment formulae I)
Let $N$ denote a Poisson process with law ${\mathcal P}(dN|\rho, \eta)$.
In addition, let $f$ denote a positive function on ${\mathcal
  S}\times{\mathcal Y}$. Suppose that $w(\mu)$ is a positive integrable
function of $\mu$, such that it is representable as $w(\mu)={\mbox
  e}^{-N(f)}$. Then,

\noindent 
{\em (i)}
$$
w(\mu){\mathcal P}(d\mu|\rho, \eta)={\mathcal P}(d\mu|{\mbox e}^{-f}\rho,
\eta)E[w(\mu)|\rho,\eta]:={\mathcal P}(d\mu|{\mbox e}^{-f}\rho,
\eta){\mathcal L}_{N}(f|\rho,\eta)
\label{keyswitch}.
$$
\noindent
{\em (ii)}
If $\kappa_{n}({\mbox e}^{-f}\rho|y)<\infty$ then for each $n$,
$$
w(\mu)\prod_{i=1}^{n}\mu(dY_i){\mathcal P}(d\mu|\rho, \eta)
={\mathcal P}
(d\mu|{\mbox e}^{-f}\rho,
\eta,{\bf Y}){\mathcal L}_{N}(f|\rho,\eta)M_{\mu}(d{\bf Y}|{\mbox e}^{-f}\rho,\eta)
$$

\noindent
{\em (iii)}
If $\kappa_{n}(\rho|y)<\infty$ holds then, 
$$
w(\mu){\mathcal P}(d\mu|\rho, \eta, {\bf Y})={\mathcal P}
(d\mu|{\mbox e}^{-f}\rho,
\eta,{\bf Y})E[w(\mu)|\rho, \eta, {\bf Y}]$$
and
$$
{\mathcal L}_{N}(f|\rho,\eta)M_{\mu}(d{\bf Y}|{\mbox
  e}^{-f}\rho,\eta):=
E[w(\mu)|\rho,
\eta, {\bf Y}]M_{\mu}(d{\bf Y}|\rho,\eta),$$

where $E[w(\mu)|\rho, \eta, {\bf Y}]:={\mathcal L}_{N}(f|\rho,\eta) 
\prod_{j=1}^{n({\bf
    p})}\[\int_{\mathcal S}{\mbox e}^{-f(s,Y^{*}_j)}\PR({\js}_{j,n}\in
ds|Y^{*}_j,\rho)\]$.
\end{prop}
\Proof
As in Lemma
2.1, it suffices to show that the two sides in statement (i) have the
same Laplace functional. An application of Lemma 2.1, combined with the
fact that $\mu$ is a functional of N implies that 
\Eq 
\int_{\mathcal M}{\mbox e}^{-\mu(g)}w(\mu){\mathcal P}(d\mu|\rho,
\eta)=\int_{\mathcal M}{\mbox e}^{-\mu(g)}{\mbox e}^{-N(f)}
{\mathcal P}(dN|\rho, \eta):=\int_{\mathcal M}{\mbox e}^{-\mu(g)}
{\mathcal P}(d\mu|{\mbox e}^{-f}\rho,
\eta)E[w(\mu)|\rho,\eta]
\EndEq
which yields statement (i). The result in statement (ii) is
obtained by replacing ${\mbox e}^{-\mu(g)}$
in~\mref{condLaplace} with ${\mbox e}^{-\mu(g)}{\mbox e}^{-N(f)}$ 
and applying Lemma 2.1 to the inner integral.
\EndProof
\begin{prop}(Updating and moment formulae II)
Suppose that $f$ in Proposition 3.1 is replaced by
$f_{n}=\sum_{i=1}^{n}{\fs}_i$, for postive integrable functions
${\fs}_i$, this implies that
$w(\mu):=\prod_{i=1}^{n}w_i(\mu)$ where $w_{i}(\mu)={\mbox
  e}^{-N({\fs}_{i})}$. Now with respect to the the joint model 
${\mathcal P}(d\mu|\rho, \eta)\prod_{i=1}^{n}w_i(\mu)\mu(dY_{i}),$
the following additional formulae are given;(all expressions are assumed to be finite)

\noindent
{\em (i)}
$$
{\mathcal L}_{N}(f_{n}|\rho,\eta):={\mathcal L}_{N}(f_{1}|\rho,\eta)\prod_{i=2}^{n}
{\mathcal L}_{N}({\fs}_i|{\mbox
  e}^{-f_{i-1}}\rho,\eta):=E[\prod_{i=1}^{n}w_{i}(\mu)|\rho, \eta], 
$$
where $f_{i}:=\sum_{j=1}^{i}{\fs}_{j}$.

\noindent
{\em (ii)} The expresssion below are equivalent.
$$
E[w_i(\mu)|{\mbox e}^{-f_{i-1}}\rho, \eta, {\bf
  Y}_{i-1})]m_{\mu}(dY_{i}|{\mbox e}^{-f_{i}}\rho,\eta,{\bf Y}_{i-1}),
$$
$$
E[w_i(\mu)|{\mbox e}^{-f_{i-1}}\rho, \eta, {\bf
  Y}_{i})]m_{\mu}(dY_{i}|{\mbox e}^{-f_{i-1}}\rho,\eta,{\bf Y}_{i-1}).
$$

\noindent
{\em (iii)}
The above statements coupled with Proposition 3.1 imply that the
marginal calculation, 
$$
\int_{\mathcal M}\[\prod_{i=1}^{n}w_i(\mu)\mu(dY_{i})\]{\mathcal
  P}(d\mu|\rho, \eta)
$$
is equivalent to the following formulae,
$$
{\mathcal L}_{N}(f_{n}|\rho,\eta)M_{\mu}(d{\bf Y}|{\mbox e}^{-f_{n}}\rho,\eta),
$$
$$
E[w(\mu)|\rho,
\eta, {\bf Y}]M_{\mu}(d{\bf Y}|\rho,\eta),
$$
$$
E[w_{1}(\mu)|\rho, \eta]m_{\mu}(dY_{1}|{\mbox e}^{-f_{1}}\rho,\eta )
\prod_{i=2}^{n}E[w_i(\mu)|{\mbox e}^{-f_{i-1}}\rho, \eta, {\bf Y}_{i-1})]m_{\mu}(dY_{i}|{\mbox e}^{-f_{i}}\rho,\eta,{\bf Y}_{i-1}),
$$
$$
E[w_{1}(\mu)|\rho, \eta, Y_{1}]m_{\mu}(dY_{1}|\rho,\eta )
\prod_{i=2}^{n}
E[w_i(\mu)|{\mbox e}^{-f_{i-1}}\rho, \eta, {\bf
Y}_{i})]m_{\mu}(dY_{i}|{\mbox e}^{-f_{i-1}}\rho,\eta,{\bf Y}_{i-1}). 
$$

\end{prop}

\Remark 

A notable special case of Proposition 3.1 (i) was established in Lo and Weng
(1989)[see also~Lo(1982) and James (2001a)] for the weighted
gamma process. In Lo and Weng (1989, Proposition 3.1), the statement proceeds
as follows 

\begin{lw}
Let ${\mathcal {G
}}_{\eta ,\beta }$ denote the law of a weighted gamma process with shape $\eta$ and
weight $\beta(\cdot)$, then for each positive $g$
$$
\int_{\mathcal{M}}g(\mu ){\mbox e}^{-\mu(f)}\mathcal{G
}_{\eta ,\beta }(d\mu ) \\
=L_{\mathcal{G}_{\eta ,\beta }}(f)
\int_{\mathcal{M}}g(\mu )\mathcal{G}_{\eta ,\b^* }(d\mu ),
$$ 
where $\b^*=\beta/(1+\beta f)$.  
\end{lw}

This result establishes the absolute continuity of weighted
gamma processes and identifies the specific densities. In other words
the result of Lo and Weng (1989, Proposition 3.1) 
includes the  quasi-invariance result for the gamma process
recently established independently in Tsilevich, Vershik and Yor
(2001, Theorem 3.1). The applications considered by Tsilevich,
Vershik and Yor (2001) are vastly different from Lo and Weng (1989)
 and it is not surprising
that this result emerges in another context. The
development of the exponential formulae used here are directly
inspired by Lo and Weng (1989, Proposition 3.1).
\EndRemark

\subsection{A simple proof for the almost sure discreteness of
  size-biased measures}
 In this section a Fubini argument is used to establish the almost sure
discreteness of $\mu$ with law ${\mathcal P}(d\mu|\rho,\eta)$. This
 will include a simple alternative proof for the class
of completely random measures as discussed in Kingman (1993, Chapter
10). Kingman' s result is based on a modification of Blackwell's
(1973) argument for the Dirichlet process. The present technique is
based on the approach of Berk and Savage~(1979) and Lo and Weng~(1989) for the Dirichlet
process and weighted gamma process respectively. The only requirement I will need is that $\mu$
admits a disintegration, has a 1st moment measure, or is absolutely
continuous with respect to the law of another measure $\mu^{*}$ which
has one. The
latter case of course will yield the result for the stable law. 
Measurability issues vanish on Polish spaces.
The idea is to
apply a 1-step disintegration of $\mu(dx){\mathcal P}(d\mu|\rho,\eta).$
For the arguments below it 
suffices to show that the result holds over all bounded sets B, i.e.
sets such that  $\eta(B)<\infty$ so without loss of generality we
can assume that $\eta$ is a finite measure. 

\begin{prop}(Almost sure discreteness of measures)
Suppose that $\mu$ is ${\mathcal P}(d\mu|\rho,\eta)$ such that 
$\mu$ has a 1st moment measure $m_{\mu}(\cdot|\rho,\eta)$. Otherwise
suppose that 
their exists a measure $\mu^{*}$ which admits a 1st moment measure and
satisfies the absolute continuity relationship,  
\Eq
{\mathcal P}(d\mu|\rho, \eta):=g(\mu^{*}){\mathcal
  P}(d\mu^{*}|\rho^{*},\eta),
\label{absolutecon}
\EndEq
for some positive integrable function $g$. Then $\mu$ is almost surely
discrete. That is,
\Eq
\int_{\mathcal M}\mu(\{x:\mu(\{x\})=0\}){\mathcal P}(d\mu|\rho,\eta):=0
\EndEq
\end{prop}

\Proof
Using a disintegration argument,
$$
\int_{\mathcal M}\mu(\{x:\mu\{x\}=0\}){\mathcal P}(d\mu|\rho,\eta):=
\int\[\int_{\mathcal M}I\{x:\mu(\{x\})=0\}{\mathcal
  P}(d\mu|\rho,\eta, x)\]m_{\mu}(dx|\rho,\eta)
$$
From Theorem 3.1  the law of $\mu$  taken
with respect to the inner term(that is given $x$) is the same as 
$$
\mu(\cdot)+h(J)\delta_{x}(\cdot)
$$
where $h(J)$ is strictly postive and $J$ has law $\propto h(s)\rho(ds|x)$ for
almost all $x$. But since $\mu$ is not negative it follows that 
$$
\mu(\{x\})+h(J)\delta_{x}(\{x\})\ge h(J)>0.
$$
for almost all $x$ with respect to $m_{\mu}$. Hence the inner term is
zero which concludes the result. If again $\mu$ does not
admit such a disintegration, apply the result to the random measure
$\mu^{*}$, with for instance L{\'e}vy measure ${\mbox e}^{-gh}\rho$, or
some other operation, and use the absolute continuity of measures. 
\EndProof

\Remark
This method also implies the almost sure discretenes of the measures
in Section 5. Beyond the mild restriction to Polish spaces, I believe
this is the most general result of this type. The absolute continuity
in~\mref{absolutecon} coupled with the existence of such $\mu^{*}$,
for instance via Proposition 3.1, seems to exhaust the possibilities.
\EndRemark

\section{Intensity rate mixture models, L{\'e}vy moving averages and  shot-noise processes}
In this section analysis for the class of mixture of hazards models
as discussed in the introduction, 
otherwise known as L{\'e}vy-Cox moving average models, is given. Very
little is known about the posterior structure of such models with a
notable exception being the case of mixtures of weighted gamma processes
which is discussed in various degrees of generality in Dykstra and
Laud (1981), Lo (1982), Lo and Weng (1989). Wolpert and Ickstadt
(1998a) and James (2001a) consider semiparametric extensions of this model. 
Full partition based posterior analysis ina a
general multiplicative setting is given in Lo and Weng (1989) and
James (2001a). One consequence of the absence of a general analysis of
this model is the unavailabilty of computational procedures which
sample from the updated or posterior based models.   
As mentioned in the introduction such models are currently used in
Spatial statistical applications and  survival analysis. An interesting
class of models which fits into this framework is the generalised
gamma process proposed in Brix (1999). Wolpert and Ickstadt (1998b) propose
models based on arbitrary L{\'e}vy processes. The analysis here includes
these models as well as mixture models based on the general
size-biased random measures described in section 3. 
Consider the random hazard or intensity rate, 
\Eq
\lambda(t|\mu)=\int_{\mathcal{Y}}K(t|v)\mu(dv),
\label{independent}
\EndEq
where $K$ denotes a known $(\tau,\eta)$-integrable kernel on a Polish space
${\mathcal X}\times {\mathcal Y}$ and $\mu$ is modelled as a
random measure with law ${\mathcal P}(d\mu|\rho, \eta)$.
The representation in~\mref{independent} defines a large class of
random measures $\lambda$  which are not independent increment
processes. 

In this section explicit posterior characteristics of $
\mu$ and hence $\lambda$ based on the multiplicative intensity likelihood,
\Eq
L({\bf X}|\mu ){\mathcal{P}}(d\mu|\rho,\eta)=\left[ \prod_{i=1}^{n}\int_{
\mathcal{Y}}K(X_{i}|Y_{i})\mu(dY_{i})\right]{\mbox e}^{-\mu(f_{K})}
{\mathcal{P}}(d\mu|\rho,\eta)
\label{likelihood}
\EndEq
are derived where,
$
f_{K}(v)=\left[ \int_{\mathcal{S}
}Y(s)K(s|v)\tau (ds)\right]. 
$

Here $X_{i}$ are
observations in a (Polish space) region $\mathcal{X}$ and $Y_{i}$ can be
viewed as missing observations on $\mathcal {Y}$. $Y(s)$ is a non-negative predictable
function which for many applications in event history analysis denotes the
number of observed individuals still at risk just before time
$s$. An important point is that the structural form of $L({\bf
  X}|\mu)$ remains the same under right censoring and left
filtering [see Jacod(1975) and Andersen, Borgan, Gill and Keiding (1993)]. If $Y(s)=1$
then the model may correspond to the likelihood of an inhomogeneous
Poisson process with intensity rate $\lambda(t|\mu)$. 
One purpose of the development of Lemma 2.1 is to handle the
exponential term in ~\mref{likelihood}. Thus mimicking the application
of Lo and Weng (1989, Proposition 3.1). Proposition 3.1  and Theorem
3.1 in this paper readily yield the desired results. Moreover the
appearance of the exponential term combined with Proposition 3.1 show
that analysis of this model only requires the weaker condition, 
\Eq
\kappa_{n}({\mbox e}^{-f_{K}}\rho|y)=\int_{0}^{\infty}
{{\mbox e}^{-f_{K}(y)s}}{h(s)}^{n}\rho(ds|y)<\infty.
\label{condhazard}
\EndEq
The condition~\mref{condhazard} will be assumed throughout this
section and now for instance admits analysis for
the stable law $\rho_{\alpha,0}$, via Proposition 3.1.

An application of Proposition 3.1 combined with Corollary 3.1 show
that  the marginal likelihood is, 
\Eq 
\int_{\mathcal {M}}L({\bf X}|\mu ){\mathcal{P}}(d\mu|\rho, \eta)
=
{\mathcal{L}}_{\mu}
(f_{K}|\rho,\eta)\
\sum_{{\bf p}}\prod_{j=1}^{n({\bf p})}\int_{\mathcal{Y}}
\[\prod_{i\in
    C_j}K(X_{i}|Y^{*}_j)\]\kappa_{e_{j,n}}({\mbox e}^{-f_{K}}\rho|Y^{*}_{j})
\eta(dY^{*}_{j}).
\label{marginal1}
\EndEq

\subsection{Posterior characterizations}
Explicit posterior characterizations are
now given which follows immediately from the Theorem 3.1, Corollary 3.1 and
Proposition 3.1. (Note that it is 
assumed throughout that all relevant integrals are finite).

\begin{thm}
The posterior
distribution of ${\bf Y},\m|{\bf X}$ based on the
model~\mref{likelihood} is representable as, 
$$
\pi(d{\bf Y},d\mu|{\bf X})\propto 
{\mathcal{P}}(d\mu|{\mbox e}^{-f_{K}}\rho, \eta, {\bf Y})
\[\prod_{i=1}^{n(\p)}\PR (dY _i^*|{\mbox e}^{-f_{K}}\rho,\eta,{\bf K})\]\,\pi(\p|{\bf X}),
$$
where
$
\pi({\bf p}|{\bf X})\propto
{\prod_{j=1}^{n(\p)}}\int_{\mathcal{Y}}
\[\prod_{i\in
    C_j}K(X_{i}|Y^{*}_j)\]
\kappa_{e_{j,n}}({\mbox e}^{-f_{K}}\rho|Y^{*}_{j})
\eta(dY^{*}_{j})
$
is a (posterior) distribution of ${\bf p}|{\bf X}$.
\end{thm}

Similar to Lo and Weng (1989, Theorem 4.2), Theorem 4.1 implies for
instance that the 
posterior expectation of the intensity,  $\lambda|{\bf Y}$, is 
\Eq
E[\lambda (t|\m )|{\bf Y}]=\int_{\mathcal Y}
K(t|v)\kappa_{1}({\mbox e}^{-f_{K}}\rho|v)\eta(dv)+\sum_{j=1}^{n(\p)}
K(t|Y^{*}_{j}){\frac{\kappa_{1+e_{j,n}}({\mbox e}^{-f_{K}}\rho |Y^{*}_{j})}
{\kappa_{e_{j,n}}({\mbox e}^{-f_{K}}\rho|Y^{*}_{j})}}
\label{h.expectation}
\EndEq
and hence the posterior expection given ${\bf X}$ is ,  
\Eq
E[\lambda (t|\m )|{\bf X}]=\sum_{\p}\(\int_{\mathcal Y}
K(t|v)\kappa_{1}(\rho_{f_K}|v)\eta(dv)+\sum_{j=1}^{n(\p)}
\int_{\mathcal Y}
K(t|v){\frac{\kappa_{1+e_{j,n}}({\mbox e}^{-f_{K}}\rho|v)}{\kappa_{e_{j,n}}(
{\mbox e}^{-f_{K}}\rho
|v)}}
\pi (dv|C_{j})\)\pi({\bf p}|{\bf X}).
\label{h.expectation}
\EndEq

\Example{Generalised gamma process}
A brief description of the results related to the usage of a  generalised gamma process with intensity
$\rho_{\alpha,b}(ds)\eta(dy)$ as in Brix~(1999) are given. Note in particular that
the result holds for the stable case $b=0$. The posterior distribution
of $\mu$ given ${\bf Y}$ is denoted as 
${\mathcal P}(d\mu|{\mbox  e}^{-f_{K}}\rho_{\alpha, b},\eta,{\bf Y}
)$. That is, $J_{j,n}$ given $Y^{*}_{j}$ are independent ${\mathcal
   G}(e_{j,n}-\alpha,b+f_{K}(Y^{*}_{i}))$, and $\mu$ is now a {\it weighted} generalised gamma process with Laplace
function
$$
\exp \left\{-{\frac{1}{\alpha}}\int_{\mathcal{Y}}[(b+f_{K}(v)+g(v))^{\alpha}-
(b+f_{K}(v))^{\alpha}]\eta(dv)\right\}
$$
The joint moment measure of ${\bf Y}$ can be expressed as, 
\Eq
M_{\mu}(d{\bf Y}|{\mbox e}^{-f_{K}}\rho_{\alpha,
  b},\eta)=\(\prod_{j=1}^{n(\p)}
\Gamma(e_{j,n}-\alpha)\)\prod_{j=1}^{n(\bf
  p)}{(b+f_{K}(Y^{*}_{j}))}^{-(e_{j,n}-\alpha)}\eta(dY^{*}_{j})
\label{genmpart}
\EndEq
which generalizes an expression for the weighted gamma process, see Lo
and Weng (1989) and James (2001a). See James (2001b) for more details
related to the generalised gamma model.
\EndExample

\Remark
Note that from a practical point of view the distribution of $J_{j,n}$
based on ${\mathcal P}(d\mu|{\mbox
  e}^{-f_{K}}\rho,\eta,{\bf Y} )$ may not always be easy
to simulate. If however the moment condition in Theorem 3.1 holds then
one can use
an alternative characterization of the posterior based on ${\mathcal
  P}(d\mu|\rho,\eta,{\bf Y})$. In that case one does not
marginalize over the exponential term but instead works with the
measure,
$$
{\mbox e}^{-\mu(f_{K})}{\mathcal
  P}(d\mu|\rho,\eta,{\bf Y})
$$
\EndRemark

\Remark
James (2001a) gives results for semi-parametric weighted Gamma process
mixture models under more complex multiplicative intensity structures. 
That is for cases where the kernel $K$ depends on a Euclidean
parameter $\psi$ and where for instance there may be several
independent Poisson processes. A careful examination of that work, coupled with
with the results given here provides an obvious way to obtain the
corresponding result for the general processes. A notable wrinkle is
that the Laplace functionals will depend on $\psi$. A discussion of this is omitted for brevity.
\EndRemark

\subsection{Simulating the posterior}
Here the algorithm discussed in Section 2.3 
is applied to this setting to demonstrate a possible approach to
approximate posterior quantities. Again MCMC based methods can also be
deduced from the algorithm below.
First set, 
$$
l(r|{\bf K})=\int_{\mathcal Y}
K_{r+1}(X_{r+1}|v)\kappa_{1}({\mbox e}^{-f_{K}}\rho|v)\eta(dv)+\sum_{j=1}^{n(\p_{r})}
\int_{\mathcal Y}
K_{r+1}(X_{r+1}|v)
{\frac{\kappa_{1+e_{j,r}}({\mbox e}^{-f_{K}}\rho|v)}
{\kappa_{e_{j,r}}({\mbox e}^{-f_{K}}\rho|v)}}
\pi (dv|C_{j,r}),
\label{h.expectation}
$$
where in particular , 
$
l(0|{\bf K})=\int_{\mathcal Y}
K_{r+1}(X_{r+1}|v)\kappa_{1}({\mbox e}^{-f_{K}}\rho|v)\eta(dv).
$

One can now use the variant of the WCR described in Section 2.3 with the
seating rule: Given $\p_{r}$, customer $r+1$ sits at
table $C_{j,r}$ with probability
$$
\PR(\p_{r+1}|\p _{r})={l (r|{\bf K})}^{-1}
\int_{\mathcal Y}
K_{r+1}(X_{r+1}|v)
{\frac{\kappa_{1+e_{j,r}}({\mbox e}^{-f_{K}}\rho|v)}
{\kappa_{e_{j,r}}({\mbox e}^{-f_{K}}\rho|v)}}
\pi (dv|C_{j,r})
$$
where $\p_{r+1}=\p_r\cup\{r+1\in C_{i,r}\}$ for $i=1,\ldots,n(\p _r)$.
Otherwise, customer $r+1$ sits at a new table with probability
$$
\PR(\p _{r+1}|\p _{r}) ={l (r|{\bf K})}^{-1}
\int_{\mathcal Y} K(X_{r+1}|v)\kappa_{1}({\mbox e}^{-f_{K}}\rho|v)\eta(dv).
$$ The completion of Step $n$ produces a
${\bf p}=\{C_1,\ldots,C_{n({\bf p})}\}={\bf p}_n$,
where now, from Lemma 2.3,  ${\bf p}$ is drawn from the density $q({\bf p}|{\bf K})$ which
satisfies,
$$
{\mathcal {I}}(\p|{\bf K})q({\bf p}|{\bf K})=
{\prod_{j=1}^{n(\p)}}\int_{\mathcal{Y}}\[\prod_{i\in
    C_j}K(X_{i}|Y^{*}_j)\]
\kappa_{e_{j,n}}({\mbox e}^{-f_{K}}\rho|Y^{*}_{j})\eta(dY^{*}_j),
$$
This fact, together with Theorem 4.1,  implies that for any integrable
function $t({\bf p})$,
\Eq
\sum_ {\p}t(\p)\pi(\p|{\bf
  X})=\frac{\sum_{\p}t(\p){\mathcal{I}}(\p|{\bf K})q(\p|{\bf
    K})}{\sum_{\p}{\mathcal{I}}(\p|{\bf K})q(\p|{\bf K})}.
\label{montecarlo}
\EndEq
The expression~\mref{montecarlo} and Theorem 4.1 now suggest a 
method to approximate the posterior law ${\mathcal{P}}\(d\mu|{\bf X}\)$
\Enumerate
\item
Using the seating algorithm above, draw B iid random partitions
$\p=\{C_1,\ldots,C_{n(\p)}\}$ from $q(\p|{\bf K})$
\item
Use the value of $\p$ to draw $Y_j^*$ independently
from $\pi(dY^*_{j}|C_j)$ for $j=1,\ldots,n(\p)$.
This yields\\ ${\bf Y}^{*}=(Y_1^*,\ldots,Y_{n(\p)}^*)$.
\item
Using the current value of $({\bf Y}^*,\p)$, approximate a draw from
the random measure
\Eq
\mu_{f_{K}}(\cdot)+\sum_{j=1}^{n(\p)}
h(J_{j,n})\delta_{Y^{*}_j}(\cdot) 
\label{post2.1},
\EndEq
which is distributed as ${\mathcal{P}}(d\mu|{\mbox
  e}^{-f_{K}}\rho
,\eta, {\bf Y})$.
\item
To approximate the posterior law of a functional $g(\mu)$, run the 
previous steps $B$ times independently obtaining values ${\mu}^{(b)}$ with
importance weights ${\mathcal I}(\p^{(b)})$, for $b=1\ldots,B$. Approximate the
law, ${\mathcal{P}}\(g(\mu)\in\cdot|{\bf X}\)$, with
\Eq
\frac{\sum_{b=1}^B
  I\{g({\mu}^{(b)})\in\cdot\}\,{\mathcal{I}}(\p^{(b)}|{\bf K})}
{\sum_{b=1}^B{\mathcal{I}}(\p^{(b)}|{\bf K})}.
\label{posteriorapprox}
\EndEq
 If one only needs to approximate moments, 
 or an integration which yields a $t({\bf p})$ in
closed form, for instance the likelihood,  then steps 2 and 3 can be
eliminated  and one can replace~\mref{posteriorapprox} with 
\Eq
\frac{\sum_{b=1}^B t({\bf p}^{(b)})\,{\mathcal{I}}(\p^{(b)}|{\bf K})}
{\sum_{b=1}^B{\mathcal{I}}(\p^{(b)}|{\bf K})}.
\EndEq
\EndEnumerate

\Remark
Note that in~\mref{post2.1} the main  difficulty is to approximate
a draw from $\mu_{f_{K}}$. Brix~(1999) discusses methods on
how to approximate a generalized gamma process ${\mathcal P}(d\mu|\rho_{\alpha,b},\eta)$. It
should be straightforward to extend this to a
${\mathcal P}(d\mu|{\mbox e}^{-f_{K}}\rho_{\alpha,b}, \eta)$.
See also Wolpert
and Ickstadt~(1998b) for some possible ideas in the general setting.
I believe that the mixture representations given in the next section
may also be useful in this regard.
\EndRemark

\section{Analysis of a Scaling operation which arises in Brownian
  Excursion theory}
The previous section describes applications of various 
  exponential change of measure operations. As shown from Lemma 2.1 this
  operation results in a change of measure from a Poisson process
  with intensity ${\nu}$ to  another Poisson
  process law with intensity ${\mbox e}^{-f}\nu$. An important aspect
  of that is one can still apply directly Lemma 2.2 to the
  transformed Poisson law, which yields the
  various results in the previous section. In particular this
  operation transforms processes $\mu$ which do not admit moment
  measures to ones which do. An important example is the stable law which is
  transformed to a form of weighted generalised gamma process.
Another important operation, besides the exponential change of
measure,  is a
  type of scaling, which arises for instance in Brownian excursion
  theory[see for instance Pitman and Yor (1992, 1997, 2001)]. This operation
  no longer preserves the Poisson nature of $N$ or similarly the
  structure of the biased models $\mu$. This in itself does not present a
  major obstacle as one could still apply Lemma 2.2 to the Poisson law first.  
The law resulting from the scaling operation may not have an
  obviously understandable form. However, indeed hidden in a scaling
  operation is an exponential form via a gamma integral
  identity. This identity has been used frequently in various
  contexts. Here, a slightly different variation will be used,
  where the exponential change of measure idea  will be applied
  internally leading to a variety of interesting consequences.
As an important special case
  we look at the ${PD}(\alpha, \theta)$
  model. In general, analysis of the simple scaling structure leads to results quite related to Pitman and Yor
  (1992, 1997, 2001). [See also Perman, Pitman and Yor (1992), Section 4].
In particular see Pitman and Yor (1992) Section
  3. 
\Remark
It will become quite clear to the experts,on excursion theory and such
matters, that the overlap with Pitman
and Yor (1992, Section 3) is hardly coincidental. Although this was
not my initial motivation. 
Some of the results
given below amplify on Pitman and Yor (1992, Theorem 3.1 and
especially Remarks 3.3 and 3.4). This
section may be viewed as extensions of their Section 3. Given the
new results I obtain for the ${PD}(\alpha,\theta)$, among
other things, the
present exposition should clearly
provide new insights, for the experts, into matters of which I myself
have no expertise. Again applications of Lemma 2.1 and Proposition 3.1
play a fundamental role.
\EndRemark

\Remark 
In this section it is assumed that $\eta$ satisifies the integrability
condition in Corollary 3.1 when ${\gs}:=1$. That is when $\rho$ is
homogeoneous then $\eta:=cH$ for some scalar $c$. For the Dirichlet
process $\eta:=\theta H$ but is not to be confused with $\theta$ used
below to mimic the scaling operation associated with $PD(\alpha,
\theta)$ for $\alpha>0$.
\EndRemark
 
The present analysis yields new mixture representations of random
measures $N$ and $\mu$. What is most important is how they arise
within the context of the scaling operation. [The reader should again note Pitman and Yor
(1992, Remark 3.4)]. This leads to an analogous representations of
random probability measures defined as,
\Eq
P(\cdot) ={\frac{\mu(\cdot)}{\mu({\mathcal{Y}})}}=
{\frac{\int_{\mathcal S}h(s)N(ds,\cdot)}
{\int_{\mathcal S}h(u)N(du,{\mathcal Y})}}:={\frac{\mu(\cdot)}{T}}
\label{randommeasure}
\EndEq
whose law is determined by an appropriate law on $N$. That is either a
Poison law or a {\it scaled} Poisson law to be described below.
In addition
various characterizations of the the posterior distributions
of $N$,$\mu$, and $P$ are obtained. 
These results are applied to obtain general identities and
representations  for
the two-parameter family with parameters $0<\alpha< 1$,
$\theta>-\alpha$. This  extends  an identity given for the case
$PD(\alpha,\alpha)$ in Pitman and Yor (2001). 
 
It is known that the general ${\mathcal PD}_{\alpha,\theta}(dP|H)$
cannot directly be defined via  normalization of an independent
  increment process. The exceptions are the Dirichlet process and the
 Stable law process. Pitman and Yor (1997),[see also Tsilevich,
Vershik and Yor (2000)] establish the following relationship. 
 Let $PD_{\alpha,
  \theta}(d{\mu}|\eta)$ denote a law of $\mu$  such that its
normalisation results in a Poisson-Dirichlet random probability
measure,with law denoted as ${\mathcal
  {PD}}_{\alpha,\theta}(P|H)$. In particular ${\mathcal {PD}}_{\alpha,
  0}(d\mu|\eta):={\mathcal P}(d\mu|\rho_{\alpha},\eta)
$ is the Stable process. 
From Pitman and Yor (1997), [see
in particular 
Tsilevich, Vershik and Yor (2000) where this form is taken], it
follows that for $0<\alpha< 1,\theta>-\alpha$, 
\Eq
{\mathcal PD}_{\alpha,\theta}(d{\mu}|\eta)=c_{\alpha,\theta}{T}^{-\theta}{\mathcal PD}_{\alpha,0}(d\mu|\eta),   
\label{pychange}
\EndEq
where $c_{\alpha,\theta}$ is a normalizing constant.
Pitman and Yor
(1997) describe various results concerning the ${PD}(\alpha,
  \theta)$ law of the $(P_i)$,  related to this fact. 
Notably, 
Proposition 21, Proposition 22, and Proposition 33. In particular they
show that if the sequence $(P_i)$ is ${PD}(\alpha, \theta)$ then 
\Eq
PD(\alpha, \theta)=\int_{0}^{\infty}PD(\alpha|t)\gamma_{\alpha,
  \theta}(dt)
\EndEq
where $PD(\alpha|t)$ denotes the conditional Poisson-Kingman[see
section 8] law of
$PD(\alpha)$ conditioned on $T$, $\gamma_{\alpha,\theta}(dt)=c^{-1}_{\alpha,
  \theta}t^{-\theta}f_{\alpha}(t)$, 
$c_{\alpha,
  \theta}=E[T^{-\theta}]=\int_{0}^{\infty}t^{-\theta}f_{\alpha}(t)$
is the normalising constant and $f_{\alpha}(t) $ denotes the density of a
stable law random variable which is the distribution of $T$. 
Tsilevich, Verhsik and Yor (2000)
using~\mref{pychange} in combination with the following gamma identity for $\theta>0$,
\Eq
T^{-\theta}:={\frac{1}{\Gamma(\theta)}}\int_{0}^{\infty}v^{\theta-1}{\mbox e}^{-vT}dv,
\label{gammaid}
\EndEq
establish quite remarkably and simply a two-parameter extension of the
Markov-Krein correspondence via the Laplace functional of a stable
law.[Their result will be extended in a general fashion in Section 6].
This identity is also used in Perman, Pitman and Yor~(1992) and
Pitman and Yor (1997) among other places. The results discussed above
are used primarily to deduce properties related to the
normalized process $\mu(\cdot)/T$. That is, equivalently the
$PD(\alpha,\theta)$ family. The interest here however is in another characterization of
the law  ${\mathcal P}_{\alpha,\theta}(d\mu|\eta)$  which allows more direct
usage of it and of course implies results for the normalized process
and synonymously $PD(\alpha,\theta)$. The method relies on
using the identity~\mref{gammaid} conditioning on various transformed
densities for $V$ rather than $T$ which will lead to a variety of interesting results. 
This analysis will be applied to general processes,
$\mu$ and $N$, subject to the same type
of scaling operation.  For all $\theta>-\alpha$, in particular the
case $-\alpha<\theta\leq 0$, the same identity~\mref{gammaid} will be used 
for each fixed n, 
\Eq
{\frac{1}{\Gamma(n)}}\int_{0}^{\infty}T^{n}v^{n-1}{\mbox e}^{-vT}dv:=1,
\label{gammaid2}
\EndEq
where the main point now is to work with $T^{n}$ rather than its
reciprocal, and otherwise use the fact that the left hand-side of
~\mref{gammaid2} is one. 
Now assuming that 
\Eq
E[T^{-\theta}|\rho,\eta]:=\int_{\mathcal M}T^{-\theta}{\mathcal
  P}(dN|\rho,\eta):=\int_{\mathcal M}T^{-\theta}{\mathcal
  P}(d\mu|\rho, \eta)<\infty
\label{negmomcond1}
\EndEq
consider the equality
of laws,
\Eq
{\ps}(dN|\rho, \eta, \theta):=
{\frac{T^{-\theta}{\mathcal P}(dN|\rho, \eta)}{E[T^{-\theta}|\rho,\eta]}},
{\mbox and }, 
{\ps}(d\mu|\rho,\eta,\theta):=
{\frac{{T}^{-\theta}{\mathcal P}(d\mu|\rho, \eta)}{E[T^{-\theta}|\rho,\eta]}}
\label{scalelawp}
\EndEq 
For each fixed $v$, let
$
E[{\mbox e}^{-vT}|\rho, \eta]:=\int_{\mathcal M}{\mbox e}^{-vT}{\mathcal
  P}(dN|\rho, \eta)
$
denote the Laplace transform of $T$ taken relative to ${\mathcal
  P}(dN|\rho, \eta)$.
 The following relationship will prove
useful and allows analysis for  the case $-\alpha<\theta\leq 0$. Suppose that for $n\ge 1$, 
\Eq
E[T^{n}|\rho,\eta, \theta+n]:=\int_{\mathcal M}T^{n}{\ps}(dN|\rho, \eta,
\theta+n)={\frac{E[T^{-\theta}|\rho,\eta]}{E[T^{-(\theta+n)}|\rho,\eta]}}<\infty
\label{mix2cond}
\EndEq
then 
\Eq
{\ps}(dN|\rho, \eta,
\theta):={\frac{T^{n}{\ps}(dN|\rho, \eta, \theta+n)}{E[T^{n}|\rho,\eta, \theta+n]}},
{\mbox and }\quad 
{\ps}(d\mu|\rho, \eta,
\theta):={\frac{T^{n}{\ps}(d\mu|\rho, \eta, \theta+n)}{E[T^{n}|\rho,\eta, \theta+n]}}
\label{doublemix2}.
\EndEq
Additionally, denote the laws of $P$ taken relative to ${\mathcal
  P}(\cdot|\rho,\eta)$ and ${\ps}(dN|\rho,\eta,\theta)$ as ${\mathcal
  P}(dP|\rho,\eta)$ and ${\ps}(dP|\rho,\eta,\theta)$ respectively. When
$\rho$ does not depend on $y$ and $\eta:=cH$, then 
\Eq
P:=\sum_{i=1}^{\infty}P_{i}\delta_{Z_{i}} 
\label{repP}
\EndEq
where the sequence $(P_i)$ is independent of $(Z_i)$ which are iid $H$. 
In that case the analysis of $P$ is in principle equivalent to the
analysis of $(P_{i})$. That is the $P$ models are special cases of
species sampling models~[See Pitman (1996), Hansen and Pitman (2001)
and Ishwaran and James (2001a)].
Hence, for instance, the results of Perman,
Pitman and Yor (1992) and Pitman (1995b), concerning the EPPF etc.; can
be applied to this setting to obtain information about $P$. When
$\rho$ depends on $y$, then one can still represent $P$ in the form
$~\mref{repP}$.  However the
independence property between $(P_{i})$ and $(Z_i)$ no longer
holds. It will become clear that if interest is simply the marginal
distributional properties of the $(P_{i})$ then the results of Pitman,
Perman and Yor (1992) and Pitman (1995b) can be applied using $\Omega$
rather than $\rho$. At any rate, a different analysis will be used
here which is based primarily on information contained in the random
measures $N$ and $\mu$ which will yield relevant information about
the $(P_i)$ etc. This is useful even in the species sampling case
where the EPPF may be intractable or does not easily convey
information about $P$. 

\Remark 
Although the emphasis seems to be on the scaled laws
${\ps}(\cdot|\rho,\eta,\theta+n)$ this is only partially the case. In
particular it should be clear that the negative moment
conditions ~\mref{negmomcond1} and~\mref{mix2cond} do not hold in
general. This is the case for a  gamma process, which could be denoted
as ${\mathcal P}_{0,\theta}(d\mu|\eta)$, and hence it is not a proper
${\ps}(d\mu|\rho,\eta,\theta)$. The forthcoming discussion
investigates properties of general $N$ and
$\mu$ based on  manipulation of both \mref{gammaid} and ~\mref{gammaid2}.
The results below related to 
${\ps}(\cdot|\rho,\eta,\theta+n)$ hold provided ~\mref{negmomcond1}
and~\mref{mix2cond} are true.
\EndRemark

\subsection{Mixture representations for general processes}

The results follow from
the proof of  Theorem 5.1 below.
\begin{prop}
The follwing
identities hold which have various implications.

\noindent
{\em (i)}For $n>0$,
\Eq
{\Gamma}(n):=\sum_{{\bf p}}\int_{0}^{\infty}
\[\prod_{j=1}^{n({\bf p})}\kappa_{e_{j,n}}({\mbox e}^{-vh}\Omega)\]v^{n-1}E[{\mbox e}^{-vT}|\rho, \eta]dv,
\EndEq
which is equivalent to,
$
\int_{0}^{\infty}E[T^{n}|{\mbox e}^{-vh}\Omega,\eta]{v}^{n-1}E[{\mbox
  e}^{-vT}|\rho, \eta]dv.
$

\noindent
{\em (ii)}
Statement (i) implies that there exist a joint distribution of ${V, {\bf
  p}}$ given by,
\Eq
\pi_{n,0}(dv,{\bf p}|\rho,\eta):=
{\frac{\[\prod_{j=1}^{n({\bf p})}\kappa_{e_{j,n}}({\mbox
      e}^{-vh}\Omega)\]
v^{n-1}E[{\mbox e}^{-vT}|\rho, \eta]dv}{\Gamma(n)}}.
\label{jointp1}
\EndEq
\noindent
{\em (iii)}
Suppose that~\mref{mix2cond} holds for $\theta+n>0$ then, 
$
\Gamma(\theta+n)E[T^{-(\theta+n)}|\rho,\eta]:=
\int_{0}^{\infty}v^{\theta+n-1}E[{\mbox
  e}^{-vT}|\rho, \eta]dv$,
which identifies a random variable  $V$ on $(0,\infty)$ with density,
\Eq
\pi_{\theta+n}(dv|\rho,\eta):=
{\frac{v^{\theta+n-1}E[{\mbox
  e}^{-vT}|\rho, \eta]dv}{\Gamma(\theta+n)E[T^{-(\theta+n)}|\rho,\eta]}}
\EndEq
Additionally, 
$E[T^{n}|\rho,\eta,\theta+n]:=\int_{0}^{\infty}E[T^{n}|{\mbox
  e}^{-vh}\Omega,\eta]\pi_{\theta+n}(dv).
$
Hence there exists a random variable $V$ such that $V,{\bf p}$ has joint density,
\Eq
\pi_{\theta+n,\theta}(dv,{\bf p}|\rho,\eta):={\frac{E[T^{-(\theta+n)}|\rho,\eta]}{E[T^{-\theta}|\rho,\eta]}}\[\prod_{j=1}^{n({\bf p})}\kappa_{e_{j,n}}({\mbox
      e}^{-vh}\Omega)\]\pi_{\theta+n}(dv)
\label{jointp2}
\EndEq
\noindent
{\em (iv)} 
The marginal distribution of ${\bf p}$ in~\mref{jointp1}
and~\mref{jointp2} is given by,
\Eq
\pi_{\theta+n,\theta}({\bf p}|\rho,\eta):={\frac{E[T^{-(\theta+n)}|\rho,\eta]}{E[T^{-\theta}|\rho,\eta]}}\int_{0}^{\infty}\[\prod_{j=1}^{n({\bf p})}\kappa_{e_{j,n}}({\mbox
      e}^{-vh}\Omega)\]\pi_{\theta+n}(dv)
\label{EPPFmix}
\EndEq
\end{prop}

\Remark 
The formula $\pi_{n,0}({\bf p})$ derived from~\mref{EPPFmix}
corresponds to the EPPF formula given in Pitman (1995b, Corollary 6
formula (32)), when $h(s):=s\in(0,\infty)$ and $\rho$ is homogeneous.  
Pitman~(1995b) uses the identity~\mref{gammaid} applied to
$T^{-n}$. Indeed it will be stated formally that~\mref{EPPFmix}
are EPPF formulas which can be seen as an extension of Pitman's result. 
\EndRemark

The next result identifies ${\ps}(dN|\rho,\eta, \theta)$ and in fact
arbitrary Poisson laws ${\mathcal P}(dN|\rho, \eta)$ as a mixture
relative to Poisson random measures $N^{*}$ with conditional(given V) laws
${\mathcal P}(dN|{\mbox e}^{-vh}\rho,\eta)$. Moreover since the conditional law on $N^{*}$ has
moments  the result also contains mixture representations of $\mu$ of the type
$N^{*}+\sum_{j=1}^{n({\bf p)}}\delta_{J_{j,n},Y^{*}_{j}}$. The
marginal distribution of $N^{*}$ actually depends on $n$.

\begin{thm}(Poisson Mixture representations)
\noindent
{\em (i)}
For $\theta>0$, the law ${\ps}(dN|\rho, \theta, \eta)$ defined in~\mref{scalelawp} is
expressible as the following mixture,
\Eq
{\ps}(dN|\rho, \theta, \eta):=\int_{0}^{\infty}{\mathcal P}(dN|{\mbox
  e}^{-vh}\rho, \eta)\pi_{\theta}(dv|\rho,\eta).
\EndEq
This indicates that $N|V=v$ has a Poisson law ${\mathcal P}(dN|{\mbox
  e}^{-vh}\rho, \eta)$  with intensity ${\mbox
e}^{-vh(s)}\rho(ds|y)\eta(dy)$ and $V$ is a random variable on
$(0,\infty)$ with density $\pi_{\theta}(dv|\rho,\eta)$. Equivalently
the Laplace function of $N$ is,
\Eq
{\mathcal L}_{N}(f|\rho,\eta,\theta):=\int_{0}^{\infty}
{\mathcal L}_{N}(f|{\mbox e}^{-vh}\rho,\eta) \pi_{\theta}(dv|\rho,\eta)
\EndEq
[Compare statement (i) with Pitman and Yor (1992, Remarks 3.3, 3.4)]

\noindent
{\em (ii)} 
For each $n>0$, a Poisson law ${\mathcal P}(dN|\rho,\eta)$ can be represented
as, 
\Eq
{\mathcal P}(dN|\rho,\eta):=\sum_{\bf p}\int_{{\mathcal
  S}^{n({\bf p})}\times{\mathcal Y}^{n({\bf p})}\times\rr^{+}}
{\mathcal P}(dN,d{\bf J},d{\bf Y}^{*}|{\mbox
  e}^{-vh}\rho,\eta)\pi_{n,0}(dv,{\bf p}|\rho,\eta)
\EndEq
That is from,~\mref{mixpoisson},the random measure $N$ can be
represented as
\Eq
N^{*}+\sum_{j=1}^{n({\bf p})}\delta_{(J_{j,n},Y^{*}_{j})},
\label{psumrep}
\EndEq 
where given
$V,{\bf p}$ with law $\pi_{n,0}(dv,{\bf p}|\rho,\eta)$, $N^{*}$ is
conditionally independent of ${\bf J},{\bf Y}^{*}$, 
with distribution ${\mathcal P}(dN|{\mbox e}^{-vh}\rho, \eta,\rho)$. The random variables $({\bf J},{\bf Y}^{*})$
given  $V,{\bf p}$ have distribution $\PR(d{\bf J},d{\bf Y}^{*}|{\mbox
  e}^{-vh}\rho,\eta)$ as in~\mref{jointlaw}. 

\noindent
{\em (iii)}
For $\theta+n>0$,this includes $-\alpha<\theta\leq 0$, the law 
${\ps }(dN|\rho,\eta,\theta)$ is equivalent to, 
\Eq
\sum_{\bf p}\int_{{\mathcal
  S}^{n({\bf p})}\times{\mathcal Y}^{n({\bf p})}\times\rr^{+}}
{\mathcal P}(dN,d{\bf J},d{\bf Y}^{*}|{\mbox
  e}^{-vh}\rho,\eta)\pi_{\theta+n,\theta}(dv,{\bf p}|\rho,\eta).
\EndEq
That is,  $N$ can be represented as~\mref{psumrep} except that the distribution of
$V,{\bf p}$ is $\pi_{\theta+n,\theta}(dv,{\bf p}|\rho,\eta)$.
When $\theta:=0$ statement (iii) reduces to statement (ii).
\end{thm}

\Proof
For the result in (i), it suffices to evaluate the Laplace functional of ${\ps}$ defined
in~\mref{scalelawp}. That is,
$$
E[T^{-\theta}|\rho,\eta]\int_{\mathcal M}{\mbox
  e}^{-N(f)}{\ps}(dN|\rho, \eta, \theta)
:=
\int_{\mathcal M}{\mbox
  e}^{-N(f)}T^{-\theta}{\mathcal P}(dN|\rho, \eta).
$$
Now apply the identity~\mref{gammaid} to see that the right hand side
is equal to
\Eq
{\frac{1}{\Gamma(\theta)}}\int_{0}^{\infty}v^{\theta-1}E[{\mbox
  e}^{-N(f)}{\mbox e}^{-vT}|\rho, \eta]dv.
\EndEq
An application of Lemma 2.1 and dividing through by $E[T^{-\theta}|\rho,\eta]$
now show that the Laplace functional of ${\ps}(dN|\rho, \eta, \theta)$ is, 
\Eq
\int_{0}^{\infty}E[{\mbox
  e}^{-N(f)}|{\mbox e}^{-vh}\rho,
\eta]{\frac{v^{\theta-1}E[{\mbox e}^{-vT}|\rho,
    \eta]dv}{E[T^{-\theta}|\rho,\eta]\Gamma(\theta)}},
\EndEq
as desired.
For (ii), use~\mref{gammaid2} as follows
\Eq
\Gamma(n){\mathcal L}_{N}(f|\rho,\eta):=\int_{0}^{\infty}
E[T^{n}{\mbox
  e}^{-N(f)}{\mbox e}^{-vT}|\rho, \eta]v^{n-1}dv.
\EndEq
Noting that $T^{n}:=\int_{{\mathcal Y}^{n}}\prod_{i=1}^{n}\mu(dY_{i})$, the result
follows immediately by an application of (ii) in Proposition 3.1 and Corollary
3.1. Statement (iii) follows by first applying (i) to
${\ps}(dN|\rho,\eta,\theta+n)$ in~\mref{doublemix2} and then applying
(ii) to ${\mathcal P}(dN|{\mbox e}^{-vh}\rho,\eta)$. 
\EndProof

\Remark
A clear distinction between a Poisson law ${\mathcal P}(dN|\rho,\eta)$ and the
law on $N$ defined as ${\ps}(dN|\rho,\eta,\theta)$ was used
above. This distinction was made because the conditions and techniques used to derive the results were slightly
different. Hereafter, for some brevity, the notation $\theta$ will be
used for all objects. When $\theta$ is set equal to zero (when
applicable) this will correspond to results for the Poisson based laws
of $N$, $\mu$, and $P$ and their corresponding mixing laws.
\EndRemark

\begin{cor}(Mixture representations for $\mu$)
The results (i),(ii), (iii) in Theorem 5.1 imply analogous results for
${\ps}(d\mu|\rho,\eta,\theta)$ and ${\mathcal P}(d\mu|\rho, \eta)$. 

\noindent
{\em (i)}For $\theta>0$, 
$$
{\ps}(d\mu|\rho,\eta,\theta):=\int_{0}^{\infty}{\mathcal P}(d\mu|{\mbox
  e}^{-vh}\rho, \eta)\pi_{\theta}(dv|\rho,\eta)
$$

\noindent
{\em (ii)} For $n\ge 1$ and  $\theta+n>0$, including $-\alpha<\theta\leq 0$,
Statement (ii) and (iii) in Theorem 5.1 implies that,
for each $n>0$ the random measures $\mu$ with distribution ${\mathcal
  P}(d\mu|\rho,\eta)$ or ${\ps}(d\mu|\rho,\eta,\theta)$ can be represented as, 
\Eq
\mu^{*}+\sum_{j=1}^{n({\bf p})}h(J_{j,n})\delta_{Y^{*}_{j}},
\label{musumrep}
\EndEq  
where  $\mu^{*}|V,{\bf p}$ is ${\mathcal P}(d\mu|{\mbox
  e}^{-vh}\rho)$, and $(J_{j,n},Y^{*}_{j})|V,p$ are conditionally
independent of $\mu^{*}$ with distribution $\PR(d{\bf J},d{\bf Y}^{*}|{\mbox
  e}^{-vh}\rho,\eta)$. The distribution of $V,{\bf p}$ is
  $\pi_{\theta+n,\theta}(dv,{\bf p}|\rho,\eta)$.

\noindent
{\em (iii)}
Equivalently these
results yield the following expressions for the respective Laplace
functionals. 
For $\theta>0$, 
\Eq
{\mathcal L}_{\mu}(g|\rho,\eta,
\theta):=\int_{0}^{\infty}{\mathcal L}_{\mu}(g|{\mbox
  e}^{-vh}\rho)\pi_{\theta}(dv|\rho,\eta), 
\EndEq
and for $\theta+n>0$,$n\ge 1$, 
\Eq
{\mathcal L}_{\mu}(g|\rho,\eta,
\theta):=\sum_{{\bf p}}\int_{0}^{\infty}
{\mathcal L}_{\mu}(g|{\mbox
  e}^{-vh}\rho,\eta,{\bf p})
\pi_{\theta+n,\theta}(dv,{\bf
  p}|\rho,\eta),
\label{laplacemix}
\EndEq
where 
$$
{\mathcal L}_{\mu}(g|{\mbox
  e}^{-vh}\rho,\eta,{\bf p}):=
\int_{{\mathcal Y}^{n({\bf p})}}{\mathcal L}_{\mu}(g|{\mbox
  e}^{-vh}\rho,\eta,{\bf Y})
\prod_{j=1}^{n({\bf p})}\PR(dY^{*}_{j}|{\mbox
  e}^{-vh}\rho,\eta)
$$
Setting $\theta=0$ in~\mref{laplacemix} yields an identity for the
Laplace functional of ${\mathcal P}(d\mu|\rho,\eta)$. 
\end{cor}

Now mixture representations for $P$ are given.
First set, 
$$
p^{*}_{n}:={\frac{\mu^{*}({\mathcal Y})}{\mu^{*}({\mathcal Y})+
\sum_{j=1}^{n({\bf p})}h(J_{j,n})}}:={\frac{T^{*}}{T^{*}+\sum_{j=1}^{n({\bf p})}h(J_{j,n})}}
$$
and define a random probability measure as 
$$
P^{*}(\cdot):={\frac{\mu^{*}(\cdot)}{\mu^{*}({\mathcal
    Y})}}:={\frac{\mu^{*}(\cdot)}{T^{*}}}.
$$
\begin{prop}(Mixture representations for random probability measures)

\noindent
{\em (i)}
If $\theta>0$ then the random probability measure,
\Eq
{\ps}(dP|\rho,\eta,\theta):=\int_{0}^{\infty}{\mathcal
  P}(dP|{\mbox e}^{-vh}\rho,\eta)\pi_{\theta}(dv|\rho,\eta)
\EndEq
\noindent
{\em (ii)} If $P$ is ${\ps}(dP|\rho,\eta,\theta)$ or ${\mathcal P}(dP|\rho,\eta)$ then for each
  $n\ge 1$ and $\theta+n>0$, it is equivalent in
  distribution to the random probability measure
\Eq
p^{*}_{n}P^{*}(\cdot)+(1-p^{*}_{n})\frac{\sum_{j=1}^{n({\bf
  p})}h(J_{j,n})\delta_{Y^{*}_{j}}(\cdot)}{\sum_{j=1}^{n({\bf p})}h(J_{j,n})}
\label{rpm}
\EndEq
with distribution specified by statement (ii)
Corollary 5.1.
\end{prop}

\subsection{Duality of mixture representations and posterior
  distributions}
 
In this section $Y_{1},\ldots,Y_{n}|P$ are iid random variables with
distribution $P$. That is, this implies the joint model for ${\bf Y}|P$, 
is $\prod_{i=1}^{n}P(dY_{i})$.
The law of $P$ is either ${\mathcal P}(dP|\rho,\eta)$
or ${\ps }(dP|\rho,\eta,\theta)$. The interest is in obtaining
posterior distributions  of $N$, $\mu$ and hence $P$ and relevant
information about the marginal structure of ${\bf Y}$. 
Define a conditional distribution of
$V|{\bf Y}$ as,
\Eq
\pi_{\theta+n}(dv|{\bf p},{\bf Y^{*}})\propto\[\prod_{j=1}^{n({\bf p})}
\kappa_{e_{j,n}}({\mbox e}^{-vh}\rho|Y^{*}_{j})\]\pi_{\theta+n}(dv)
\label{postvy}
\EndEq
when $\rho$ does not depend on $y$, then the distribution of $V|{\bf
  Y}$ only depends on ${\bf p}$ and~\mref{postvy} reduces to 
\Eq
\pi_{\theta+n}(dv|{\bf p},{\bf Y^{*}})\propto
\[\prod_{j=1}^{n({\bf p})}
\kappa_{e_{j,n}}({\mbox e}^{-vh}\rho)\]\pi_{\theta+n}(dv).
\EndEq

\begin{thm}
\noindent
{\em (i)}
The marginal distribution of ${\bf Y}$ is 
\Eq
\pi(d{\bf Y}):=\int_{0}^{\infty}\left[\prod_{j=1}^{n({\bf p})}\PR(dY^{*}_{j}|{\mbox
  e}^{-vh}\rho,\eta)\right]\pi_{\theta+n,\theta}(dv,{\bf p}|\rho,\eta).
\EndEq
This implies that the quantity $\pi_{\theta+n,\theta}({\bf
  p}|\rho,\eta)$ is an EPPF. When $\rho$ does not depend on $y$
then the marginal distribution of ${\bf Y}$ is expressible as 
\Eq
\pi(d{\bf Y}):=\pi_{\theta+n,\theta}({\bf
  p}|\rho,\eta)\prod_{j=1}^{n({\bf p})}H(dY^{*}_{j}).
\EndEq

\noindent
{\em (ii)}
Statement (i) combined with Theorem 5.1 imply that the posterior
distribution of $N$, $\mu$ given ${\bf Y}$ is
identical to the mixtures,
\Eq
\int_{{\rr}^{+}\times{{\mathcal S}^{n({\bf p})}}}{\mathcal P}(dN,d{\bf J}|{\mbox e}^{-vh}\rho,\eta,{\bf Y}^{*})
\pi_{\theta+n}(dv|{\bf p},{\bf Y^{*}});
\int_{{\rr}^{+}}{\mathcal P}(d\mu|{\mbox e}^{-vh}\rho,\eta,{\bf Y}^{*})
\pi_{\theta+n}(dv|{\bf p},{\bf Y^{*}})
\label{moopost}
\EndEq
respectively.

\noindent
{\em (iii)}
Statement (ii) implies that the posterior distribution of $P$ is
determined by either of the laws in ~\mref{moopost}. Combined
with the mixture representations this implies that the distribution of
$P$ given ${\bf Y}$ is equivalent to the distribution of the random
measure
\Eq
p^{*}_{n}P^{*}(\cdot)+(1-p^{*}_{n})\frac{\sum_{j=1}^{n({\bf
  p})}h(J_{j,n})\delta_{Y^{*}_{j}}(\cdot)}{\sum_{j=1}^{n({\bf p})}h(J_{j,n})}
\label{rmpost}
\EndEq
where the distributions of $\mu^{*}$ and $(J_{j,n})$ given $V$,${\bf
  Y}$ is specified by ${\mathcal P}(d\mu|{\mbox e}^{-vh}\rho,\eta,{\bf Y}^{*})$
and $V|{\bf Y}$ is $\pi_{\theta+n}(dv|{\bf p},{\bf Y^{*}})$ 
\end{thm}
\Proof
Given the mixture representations in Theorem 5.1, the result follows by an 
 appeal to Fubini's theorem
which identifies the posterior laws of $N$,${\mu}$, $P$ as mixtures relative to the
marginal distribution of ${\bf Y}$. That is, for instance ${\mathcal
  P}(dN):=\int_{{\mathcal Y}^{n}}{\mathcal P}(dN|{\bf Y})\pi(d{\bf
  Y})$. (More formally one could evaluate the Laplace functional of $N$ on
both sides). Hence it suffices to identify the marginal distribution of ${\bf
  Y}$ and then apply a simple algebraic rearrangement in the mixture
representations of $N$ given in statements (ii) and (iii) of Theorem
5.1. The
identification of $\pi(d{\bf Y})$ is straightforward.
\EndProof

The conclusion that $\pi_{\theta+n,\theta}({\bf
p}|\rho,\eta)$ is an EPPF perhaps requires further discussion as the
technique I used may be a bit unfamiliar. Essentially from the 
theory of exchangeability if
$Y_1,\ldots,Y_{n}|P$ are iid $P$ then the marginal distribution of ${\bf
  Y}=({\bf Y}^{*},{\bf p})$ is exchangeable. Hence once the unique
  values ${\bf Y}^{*}$ are exposed an integration with respect to $\eta$
  leaves only a marginal dsitribution of ${\bf p}$ which must be 
an EPPF regardless of
whether or not the unique 
$\{Y^{*}_{1},\ldots,Y^{*}_{n({\bf p})}\}$ are iid $H$. That is
regardless of whether or not $P$ is a
species sampling model as  described in Pitman (1996). What is
    lost is the 1-1 correspondence between $({\bf p},H)$ and the random
    probability measure $P$. For instance, it is conceivable that the
    EPPF ${PD}({\bf p}|\alpha,\theta)$ could be embedded in
a model $P$ which is not ${\mathcal P}_{\alpha, \theta}(dP|H)$. In
    that case the ${Y^{*}_{1},\ldots,Y^{*}_{n({\bf
        p})}}$ cannot be iid $H$.
In other words an EPPF can
  always be found by working with a $P$ model, finding the joint
  marginal distribution of ${\bf
  Y}$, and then marginalizing over the unique values. This is obvious
when $P$ is a species sampling model. However, it is a simple matter to
verify the addition rules given in Pitman (1995a,b, 1996) for an EPPF by
applying the Bayesian idea of a prediction rule. In particular for
each $n$ evaluate,
$$
E[P({\mathcal Y})|{\bf Y}]:=E[p^{*}_{n}|{\bf Y}]+E[(1-p^{*}_{n})|{\bf Y}]:=1.
$$
Proper manipulation of the middle expression will yield the obvious rules.
 
\subsection{Results for ${PD}(\alpha,\theta)$}

A description of the ${\mathcal
  PD}_{\alpha, \theta}(d\mu|\eta)$ laws is now given. 
Throughout this section the notation $\mu_L$, $T_L$ etc will be used
to denote the (conditional law) of $\mu$, $T$ depending on a random variable $L$.
\subsubsection{Distributional properties of ${\mathcal {PD}}_{\alpha,\theta}(d\mu|\eta)$} 
\begin{cor}
\noindent
{\em (i)}
If $N$ is ${\mathcal P}(dN|\rho_{\alpha},\eta)$ and $h(s):=s$, then
for $\theta>0$, 
${\ps}(dN|\rho_{\alpha},\eta, \theta)$ is such that $N|V=$ is a Poisson
random measure with intensity,
\Eq
\rho_{\alpha,v}(ds)\eta(dy):={\mbox e}^{-vs}\rho_{\alpha}(ds)\eta(dy)
\EndEq
corresponding to the L{\'e}vy measure of a generalised gamma process. The
density of $V$ is 
\Eq 
\tau_{\theta}(dv|\rho):={\frac{1}{c_{\alpha, \theta}\Gamma(\theta)}}v^{\theta-1}{\mbox
  e}^{-Kv^{\alpha}}dv.
\EndEq
By a change of variable the distribution of $L:=V^{\alpha}$ has a gamma distribution with parameters
$({\frac{\theta}{\alpha}}, K)$, That is, $L$ is ${\mathcal G}({\frac{\theta}{\alpha}},K)$.
The factor $K$ is determined in part by the total mass of $\eta$. It
can be dispensed with by re-scaling.
\end{cor}

\Remark
Now the connection to Pitman and Yor (1992, section 3, p. 335-336) should be more
transparent.[See also Pitman and Yor (2001, Theorem 3)]. The exponential law arises by setting $K=1$(or by
rescaling $L$) and the
choice of $\theta=\alpha$. 
\EndRemark

\begin{prop}
Corollary 5.2 implies that the distribution of ${\mathcal {PD}}_{\alpha, \theta}(d\mu|\eta)$ with
respect to the mixing distribution of $L$ is,
\Eq
{\mathcal {PD}}_{\alpha, \theta}(d\mu|\eta):=\int_{0}^{\infty}{\mathcal
  P}(d\mu|\rho_{\alpha, {L}^{1/\alpha}},\eta){\mathcal G}(dL|
{\theta/\alpha},K)
\EndEq
In other words $\mu|L$ is a generalized gamma random measure
with L{\'e}vy measure 
$\rho_{\alpha,L^{1/\alpha}}(ds)\eta(dy)$ and $L$ is a gamma random
  variable with parameters $({\frac{\theta}{\alpha}},K)$.
When $K=1$ and $\theta=\alpha$, $L$ is a standard  exponential random
variable. The Laplace functional of ${\mathcal {PD}}_{\alpha, \theta}(d\mu|\eta)$ can be expressed as, 
\Eq
\int_{\mathcal M}{\mbox e}^{-\mu(g)}{\mathcal
  {PD}}_{\alpha, \theta}(d\mu|\eta):=\int_{0}^{\infty}
{\mathcal L}_{\mu}(g|L,\alpha)
{\mathcal G}(dL|{\theta/\alpha},K)
\label{laplacePD}
\EndEq
where 
$$
{\mathcal L}_{\mu}(g|L,\alpha):=
{\mbox e}^{-{\frac{1}{\alpha}}\int_{\mathcal Y}
{(g(y)+{L}^{{\frac{1}{\alpha}}})}^{\alpha}\eta(dy)}{\mbox
  e}^{KL}
$$
is the conditional Laplace functional of $\mu$ given $L$. 

\noindent
{\em (ii)} Furthermore, suppose that conditioned on $L$, $\mu$ is
  multiplied by $L^{1/\alpha}$. Then the unconditional law of
  $L^{1/\alpha}\mu_{L}$ is given by its Laplace functional,
\Eq
E[{\mbox e}^{-L^{1/\alpha}\mu_{L}(g)}]:=\int_{0}^{\infty}{\mathcal
  L}_{\mu}(L^{1/\alpha}g|L,\alpha){\mathcal G}(dL|{\theta/\alpha},K)
:={\[\int_{\mathcal Y}
{(g(y)+1)}^{\alpha}H(dy)\]}^{-\theta/\alpha}
\EndEq

\noindent
{\em (iii)} 
Statement (ii) implies that the distribution of $L^{1/\alpha}T_{L}$
is ${\mathcal G} (\theta)$.
\end{prop}

\Remark 
Statement (iii) should be compared with $T$ of Pitman and Yor (1997,
Proposition 21). Note that the explicit expression for the
Laplace functional in (ii) is obtained via Brix's (1999) expression
for the generalised gamma measure.
\EndRemark

Now noting that,
\Eq
\prod_{j=1}^{n({\bf p})}\kappa_{e_{j,n}}({\mbox
      e}^{-vh}\Omega):=v^{-(n-n({\bf p})\alpha)}{\eta({\mathcal
      Y})}^{n({\bf p})}
\prod_{j=1}^{n({\bf p})}\Gamma(e_{j,n}-\alpha),
\EndEq 
the results below are easily deduced.
\begin{prop}
\noindent
{\em (i)}
For all $\theta>-\alpha$ and $n\ge 1$, the ${\mathcal PD}_{\alpha, \theta}(d\mu|\eta)$ law is
representable as the random measure, 
\Eq
\mu^{*}(\cdot)+\sum_{j=1}^{n({\bf p})}J_{j,n}\delta_{Y^{*}_{j}}(\cdot)
\label{mixmeasure}
\EndEq
where $\mu^{*}$ given a random variable $L$ is a generalised gamma random measure with 
intensity $\rho_{\alpha, {L}^{1/\alpha}}(ds)\eta(dy)$. The $(J_{j,n})$ given
$(L,{\bf Y}^{*})$  are independent of $\mu^{*}$ with respective Gamma
distributions ${\mathcal G}{(e_{j,n}-\alpha,
  {L}^{1/\alpha}})$. ${L}$ given ${\bf p}$ is ${\mathcal G}(n({\bf p})+{\theta/\alpha}, K)$
for some constant $K$. In particular by cancellation one can set
$K=\alpha$ or $K=1$. 
Conditionally independent of ${L}$, $\{Y^{*}_{1},\ldots, Y^{*}_{n({\bf p})}\}$ are iid
$\eta(\cdot)/\eta({\mathcal Y}):=H(\cdot)$. 

\noindent
{\em (ii)}
The distribution of $\bf
p$ is the EPPF, $ {PD}({\bf p}|\alpha,\theta)$. In
addition the marginal law of $\mu^{*}|{\bf p}$  is ${\mathcal PD}_{\alpha,
  \theta+n({\bf p})\alpha}(d\mu|\eta)$, which follows since 
similar to~\mref{laplacePD} its Laplace functional is
\Eq
\int_{0}^{\infty}
{\mathcal L}_{\mu}(g|L,\alpha)
{\mathcal G}(dL|{n({\bf p})+{\theta/\alpha}},K)
\label{laplacePD2}
\EndEq
\noindent
{\em (iii)}
Hence
the random probability measure
$$
\frac{\mu^{*}(\cdot)}{\mu^{*}({\mathcal Y})}
$$
given ${\bf p}$ is ${\mathcal PD}_{\alpha,
  \theta+n({\bf p})\alpha}(dP|H)$. Denote the random probability measure with this
law as $P_{\alpha,\theta+n({\bf p})\alpha}$.

\noindent
{\em (iv)}
Given $L,{\bf p}$, the random variables $G_{j,n}:=L^{1/\alpha}J_{j,n}$ are
independent ${\mathcal G}(e_{j,n}-\alpha)$ independent of $L$ and
$\mu^{*}$. Moreover given ${\bf p}$ the distribution of
$L^{1/\alpha}\mu^{*}_{L}$ is given by the Laplace functional
$$
{\[\int_{\mathcal Y}
{(g(y)+1)}^{\alpha}H(dy)\]}^{-(\theta+n({\bf p})\alpha)/\alpha}
$$
and $L^{1/\alpha}T^{*}_{L}$ given ${\bf p}$ has distribution
${\mathcal G}(\theta+n({\bf p})\alpha)$
\end{prop}
 
\Remark
Suppose that $K=1$, and $n=1$, then in particular for the stable case ${\mathcal {PD}}_{\alpha,0}(d\mu|\eta)$,
it follows that $L$ is exponential $(1)$. For the ${\mathcal
  PD}_{\alpha,\alpha}(d\mu|\eta)$, $L$ is ${\mathcal
  G}{(2, 1)}$. 
\EndRemark

\subsubsection{Identities for ${PD}(\alpha,\theta)$}
The propositions above are  now applied to derive an alternate representation for the distribution of
the general ${PD}(\alpha, \theta)$ and related models. This, in particular generalizes the
right-hand side of the construction of a ${PD}(\alpha, \alpha)$ model
given in Pitman and Yor (2001, Example 8, eq. (33)) to previously
unknown ones for the general ${PD}(\alpha, \theta)$ model. 
 
Notice that using the change of variable $u=vs$,  ${\mbox
  e}^{-vs}\rho_{\alpha}(ds)$ transforms to the Levy measure
\Eq
v^{\alpha}{\mbox e}^{-u}\rho_{\alpha}(du)
\EndEq
which yields the equivalence of the sets
\Eq
\{y: \int_{y}^{\infty} {\mbox
  e}^{-vs}\rho_{\alpha}(ds)\leq x\}:=\{u:\int_{u}^{\infty} {\mbox
  e}^{-s}\rho_{\alpha}(ds)\leq x/{{v}^{\alpha}}\}
\label{inverseid}.
\EndEq
Define,
\Eq
\Lambda^{-1}(x):= \inf\{u:\int_{u}^{\infty} {\mbox
  e}^{-s}\rho_{\alpha}(ds)\leq x\},
\EndEq
and set $\Gamma_{j}:=\sum_{i=1}^{j}E_{i}$ for $(E_i)$ a collection of
independent standard exponential random variables. In additon define
for all $\theta>0$
\Eq
\Sigma_{\theta/\alpha}:=\sum_{j=1}^{\infty}\Lambda^{-1}(\Gamma_{j}/L)
\EndEq
where $L$ is a ${\mathcal
  G}(\theta/\alpha)$ random variable, independent of $(E_{i})$.
Now using the change of variable $L=V^{\alpha}$,~\mref{inverseid} and
Proposition 5.3, 5.4, it is easy to see from an application of
Khintchine's
 (1937) Inverse L{\'e}vy method,[see also
Ferguson and Klass~(1972), Wolpert and Ickstadt (1998b), Sato (1999), Rosinski
(2001), and Banjevic, Ishwaran and Zarepour (2002)], that the 
following identities hold;
\begin{prop}(Distributional representations for $PD(\alpha,\theta)$)
Choose $0<\alpha< 1$,

\noindent
{\em (i)}
then for $\theta>0$,the distribution of the sequence 
\Eq
(\Lambda^{-1}(\Gamma_{j}/L)/{\Sigma_{\theta/\alpha}};\quad j=1,2,\ldots),
\EndEq
is ${PD}(\alpha, \theta)$.

\noindent
{\em (ii)}
Let $(Z_j)$ denote an
iid sequence with distribution H chosen idependently of $L$ and
$(E_j)$. 
If $\theta>0$, then equivalent to (i), the random probability measure, 
\Eq
P_{\alpha, \theta}(\cdot):=\sum_{j=1}^{\infty}{\frac{\Lambda^{-1}(\Gamma_{j}/L)}{\Sigma_{\theta/\alpha}}}\delta_{Z_{j}}(\cdot)
\EndEq
is ${\mathcal PD}_{\alpha,\theta}(dP|H)$.

\noindent
{\em (iii)}
For all $\theta>-\alpha$ and $n\ge  1$, a ${\mathcal
  {PD}}_{\alpha,\theta}(dP|H)$ random probability measure is
representable as, 
\Eq
P_{\alpha, \theta}(\cdot):={\frac{\Sigma_{n({\bf p})+\theta/\alpha}}{
\Sigma_{n({\bf p})+\theta/\alpha}+\sum_{j=1}^{n({\bf p})}J_{j,n}}}
P_{\alpha,\theta+n({\bf p})\alpha}(\cdot)
+{\frac{\sum_{j=1}^{n({\bf p})}J_{j,n}\delta_{Y^{*}_{j}}(\cdot)}{
\Sigma_{n({\bf p})+\theta/\alpha}+\sum_{j=1}^{n({\bf p})}J_{j,n}}}
\EndEq
where, for clarity,
$\Sigma_{n({\bf p})+\theta/\alpha}:=\sum_{j=1}^{\infty}\Lambda^{-1}(\Gamma_{j}/L)$. 
$L|{\bf p}$ is gamma distributed with parameters $(n({\bf p})+{\theta/\alpha},1)$, $(J_{j,n})|L,{\bf p}$ are respectively independent ${\mathcal
  G}{(e_{j,n}-\alpha, {L}^{1/\alpha}})$ and  $(Y^{*}_{j})|{\bf p} $ are
iid  $H$. Note that $(J_{j,n})$ are not independent of $\Sigma_{n({\bf p})+\theta/\alpha}$ 
\end{prop}

Now to obtain another representation of $P_{\alpha,\theta}$, which also
serves to directly recover Pitman's~(1996) description of the posterior
distribution. Notice from Proposition 5.3 and 5.4 that given ${\bf p}$ the
following equivalence in distribution holds for each $n\ge 1$;
\Eq
\frac{L^{1/\alpha}T_{L}}{L^{1/\alpha}T_{L}+
L^{1/\alpha}\sum_{j=1}^{n({\bf p})}J_{j,n}}
:={\frac{G_{\theta+n({\bf
    p})\alpha}}{G_{\theta+n({\bf p})\alpha}+\sum_{j=1}^{n({\bf
        p})}G_{j,n}}}.
\EndEq
[This is also true for n=0 by Proposition 5.3]
The key point is that the quantity above is  independent of the
mixing distribution on $L$ for all $n$. Moreover, ${L^{1/\alpha}}T_{L}:=G_{\theta+n({\bf
    p})\alpha}$ is ${\mathcal
  G}(\theta+n({\bf p})\alpha)$ and independent of the gamma random
varibles $(G_{j,n})$ as defined in
Proposition 5.4. Hence the following result, 
\begin{prop}
For all $\theta>-\alpha$, 
\Eq
P_{\alpha,\theta}(\cdot):=p_{n}P_{\alpha,\theta+n({\bf
    p})\alpha}(\cdot)+(1-p_{n})
{\frac{\sum_{j=1}^{n({\bf
        p})}G_{j,n}\delta_{Y^{*}_{j}}(\cdot)}{\sum_{j=1}^{n({\bf
        p})}G_{j,n}}}
\label{priorpost}
\EndEq
where 
$$
p_n:={\frac{G_{\theta+n({\bf
    p})\alpha}}{G_{\theta+n({\bf p})\alpha}+\sum_{j=1}^{n({\bf
        p})}G_{j,n}}}.
$$
Given ${\bf p}$ the quantity above does not depend on the mixing
distribution $L$. Hence given ${\bf Y}=({\bf Y}^{*},{\bf p})$
the posterior distribution of a $P$ which is ${\mathcal PD}_{\alpha,
  \theta}(dP|H)$ is immediately seen to be equivalent to the random
measure on the right of~\mref{priorpost} when $(Y^{*}_{j})$ and ${\bf
  p}$ are fixed. This corresponds exactly with the posterior
distribution described in Pitman (1996)  
\end{prop}

\Remark
It is certainly obvious that one could use Pitman's~(1996)
posterior characterization to obtain the simple mixture representation in
proposition 5.6. The main point however is really how the previous
descriptions (which are less obvious)
led up to this result. Moreover, none of the arguments appealed to the
stick-breaking representation of ${\mathcal PD}_{\alpha, \theta}(dP|H)$. 
\EndRemark

\Remark
The analysis  of the ${PD}(\alpha,\theta)$ models revealed various
independence structures via a simple transformation. In general this will
not be the case but there are certainly many instances where an appropriate
transformation of the $(J_{j,n})$ will render them independent of
$\mu^{*}$ and the mixing distribution on $V$. Such an operation should
simplify the analysis. One might try this with the intensities described
in Pitman and Yor (2001).
\EndRemark

\Remark
I wonder what if any interpretation does an adjustment to 
the left hand-side of Pitman and Yor (2001, Theorem 3, eq. (19)) 
have when $\epsilon_{0}$ is
replaced by what one might guess from Proposition 5.3 to be a
${\mathcal G}(\theta/\alpha)$ random variable
\EndRemark

\subsection{Results for the Dirichlet Process and generalised gamma
process}
Results for the gamma process with shape
$\eta(\cdot):=\theta H(\cdot)$, that is ${\mathcal
  PD}_{0,\theta}(d\mu|\eta)$, follow by using $\rho_{0,1}$ in place of
$\rho_{\alpha}$. In this case $T$ is ${\mathcal G}(\theta)$, and 
\Eq
\prod_{j=1}^{n({\bf p})}\kappa_{e_{j,n}}({\mbox
      e}^{-vh}\Omega):=(1+v)^{-n}{\theta}^{n({\bf p})}
\prod_{j=1}^{n({\bf p})}\Gamma(e_{j,n}).
\EndEq
Which yields readily,
\begin{prop}
\noindent
\Eq
\pi_{n;0}(dv,{\bf p}):= 
{PD}({\bf p}|\theta)\tau_{\theta;n}(dv)
\EndEq
where 
\Eq
\tau_{\theta;n}(dv):=\frac{\Gamma(\theta+n){(1+v)}^{-(n+\theta)}v^{n-1}dv}
{\Gamma(\theta)}
\label{gamdensity}
\EndEq
which implies that $V$ and ${\bf p}$ are independent. Note the density
$\tau_{\theta;n}(dv)$ is well defined provided $\theta>0$. 

\noindent
{\em (ii)}
The $(J_{j,n})|V,{\bf p}$ are independent ${\mathcal G}(e_{j,n},1+V)$ and the
distribution of $\mu^{*}|V,{\bf p}$ is a (simple) weighted gamma process,
i.e. has intensity $\rho_{0,1+V}\theta H$, and
is independent of ${\bf p}$.  

\noindent
{\em (iii)}
Hence it follows that given $V$ and ${\bf p}$
\Eq
{(V+1)}\[\mu^{*}(\cdot)+\sum_{j=1}^{n(\bf
  p)}J_{j,n}\delta_{Y^{*}_{j}}\]
\label{gamN}
\EndEq
is a {\it mixture} of gamma processes independent of $V$. That is
additionally given ${\bf Y}^{*}$, the measure in ~\mref{gamN} is a gamma process
with shape  
$\theta H(\cdot)+\sum_{j=1}^{n({\bf
    p})}e_{j,n}\delta_{Y^{*}_{j}}(\cdot)$. This fact serves to recover
the well-known result of Ferguson (1973) for the posterior
distribution of the Dirichlet process. 
\end{prop}

\Remark
It is not so surprising that the gamma/Dirichlet model is such that
the mixing distribution $V$
and ${\bf p}$ are independent. It is also perhaps true, given the
properties of ${PD}(\theta)$, that this is
the only species sampling model with this property.
\EndRemark

The arguments above may be applied to the generalised gamma model with
intensity~$\rho_{\alpha, b}$.  In this case, from section 3, it follows
that,
\Eq
\prod_{j=1}^{n({\bf p})}\kappa_{e_{j,n}}({\mbox
      e}^{-vh}\Omega):=(v+b)^{-(n-n({\bf p})\alpha)}
{\theta}^{n({\bf p})}\prod_{j=1}^{n({\bf p})}\Gamma(e_{j,n}-\alpha)
\EndEq 
where it is assumed that $\eta({\mathcal Y}):=\theta$.
This leads to, a joint density of $V,{\bf p}$ specified as 
\Eq
{\theta}^{n({\bf p})}\left[\prod_{j=1}^{n({\bf p})}\Gamma(e_{j,n}-\alpha)\right]
\frac{{(v+b)}^{-n+n({\bf p})\alpha}v^{n-1}
{\mbox e}^{-[(v+b)^{\alpha}-b^{\alpha}]K}
dv}{\Gamma(n)}
\EndEq
In addition the 
$(J_{j,n})|V,{\bf p}$ are ${\mathcal G}(e_{j,n}-\alpha, b+V)$
and $\mu^{*}|V,{\bf p}$ is a generalised gamma process with L{\'e}vy
measure ${\rho_{\alpha,b+V}}$. Hence $(b+V)J_{j,n}$ are ${\mathcal
  G}(e_{j,n}-\alpha)$. The Laplace functional of
$(V+b)\mu^{*}|V,{\bf p}$ is 
\Eq
{\mbox e}^{-{\frac{1}{\alpha}}{(v+b)}^{\alpha}\int_{\mathcal Y}
[{(g(y)+1)}^{\alpha}-1]\eta(dy)}
\EndEq

\Remark 
In order to incorporate larger classes of models for $P$ one could use 
a weighted Poisson law 
\Eq
{\mathcal Q}(dN|\rho,\eta):={\frac{w(N){\mathcal
      P}(dN|\rho,\eta)}{E[w(N)]}}
\EndEq
for an arbitrary integrable function $w$. This will be used in Section
8. 
\EndRemark

\section{Distributions of joint linear functionals of P; variations of the
  Markov-Moment problem}
This section is a continuation of the previous one. Here, it is shown
that the joint Cauchy-Stieltljes transform of linear functionals of
$P$, which are ${\ps}(dP|\rho, \eta, \theta)$ and ${\mathcal P}(dP|\rho, \eta)$, 
is equivalent to expressions involving
the Laplace functional of random measures $V\mu_{V}$. The
precise meaning of $V\mu_{V}$ will be clear from the context below.
 The method of proof, given the
results in the the previous section, is an easy
extension of the beautiful approach used by Tsilevich, Vershik and Yor (2000)
for the Dirichlet process and the general ${\mathcal P}_{\alpha,
  \theta}(dP|H)$ family.  
The results given here  represents the most general ones that I know
of. More importantly the explicit relationship between $P$ and
${V\mu_{V}}$ is a new insight. See Kerov (1998) for many implications
of this type of result. 

Let $f_l$ denote real-valued functions on $\yy$ and define
$Pf_l=\int_\yy f_l(y)P(dy)$ for $l=1,\ldots,q$. In addition let $z_{l}$
for $l=1,\ldots,q$
denote non-negative scalars. In this section the
calculation of the following transform(in relation to Laplace
functionals of $V\mu_{V}$) is discussed; 
\Eq
\int_{\mathcal M^{*}}\frac{1}{1+\sum_{l=1}^{q}z_{l}Pf_{l}}
{\ps}(dP|\rho,\eta,\theta)
\label{Mellintransform}
\EndEq
for all $\theta>-\alpha$. Again when $\theta=0$ this coincides with
the ${\mathcal P}(dP|\rho,\eta)$ laws.
The quantity~\label{Mellintransform} characterizes the joint 
distribution of $(Pf_1,\ldots,Pf_q)$. 
Kerov and Tsilevich (1998) in the case
of ${\mathcal P}_{\alpha, \theta}(dP|H)$ used
combinatorial arguments to obtain the moment expressions in the case of
the Dirichlet and two-parameter models to yield extensions of the
Markov-Krein identity for $(Pf_1,\ldots,Pf_q)$.  Their results extend the work
of Cifarelli and Regazzini~(1990) for the case of the distribution of
the mean functional, ${\int y P(dy)}$, when $P$ has a Dirichlet
process law. The mean case is also discussed in Diaconis and
Kemperman~(1996) where in addition the result for the joint
distribution of functionals like $(Pf_1,\ldots,Pf_q)$ was proposed as
an open problem. Tsilevich~(1997) establishes the case for the mean
with respect to the general two-parameter processes. These results
used hard analytic techiques which would not be easily extendable to a
general scenario. However, recently Tsilevich, Vershik and
Yor (2000) devise a beautiful simple proof of the corresponding result in the
case of the Dirichlet process and the general two-parameter extension
via Laplace functionals. Given the results in section 5.1 it is a simple
matter to extend their result to the general class of probability
measure ${\ps}(dP|\eta,\rho,\theta)$. That is, following closely their
approach, relationships between~\label{Mellintransform} and the
Laplace functional of $V\mu_{V}$ are established.
The results below follow by rewriting
$$
\frac{1}{1+\sum_{l=1}^{q}z_{l}Pf_{l}}:=\frac{T}{T+\sum_{l=1}^{q}z_{l}\mu(f_{l})}
$$
and applying Corollary 5.1.
\subsection{Joint Cauchy-Stieltjes transforms and Laplace functionals}
\begin{prop}
For $\theta>0$, 
\Eq
\int_{\mathcal M^{*}}{\(1+\sum_{l=1}^{q}z_{l}Pf_{l}\)}^{-\theta}
{\ps}(dP|\rho,\eta,\theta)
:={\mathcal L}_{V\mu_{V}}(\sum_{l=1}^{q}z_{l}f_{l}|\rho,\eta,\theta)
\EndEq
where 
$$
{\mathcal L}_{V\mu_{V}}(\sum_{l=1}^{q}z_{l}f_{l}|\rho,\eta,\theta)
:=
\int_{0}^{\infty}
{\mathcal L}_{\mu}(v\sum_{l=1}^{q}z_{l}f_{l}|{\mbox e}^{-vh}\rho,\eta)
\pi_{\theta}(dv|\rho,\eta)
$$
\end{prop}
\begin{prop}
For $\theta+n>0$ and $n\ge 1$,

\noindent
{\em (i)}
$$
\int_{\mathcal M^{*}}{\(1+\sum_{l=1}^{q}z_{l}Pf_{l}\)}^{-(\theta+n)}
{\ps}(dP|\rho,\eta,\theta)
:=\sum_{{\bf p}}\int_{0}^{\infty}
{\mathcal L}_{\mu}(v\sum_{l=1}^{q}z_{l}f_{l}|{\mbox
  e}^{-vh}\rho,\eta,{\bf p})
\pi_{\theta+n,\theta}(dv,{\bf
  p}|\rho,\eta)
$$
\noindent
{\em (ii)}
The expressions in (i) are equal to;
$$
\sum_{{\bf p}}\int_{{\mathcal Y}^{n({\bf p})}}\[\int_{0}^{\infty}{\mathcal L}_{\mu}(v\sum_{l=1}^{q}z_{l}f_{l}|{\mbox
  e}^{-vh}\rho,\eta,{\bf p},{\bf Y}^{*})\pi_{\theta+n,\theta}
(dv|\rho,\eta,{\bf p},{\bf Y}^{*})\]\pi(d{\bf Y})
$$
which indicates that the posterior Cauchy-Stieltjes transform, 
$$
\int_{\mathcal M^{*}}{\(1+\sum_{l=1}^{q}z_{l}Pf_{l}\)}^{-(\theta+n)}
{\ps}(dP|\rho,\eta,\theta,{\bf Y})
$$
is equal to,
\Eq
\int_{0}^{\infty}{\mathcal L}_{\mu}(v\sum_{j=1}^{q}z_{l}f_{l}|{\mbox
  e}^{-vh}\rho,\eta,{\bf p},{\bf Y}^{*})\pi_{\theta+n,\theta}
(dv|\rho,\eta,{\bf p},{\bf Y}^{*})
\EndEq
Note again the various relationships to the (posterior) laws of $V\mu_{V}$.
\end{prop}
Now setting $g(y):=\sum_{l=1}^{q}z_{l}f_{l}(y)$ in the Laplace
transform of $L^{1/\alpha}\mu_{L}$ in statement (ii) of Proposition
(5.3) yields, 
$$ 
{\[\int_{\mathcal Y}
{(\sum_{l=1}^{q}z_{l}f_{l}(y)+1)}^{\alpha}H(dy)\]}^{-\theta/\alpha}
$$
Hence the  result of Tsilevich, Vershik and Yor (2000) for the
${\mathcal P}_{\alpha,\theta}(dP|H)$ family is recovered. 
However the relationship to $L^{1/\alpha}\mu_{L}$ is not noted in their
work. 

\subsection{A remark on moment calculations} 
As mentioned previously, Kerov and Tsilevich (1998) used nontrivial
combinatorial arguments to calculate the joint moments of
$(Pf_1,\ldots,Pf_q)$ in the case of ${\mathcal
  P}_{\alpha,\theta}(dP|H)$. Here similar to Ishwaran and James (2001a)
for species sampling models it is demonstrated that one can easily
obtain the relevant moment calculations by using Theorem 5.2. This
calculation will only be presented for the case where $\rho$ does not
depend on $y$, 
The task is to calculate
$$
E\[\prod_{l=1}^{q}(Pf_l)^{n_l}\]=\int_\mm
\[\prod_{l=1}^{q}\prod_{i=1}^{n_l}\int_\yy f_l(y_{i,l})\,P(dy_{i,l})\]{\ps}(dP|\rho,\eta,\theta)
$$
Now
analogous to~\mref{momcalc1} an application of Theorem 5.2 yields the result
\Eq
E\[\prod_{l=1}^{q}{(P(f_{l}))}^{n_l}\]=
\sum_{\bf p}{\pi_{\theta+n,\theta}({\bf p}|\rho,\eta)}\prod_{j=1}^{n({\bf p})}
\int_{{\mathcal Y}}
\[\prod_{l=1}^{q}{f_{l}}^{e^{l}_{j,n}}(u)\]H(du)
\label{momcalc2}
\EndEq
 
\section{Posterior Calculus for Extended Neutral to the Right
  processes}
In this section I focus on the concept of neutral to the right processes(NTR)
originally proposed in Doksum (1974). 
The Dirichlet process is the
most notable member of this class. Here, a new natural extension of the NTR concept to more
abstract spaces is given. It is then shown
how Proposition 3.1 can be used to  yield the most transparent and simplest
posterior analysis of such models. This includes the case of survival data models subject to right censoring
when a NTR prior is used  or synonymously when L{\`e}vy process
priors are used to model the cumulative hazard. Additionally, using
Proposition 3.1 a change of measure formula is established which relates 
Beta/Dirichlet processes to their more complex Beta/Beta-Neutral(Stacy) generalizations.

\Remark 
Doksum (1974, Theorem 3.1) establishes essentially 
a 1-1 correspondence between NTR processes
and exponential functions of subordinators. See below for explict
details. This fact seems not to be
widely noticed by probabilists investigating problems where models
under the  latter description arise. One consequence is that the
calculus that is described below for NTR's can be exploited in other
areas besides Bayesian nonparametrics. 
Here I will omit the
drift component. It is a simple matter to make adjustments starting
from the obvious modification of Lemma 2.1. (see Remark 2).
Doksum (1974, Corollary 3.2), establishes the almost
sure discreteness of NTR's under the condition that the drift
component is zero.  
\EndRemark 

\Remark
The notation $\Lambda$ will be used in this section to denote a random
cumulative hazard measure. The dependence of quantities $F_{0}$, $A_{0}$
on $\rho$, $\eta$ will be supressed. The arguments (s),(t) will be used
to denote time as is usual in survival analysis. The argument (u)
plays the role of (s) in the previous sections. 
\EndRemark

First the orginal definition proposed by Doksum is given
\Definition(Doksum~(1974))
A random distribution function $F$ on the positive real line is said to
be {\it neutral to the right} if for each $k>1$
$t_1<t_2\ldots<t_k$, there exists non-negative
independent random variables  $V_1,\ldots, V_k$ such that the
vectors satisfy,  
\Eq
{\mathcal L}\{(F(t_{1}),F(t_2)-F(t_1), F(t_{k})-F(t_{k-1}))\}={\mathcal
  L}\{(V_1,V_{2}(1-V_{1}),\ldots, V_{k}\prod_{j=1}^{k}(1-V_i))\},
\EndEq
where ${\mathcal L}$ denotes the law.
\EndDefinition
In the special case where $F$ is a Dirichlet process with shape
$\theta F_{0}(\cdot)$ then each increment 
$F(t_k)-F(t_k-1)$ is ${\mathcal B}(\theta F_{0}((t_{k-1},t_{k}]);
\theta[1-F_{0}((t_{k-1},t_{k}])$.
Doksum discusses various equivalences and implications of this
definition. From Theorem 3.1 of Doksum it follows that $F$ is a
NTR process if and only if for $t\ge 0$,
\Eq
S(t)=1-F(t)={\mbox e}^{-Z(t)},
\label{NTR1} 
\EndEq
where $Z$ is an increasing Levy process satisfying
$Z(0)=0$ and $\lim_{n\rightarrow \infty}Z(t)=\infty$ . The analysis
here will consider subordinators $Z$ without a drift component. 
In other words $Z$,
is a completely random measure on $(0,\infty)$ with associated intensity
$\rho_{z}(du|y)\eta(dy)$ for $(u,y)\in (0,\infty)\times(0,\infty)$.
Ferguson (1974) shows that a Dirichlet process with finite shape
measure, $\theta F_{0}(\cdot)$, results if 
\Eq
\rho_{z}(du|y)\eta(dy)=\frac{1}{1-{\mbox e}^{-u}}{\mbox
  e}^{-u\theta F_{0}(y,\infty)}du\theta F_{0}(dy)
\label{Dirichlet}
\EndEq
It follows from the theory of
product integration that an NTR process can also be represented as 
\Eq
S(t)=\Prodi_{u\leq t}\(1-{\Lambda(du)}\)
\label{product}
\EndEq
where ${\Lambda}$ denotes a cumulative hazard which is further
modelled as a completely random measure with intensity 
$\rho_{{\Lambda}}(u|y)ds\eta(dy)$ for $(u,y)\in
(0,1]\times(0,\infty)$.  
The symbol 
$\prodi ,$
denotes the product integral which has played
a primary role in survival analysis. In
particular the Kaplan-Meier estimator for $S $, Kaplan and Meier 
(1958), is obtained by replacing  $\Lambda $ by its empirical
counterpart, the Nelson-Aalen estimator. See the text by Andersen,
Borgan, Gill and Keiding (1993) for further elaboration.
Gill and Johansen (1990) discuss in detail
the properties of the product integral. The product integral is
also expressible as, 
\EqArray
\Prodi_{[0,t]}\left( 1-\Lambda (dv,\yy)\right)=\exp (-\Lambda^{c}(t))\prod
_{[0,t]}\left( 1-\Lambda_{d}({v})\right),
\EndEqArray
where $\Lambda^{c}$ denotes the
continuous part of  $\Lambda$.
Suppressing the dependence on
$\rho$, $\eta$,
$E[{\Lambda(t)}]=A_{0}(t)$ where $A_{0}$ denotes a prior cumulative hazard
specification. It follows that, 
\Eq
E[S(t)]:=1-F_{0}(t)=\Prodi_{u\leq t}\(1-E[{\Lambda(du)}]\):={\mbox e}^{-A_{o}(t)}. 
\label{expected}
\EndEq 

The restriction of the jumps  of the process
to $[0,1]$ ensures that $\Lambda$ is an element in the space of
cumulative hazards and hence $S$ is a proper survival function. 
Hjort
(1990) first proposed working directly with L{\'e}vy priors on the space of
cumulative hazards which is more in
line with the frequentist counting process analysis of event-history
models [See Aalen(1975, 1978) and Andersen, Borgan, Gill and Keiding (1993)]. Hjort (1990) shows that if $\Lambda$ is specified to be a Beta
process then it is a conjugate model with respect to right censoring. 
Hjort (1990, section 7A), under a Beta process specification for ${\Lambda}$
in~\mref{product},  also defines a class of {\it generalised Dirichlet
  processes} on ${\mathcal R^{+}}$ . He shows that the Dirichlet process is a special case of
this model by setting $c(s)=\theta F_{0}{([s,\infty))}$. 
 In summary,  Bayesian nonparametric methods for the simple survival setting subject
to censoring have been
discussed following the framework of Ferguson's (1973, 1974) (see also 
Freedman (1963) and Fabius (1964)) Dirichlet
process in the works of Doksum (1974), Susarla and van Ryzin (1976), Blum
and Susarla (1977), Ferguson and Phadia (1979), Lo (1993), Doss (1994), and Walker and
Muliere (1997) among others. The methods discussed above operate by placing
a Dirichlet or more general NTR prior on the unknown survival or
distribtuion function. An alternative but essentially equivalent approach 
involves working with priors on the cumulative hazard measure
discussed in Hjort (1990), Lo (1993) and most recently Kim (1999). However,
unlike the simplicity of the Dirichlet process discussed in Ferguson (1973, 1974) for complete
data models, the technical aspects of these models appear to be
formidable. 
Moreover, the technical arguments used do not easily extend to more complex settings. In
particular,  the prior to posterior characterizations
given in Ferguson and Phadia (1979)(see also Doksum (1974)), are only
developed for distribution functions  on the real line. In addition
very little is known about the marginal and partition based structures.
It will be shown that an alternate representation makes the calculus for NTR
processes indeed straightforward.  

Note importantly that there is a 1-1 correspondence between each $Z$
and $\Lambda$. Formally, the
L{\'e}vy measure of $Z$ arises as the image of
$\rho_{{\Lambda}}(du|y)\eta(dy)$ via the map $(u,y)$ to
$(-\log(1-u),y)$. For further discussion see Dey~(1999) and Dey, Erickson and
Ramamoorthi~(2000).  An important consequence, which has not been exploited in
this context, is the following identity, which holds in distribution for each $n$ where $N$
is ${\mathcal P}(dN|\rho_{\Lambda}, \eta)$. Define for $v>0$,
 \Eq
{\fs}_{v-}(u,y):=-I\{v>y\}\log(1-u)
\EndEq
then, 
\Eq
S(v-):={\mbox e}^{Z(v-)}:={\mbox e}^{-N({\fs}_{v-})}
\label{identity},
\EndEq
where the law of $N$ is ${\mathcal P}(dN|\rho_{\Lambda}, \eta)$. 

\subsection{Definition of Extended NTR processes}
Suppose that $(T_i,X_i)$ denotes a marked pair of random variables
on $\rr^{+}\times{\mathcal X}$ with distribution $F(ds,dx)$. In this
section an answer is provided to the open question of how to extend an
NTR process from $\rr^{+}$ to more general marked Polish spaces. This
provides for instance a new class of Bayesian models for multivariate
survival models.  While indeed it is easy to extend $Z$ or $\Lambda$
to more abstract spaces, the representation in~\mref{NTR1} or~\mref{product} do not
immediately suggest an obvious extension for $F$. The Dirichlet process which is defined over abitrary
spaces is a notable exception. However, James and Kwon (2000) recently propose a method which extends
the Beta-Neutral prior of Lo (1993), and by virtue of the
equivalences, the Beta-Stacy process in Muliere and Walker (1997) and Beta distribution function discussed
in Hjort~(1990, section 7A), to a spatial setting. A general definition
can be deduced from elements of their construction. 
A definition for $F$ on ${\mathcal R}^{+}\times{\mathcal X}$ is
facilitated by the usage of its associated {\it hazard measure } on ${\mathcal
  R}^{+}\times{\mathcal X}$. From Last and Brandt~(1995, A5.3), it
follows that such a measure always exists and is defined by,
\Eq
\Lambda(ds,dx):=I\{t>0\}\frac{F(ds,dx)}{S(s-)}.
\EndEq
In particular, $\Lambda(ds,{\mathcal X}):= \Lambda(ds)$ and hence
\Eq
S(s-):=\Prodi_{u< s}\(1-{\Lambda(du,{\mathcal X})}\).
\EndEq

An extended NTR is defined as, 
\Definition(Extended Neutral to the Right Process)
Let $\Lambda$ denote a completely random measure with intensity
$\rho_{{\Lambda}}(du|s)\eta(ds,dx)$ for
$(u,s,x)\in (0,1]\times(0,\infty)\times{\mathcal X}$. Furthermore, the
intensity measures is chosen such that 
$$
A_{0}(ds,dx):=E[\Lambda(ds,dx)|\rho_{\Lambda},\eta]:=\[\int_{0}^{1}u\rho_{{\Lambda}}(du|s)\]
\eta(ds,dx)
$$
is a {\it hazard measure}. [Denote the marginal cumulative hazard
$A_{0}(ds,{\mathcal X}):=A_{0}(ds)$].
Then an {\it Extended Neutral to the Right} 
process on ${{\mathcal R}^{+}}\times{\mathcal X}$ is defined for $t>0$ and
each B,an arbitrary measurable set in ${\mathcal X}$,  
\Eq
F(t,B):=\int_{0}^{t}S(s-)\Lambda(ds,B):=\int_{0}^{t}\Prodi_{u<
  s}\(1-{\Lambda(du)}\)\Lambda(ds,B)
\EndEq
In particular, $F(ds,dx):=S(s-)\Lambda(ds,dx)$. The law of $F$ is
denoted ${\mathcal PN}(dF|\rho_{\Lambda},\eta)$. The random quantities
$S(s-)$ and $\Lambda(ds,dx)$ are independent for each $s$ and arbitrary $x$ and 
\Eq
E[F(ds,dx)]:=E[S(s-)]E[\Lambda(ds,dx)]:={\mbox e}^{-A_{0}(s)}A_{0}(ds,dx)
:=F_{0}(ds,dx).
\EndEq
\EndDefinition

\Remark The definition of an extended neutral to the right process
yields, as a special case, a class of random probability measures on arbitrary
spaces ${\mathcal X}$, defined as 
\Eq
F(dx):=\int_{0}^{\infty}S(s-)\Lambda(ds,dx).
\EndEq
For instance this expression offers another identity for a Dirichlet
process on ${\mathcal X}$. 
\EndRemark

\Remark 
The definition is expressed in terms of $\Lambda$ rather than $Z$
extended to $\rr^{+}\times{\mathcal X}$ due
to the natural interpretation of a hazard measure. For instance a
description of $F(ds,dx)$ is not easily seen using $Z$. Nonetheless
for each $\Lambda$ and $F$ in $\rr^{+}\times{\mathcal X}$ one can associate
a $Z$ on $\rr^{+}\times{\mathcal X}$ with again the 
Levy measure of $Z$ arising as the image of
$\rho_{{\Lambda}}(du|y)\eta(dy,dx)$ via the map $(u,y)$ to
$(-\log(1-u),y)$. 
\EndRemark

\Remark
When $B={\mathcal X}$ it is obvious that $F(\cdot,{\mathcal X})$ is
an NTR. In addition, due to the complete randomness properties of
$\Lambda$, $F$ satisfies, 
$$
{\mathcal L}\{(F(t_{1},B),F(t_{2},B)-F(t_{1},B), F(t_{k},B)-F(t_{k-1},B))\}={\mathcal
  L}\{(V_{1,B},V_{2,B}(1-V_{1}),\ldots, V_{k,B}\prod_{j=1}^{k}(1-V_i))\},
$$
where for each $j$, $V_j:=V_{j,B}+V_{j,B^{c}}$ and $V_{i:B}$ is independent of
$V_{j,C}$ for $i\neq j$ and $B, C$ arbitrary. A posterior process will
be called an extended NTR process if the NTR properties are preserved.
\EndRemark

\subsection{Posterior distributions and moment formulae} 
Now suppose that one observes n-iid observations $(T_i, X_i)$ from   
$\prod_{i=1}^{n}F(dT_{i},dX_{i})$ and consider the following joint
models,
\Eq
\[\prod_{i=1}^{n}F(dT_{i},dX_{i})\]{\mathcal PN}(dF|\rho_\Lambda, \eta)
\quad{\mbox and}\quad
\[\prod_{i=1}^{n}F(dT_{i},dX_{i})\]{\mathcal P}(d\Lambda|\rho_\Lambda,
  \eta).
\EndEq
Here we will work with $\Lambda$, and the equivalent expression
\Eq
\[\prod_{i=1}^{n}S(T_{i}-)\Lambda(dT_{i},dX_{i})\].
\label{NTRlike}
\EndEq
If one assumes the classical  univariate right censoring applied to
the marked data as in Huang and Louis (1998) then the likelihood
under censoring takes the form
\Eq
\[\prod_{l=1}^{m}S(C_{l}-)\]\[\prod_{i=1}^{n}S(T_{i}-)\Lambda(dT_{i},dX_{i})\]
\label{NTRlike2}
\EndEq
where $C_1,\ldots,C_m$ denote $m$ independent censoring times which
indicate that there are random variables $T_{n+1},\ldots,T_{m+n}$
where it is only known that they exceed the respective censored
times. Under this assumption no information for the marks associated
with the censored points $T_{n+1},\ldots,T_{m+n}$ is
available.   

The primary focus will be to deduce the posterior distribtion of
both $F$ and $\L$ and related characteristics of the marginal
distribution of $(T_1,X_1),\ldots, (T_n, X_n)$. This will complete the
necessary disintegration which will allow one to apply both the
(extended) models for $\Lambda$ and $F$ to  a large class of data
structures beyond univariate right censoring. In fact it will become
clear from the form of~\mref{NTRlike}  that
analysis of right censoring data for NTR is really the same affair as
analysis of the complete data model. 

First, for $i=1,\ldots,n$ define 
${\ys}_{T_{i}-}(s):=I\{T_i>s\}$ and similarly for $l=1,\ldots, m$ define
${\ys}_{C_{l}-}(s):=I\{C_l>s\}$. In addition for $i=1\ldots,n$ define
${\fs}_{T_{i}-}$ satisfying, 
\Eq
(1-u)^{{\ys}_{T_{i}-}(s)}:={\mbox e}^{-{\fs}_{T_{i}-}(u,s)},
\label{fdef}
\EndEq
and for $l=1,\ldots,m$,let ${\fs}_{C_{l}-}$ be defined similarly.
Then for each $(n,m) \ge 0$ 
\Eq
(1-u)^{Y_{n,m}(s)}:={\mbox e}^{-f_{n,m}(u,s)}
\EndEq
where $Y_{n,m}(s):=\sum_{i=1}^{n}{\ys}_{T_{i}-}(s)+
\sum_{l=1}^{m}{\ys}_{C_{l}-}(s)$. 
The quantities $f_{n,m}$,
${\fs}$ are special cases of the functions defined in Proposition
3.1 and 3.2. 
It follows from the identity in~\mref{identity} that 
\Eq
S(T_{i}-):={\mbox e}^{-N({\fs}_{T_{i}-})}\quad{\mbox and}\quad
\prod_{i=1}^{n}S(T_{i}-)={\mbox e}^{-N({f}_{n})}
\label{keyresult}
\EndEq
where $N$ is ${\mathcal P}(dN|\rho_{\Lambda},\eta)$. Now the
likelihood~\mref{NTRlike2} can be rewritten as 
\Eq
{\mbox e}^{-N({f}_{n,m})}\prod_{i=1}^{n}\Lambda(dT_{i},dX_{i})
\label{identity2}
\EndEq
This representation,~\mref{identity2}, in combination with Proposition 3.1 and Theorem
3.1 yields the posterior distributions for $\Lambda$, $F$, $Z$ subject
to possible right censorship.
\begin{thm}
The posterior distribution of $\Lambda$ given the
model~\mref{NTRlike2} is ${\mathcal P}(d\Lambda| {\mbox
 e}^{-f_{n,m}}\rho_{\Lambda}, \eta, {\bf T}, {\bf X})$. 
That is, $\Lambda$ is equivalent in distribution to the random measure
\Eq
\Lambda(\cdot|Y_{n,m})+\sum_{j=1}^{n(\bf
  p)}{J_{j,n}}\delta_{T^{*}_{j},X^{*}_{j}}(\cdot),
\label{posthazard}
\EndEq  
where the law of $\Lambda(\cdot|Y_{n,m})$ is ${\mathcal P}(d\Lambda| {\mbox
  e}^{-f_{n,m}}\rho, \eta)$ indicating that its intensity measure is,
\Eq
{(1-u)}^{Y_{n,m}(s)}\rho_{\Lambda}(du|s)\eta(ds,dv)
\label{hazardupdate}.
\EndEq
The $J_{j,n}$ are (conditionally) mutually independent random variables
with distribution, for each $j$, depending on $T^{*}_j,
Y_{n,m}(T^{*}_{j})$, defined as in (32) as, 
\Eq
\PR(J_{j,n}\in
du|{\mbox e}^{-f_{n,m}}\rho_{\Lambda},T^{*}_{j}):=\frac{u^{e_j,n}{(1-u)}^{Y_{n,m}(T^{*}_{j})}\rho_{\Lambda}(du|T^{*}_j)}{\kappa_{e_j,n}({\mbox e}^{-f_{n,m}}\rho_{\Lambda}|T^*_j)}, 
\EndEq
and are conditionally independent of $\Lambda(\cdot|Y_{n,m})$.

\noindent
{\em (ii)} 
The posterior distribution of $F$ is still an extended NTR
with distribution ${\mathcal PN}(dF|{\mbox
  e}^{-f_{n,m}}\rho_{\Lambda},\eta,{\bf T},{\bf X})$
determined by replacing the random measure $\Lambda$
with~\mref{posthazard}.

\noindent
{\em (iii)} 
A posterior distribution of $Z$ is equivalent to the law of the
random measure 
$$Z(\cdot|Y_{n,m})+\sum_{j=1}^{n({\bf p})}(1-{\mbox
  e}^{-J_{j,n}})\delta_{T^{*}_{j},X^{*}_{j}}(\cdot),$$
where the L{\'e}vy measure
of  $Z(\cdot|Y_{n,m})$ arises as the image of
${(1-u)}^{Y_{n,m}(s)}\rho_{{\Lambda}}(du|s)\eta(ds,dx)$ via the map $(u,s)$ to
$(-\log(1-u),s)$
\end{thm}
When $m:=0$, the results correspond to a complete data model.

\Proof
First set $\mu:=\Lambda$, and 
$w(\Lambda):={\mbox e}^{-N(f_{n,m})}$. Now apply statement (ii) of
Proposition 3.1. 
\EndProof

\Remark
Theorem 7.1 serves to extend the results for the univariate setting to
a spatial setting. The previous works for the univariate setting do
not use explicitly a partition based representation. More importantly
the method of proof, which is new, is quite transparent. 
\EndRemark

\Remark
Note also
that due to the non-atomic nature(continuity) of $\eta(ds)$ the quantity 
$Y_{n,m}(s)$ in~\mref{hazardupdate} can be replaced by
$$
Y^{+}_{n,m}(s):=\sum_{i=1}^{n}I\{s\leq T_{i}\}+\sum_{l=1}^{m}I\{s\leq C_{l}\}.
$$
In other words calculations with respect to $\Lambda(\cdot|Y_{n,m})$
should be understood to be equivalent to those with respect to 
$\Lambda(\cdot|Y^{+}_{n,m})$. This does not apply to the distribution
of the jumps $(J_{j,n})$.
\EndRemark

Little is known in general about the {\it explicit} joint moment structure
of Neutral to the Right models. The results in Doksum (1974) are
rather vague.  In recent works expressions for the mean and variance
are given.
The formulae in Proposition 3.2 can be used to easily
obtain various equivalent expressions which goes well beyond a
variance calculation. This is seen by setting 
$w_{i}(\Lambda):={\mbox e}^{-N({\fs}_{T_{i}-)}}$ for $i=1,\ldots, n$
, and other obvious equivalences. For brevity, I will only present a
result which yields the relevant joint structure and EPPF for these
models. Such results do not appear in the literature mentioned above.

Define,
\Eq
{\tilde
  A}_{n,m}(\infty):=\int_{0}^{\infty}\int_{0}^{1}{(1-{(1-u)}^{Y^{+}_{n,m}(s)})}\rho_{\Lambda}(du|s)\eta(ds)
\EndEq
In addition
for $i=1,\ldots, n$ and $m\ge 0$ define, 
\Eq
A_{i-1,m}(t):=\int_{0}^{t}\int_{0}^{1}u{(1-u)}^{Y^{+}_{i-1,m}(s)}\rho_{
\Lambda}(du|s)\eta(ds). 
\EndEq
When $m=0$, denote $A_{i-1,m}$ as $A_{i-1}$. Now recall that,  
\Eq
\kappa_{e_{j,n}}({\mbox
  e}^{-f_{n,m}}\rho_{\Lambda}|T^{*}_{j})=\int_{0}^{1}
u^{e_{j,n}}{(1-u)}^{Y_{n,m}(T^{*}_{j})}\rho_{\Lambda}(du|T^{*}_{j}),
\EndEq

\begin{prop}

\noindent
{\em (i)}
For $m>0$, the (prior) mean calcuation for $\[\prod_{l=1}^{m}S(C_{l}-)\]\prod_{i=1}^{n}S(T_{i}-)$ is
\Eq
E[{\mbox e}^{-N({f}_{n,m})}|\rho_{\Lambda},\eta]=
{\mbox  e}^{-{\tilde A}_{n,m}(\infty)}={\mbox e}^{-A_{o}(T_{1})}\left[\prod_{i=2}^{n}
{\mbox e}^{-A_{i-1}(T_{i})}\right]\prod_{l=1}^{m}{\mbox e}^{-A_{n,l-1}(C_{l})}
\EndEq
\noindent
{\em (ii)}
When $m=0$, the joint marginal distribution of
$((T_i,X_i))$ can be expressed as 
\Eq
\prod_{j=1}^{n({\bf p})}\[\prod_{i\in C_{j}}{\mbox
  e}^{-A_{i-1}(T^{*}_{j})}\]\kappa_{e_{j,n}}({\mbox
  e}^{-f_{n}}\rho_{\Lambda}|T^{*}_{j})\eta(dT^{*}_{j},dX^{*}_{j})
\EndEq
Adjustments for the censored case are obvious.

\noindent
{\em (iii)}
From statement (i) it follows that the corresponding EPPF has the
form, 
\Eq
\int_{{\mathcal T}^{n({\bf p})}}
\prod_{j=1}^{n({\bf p})}\[\prod_{i\in C_{j}}{\mbox
  e}^{-A_{i-1}(T^{*}_{j})}\]\kappa_{e_{j,n}}({\mbox
  e}^{-f_{n}}\rho_{\Lambda}|T^{*}_{j})\eta(dT^{*}_{j})
\EndEq
\end{prop}
\begin{prop} 
\noindent
Using an algebraic rearrangement the joint marginal distribution and
EPPF can be rewritten respectively as,
\Eq
\pi({\bf p}|\rho_{\Lambda},{\bf T}^{*})\prod_{j=1}^{n({\bf
      p})}F_{0}(dT^{*}_{j},dX^{*}_{j})
{\mbox and}\quad
\int_{{\mathcal T}^{n({\bf p})}}
\pi({\bf p}|\rho_{\Lambda},{\bf T}^{*})\prod_{j=1}^{n({\bf
      p})}F_{0}(dT^{*}_{j}).
\EndEq
Where,
\Eq
\pi({\bf p}|\rho_{\Lambda},{\bf T}^{*}):=\prod_{j=1}^{n({\bf p})}
{\frac{\[\prod_{i\in C_{j}}{\mbox
  e}^{-A_{i-1}(T^{*}_{j})}\]\kappa_{e_{j,n}}({\mbox
  e}^{-f_{n}}\rho_{\Lambda}|T^{*}_{j})}{{\mbox
  e}^{-A_{0}(T^{*}_{j})}\kappa_{1}(\rho_{\Lambda}|T^{*}_{j})}}
\EndEq
\end{prop}

\Remark 
Naturally if one were interested in actually generating such
partitions etc, then a modification of the algorithm in Section 2.3
can be used. For this one could use expressions for the prediciton
rule or conditional moment measures which are readily obtainable from
an application of Proposition 3.2 combined with Theorem 7.1. 
\EndRemark

\Remark
The propositions above combined with Theorem 7.1 yield expressions for
the posterior disintegrations. It is now straightforward to 
obtain posterior characterizations for mixtures of {\it extended } 
NTR models based on kernels $(K_{i})$. This framework allows for much
more complex structures than right censoring. I have not seen general
mixtures of NTR models proposed in the literature. 
\EndRemark
\subsection{Absolute continuity of general Beta,Beta-Neutral/Stacy models to
a canonical Beta processs or Dirichlet process}

The general construction of the extended NTR models is an extension of
(presently unpublshed work) James
and Sehyug Kwon (2000). In that work the authors extend Lo's~(1993) Beta-Neutral
survival and cumulative hazard processes to the spatial setting. 
Lo (1993) derives these based on the following explict construction for the hazard
\Eq
\Lambda_{\tau,\beta}(t):=\int_{0}^{t}{\frac{\mu_{\tau}(ds)}{\mu_{\tau}
([s,\infty))+\mu_{\beta}([s,\infty))}},
\EndEq
where $\mu_{\tau}$,$\mu_{\beta}$ are independent gamma processes with
shape measure $\tau$ and $\beta$ on ${\rr}^{+}$, which as noted in
Lo (1993) yields an explicit contruction of Hjort's (1990) Beta
cumulative hazard process. James and Kwon (2000) extend this
definition by simply extending the gamma processes to a spatial setting.
Moreover they show that such models can be always derived from a
Dirichlet process on an even larger space. In other words take a two
parameter gamma process, say $\mu_{\tau,\beta}$,
on ${\rr^{+}}\times{\mathcal X}\times\{0,1\}$, such that 
$\mu_{\tau,\beta}(ds,dx,\{1\}):=\mu_{\tau}(ds,dx)$ etc. 
A corresponding Dirichlet process can be defined as 
\Eq
P_{\tau,\beta}(ds,dx,{\Delta}):={\frac{\mu_{\tau,\beta}(ds,dx,
    \Delta)}{\mu_{\tau, \beta}({\rr^{+}}\times{\mathcal
      X}\times\{0,1\})}}.
\EndEq
Then the extension of
James and Kwon (2000) can be deduced from the extended Beta-Neutral hazard measure,
\Eq
\Lambda_{\tau,\beta}(ds,dx):={\frac{\mu_{\tau,\beta}(ds,dx,\{1\})}
{\mu_{\tau,\beta}([s,\infty)\times{\mathcal
      X}\times\{0,1\})}}
\label{BNhazard}
\EndEq
where $\tau$ and $\beta$ are now measures on ${\rr}^{+}\times{\mathcal
  X}$. If $\beta$ is set to zero in~\mref{BNhazard} then the corresponding
(extended) Neutral to the Right process is a Dirichlet process with
shape parameter $\tau$ which, without loss of generality, is set to
$\theta F_{0}$. 
By virtue of the essential equivalence of the
Beta-Neutral process to the Beta-Stacy process of Walker and Muliere (1997)
and the Beta distribtion function of Hjort (1990, Section 7A), the
procedure of James and Kwon (2000) includes these models as well. In
addition they showed that the (posterior) conjugacy of the Beta-type
models on ${\rr}^{+}$ is preserved under the right censored data
spatial model discussed earlier. However their
technique relied on very special properties of the Dirichlet process
which is quite different than what has been presented in the previous
sections. Here a new result
is established which shows how one may transform a Dirichlet process
or simple Beta process to a more general one via a change of measure.
A consequence is that
the calculus for such models follows from the calculus for the
Dirichlet process plus an application of Proposition 3.1. 

Recall that a Beta process on ${\rr}^{+}$ with parameters $c(s)$ and
$A_{0}(ds,dx)$  yields a Dirichlet process if and only if $c(s):=
\theta F_{0}([s,\infty))$. Such a Beta process can be thought of 
as a canonical Beta process.

\begin{prop}   
Let ${\mathcal P}(d\Lambda|c,A_{o})$ denote the law of a Beta
process with parameters $c$ and $A_{0}(ds,dx)$. Let $Z$ denote the
corresponding L{\'e}vy process defined via the map $(u,y)$ to
$(-\log(1-u),y)$ and define a decreasing function $\beta$  
on $\rr^{+}$ such that, 
$$
T_{\beta}:=\int_{0}^{\infty}\beta(v)Z(dv)<\infty. 
$$
Then the following disintegration holds
\Eq
{\mbox e}^{-T_{\beta}}{\mathcal P}(d\Lambda|c,A_{0}):=
{\mathcal P}(d\Lambda|c+\beta,A_{0})E[{\mbox e}^{-T_{\beta}}]
\EndEq
where $-\log E[{\mbox e}^{-T_{\beta}}]$ is
\Eq
\int_{0}^{\infty}\int_{0}^{1}(1-{(1-u)}^{\beta(s)}){u}^{-1}(1-u)^{c(s)-1}duA_{0}(ds).
\EndEq
\end{prop}

The proof follows by using the alternate representation
\Eq
T_{\beta}:={\mbox e}^{-N(f_{\beta})}
\EndEq
where $f_{\beta}(u,s):=-\beta(s)\log(1-u)$ and $N$ 
has a Poisson law corresponding to ${\mathcal P}(d\Lambda|c,A_{o})$.
An interesting feature of this result is that one can obtain quite
easily an alternate expression for the EPPF of a Beta-Neutral/Stacy model by 
first applying the result for a Dirichlet process. That is, 

\begin{prop}
Suppose that $F$ is a NTR process determined by the Beta process with
parameters $c^{*}:=\theta F_{0}+\beta$, then the EPPF is given by,
\Eq
{PD}({\bf p}|\theta)\prod_{j=1}^{n(\bf
  p)}\int_{0}^{\infty}
\[{\frac{\Gamma(e_{j,n}+\theta F_{0}([y,\infty)))
\Gamma(\theta F_{0}([y,\infty))+\beta(y))}{\Gamma(\theta F_{0}([y,\infty)))
\Gamma(e_{j,n}+\theta F_{0}([y,\infty))+\beta(y))}}\]F_{0}(dy).
\EndEq
\end{prop}

\Remark
The change of measure in  Proposition 7.3 can of course be
extended easily to non-Beta processes. In general, the updated law
corresponds to a random hazard measure with L{\'e}vy measure ${\mbox
  e}^{-f_{\beta}}\rho$. Note that much more general choices of $f_{\beta}$ 
can be used via Proposition 3.1.
\EndRemark

\Remark
The correspondence between the Beta-Stacy process of Walker and
Muliere~(1997) and Hjort's (1990, section 7A) process is noted explicitly in Dey~(1999) and Dey, Ericson,
and Ramamoorthi~(2000).
The equivalence between the Beta-Neutral process and Hjort's(1990) process
was noted in Lo (1993). Given the gaps in the literature it is apparent
that Beta-Neutral processes are not
as well known as their equivalent counterparts. This seems to be
caused by the title of Lo (1993) which concerns a Censored Data
Bayesian Bootstrap.  
\EndRemark

\Remark
NTR processes seem to arise naturally in coalescent theory. See in
particular Pitman(1999, Proposition 26) which is not a Dirichlet
process. Similar types of processes with drift appear in Bertoin (2001).
\EndRemark 

\section{Posterior Distributions of Normalised processes and Poisson-Kingman models}
In this section I briefly discuss calculations for probability
measures $P$ defined using the weighted Poisson distribution
$Q(dN|\rho, \eta)$ which are more in line with the results and methods
used in Perman, Pitman
and Yor (1992) and Pitman and Yor (1992) and Pitman (1995b).  Here
descriptions of pertinent quantities will be given in terms of the
biased jumps denoted as ${\js}$. One will see that the forms of the
results appear quite different than Section 5. For completion the
definition of Poisson-Kingman models based on {\it length biased sampling
} presented in Pitman (1995b) is given, 

\begin{defin}(Pitman (1995b))
Let $P_i=(J_i/T)$ be a ranked discrete distribution derived
from the ranked points of a  Poisson Process with L{\'e}vy density
$\rho$ (not depending on y) of random lengths
$J_{1}\ge J_{2}\ge \cdots\ge 0$ by normalizing their lengths by their
sum which is $T$. Let (\~P$_{j}$) be a size-biased permutation of (P$_{i}$)
and let \~J$_{j}=(T$\~P$_{j}$) be the corresponding size-biased
permutation of the ranked lengths $(J_{i})$. The law of the sequence
($P_{i}$)  will be called the {\it Poisson-Kingman distribution with L{\'e}vy density
$\rho$}, and denoted  $PK(\rho)$. Denote by $PK(\rho|t)$ the regular
conditional distributioon of (P$_{i}$) given $(T=t)$ constructed
above. For a probability distribution $\gamma$ on $(0,\infty)$, let 
\Eq
PK(\rho,\gamma):=\int_{0}^{\infty}PK(\rho|t)\gamma(dt)
\label{poissonkingman}
\EndEq
be the distribution on the space of (P$_{i}$),[which is the space of
decreasing sequences of positive real numbers with sum 1]. Call
$PK(\rho,\gamma)$ the {\it Poisson-Kingman distribution with L{\'e}vy density
$\rho$ and mixing distribution $\gamma$ }.
\end{defin}
Pitman (1995b), points out that knowledge of the conditional law $PK(\rho|t)$
allows one to generate explicit results for distributions $PK(\rho,
\gamma)$, in particular the corresponding EPPF,  
by simply mixing over different candidate densities, $\gamma(dt)$, for
$T$. The case of the two parameter Poisson-Dirichlet family is
explained in detail in Pitman (1995b) as mentioned in section 5.

The introduction of the measure $Q(dN|\rho,\eta)$ serves to
incorporate fully the Poisson-Kingman idea while maintaining the
approach of Poisson calculus at the level of $N$. 
Some properties of this class, which shall become clear in the next
section are now described.
When $\rho$ is homogeneous and
$w(N):=w(N(\cdot,{\mathcal Y}))$ then $P$ is a species sampling
model. Additionally when $w(N):=g(T)$ and $h(s):=s\in(0,\infty)$ then
the random atoms of $P$, $(P_i)$ are $PK(\rho,\gamma)$. In that case 
\Eq
\gamma(dt):=\frac{g(t)f_{T}(dt)}{E[g(T)|\rho,\eta]}
\EndEq
where $g$ is nonegative and integrable but otherwise arbitrary, $f_{T}$
denotes the density of $T$ with respect to ${\mathcal P}(dN|\rho,\eta)$.
In other words the change of measure at $N$ via ${\mathcal
Q}(dN|\rho,\eta)$ induces the appropriate
change of measure at the level of the PK measure etc.
The results
below follow from a straightforward application of Lemma 2.2, details
are omitted.

\subsection{Posterior characterizations}
In the theorem below it is stated that
the posterior law of $P|{\bf Y}$, denoted as ${\mathcal P}(dP|{\bf
  Y})$,  corresponds to the law of a random measure
defined as, 
\Eq
P^{*}_{n}(\cdot):=
R_{n({\bf p})}P(\cdot)+
(1-R_{n({\bf p})})\sum_{j=1}^{n({\bf p})}{\frac{h({\mbox {\~J}}_{j})}
{\sum_{j=1}^{n({\bf p})}h({\mbox {\~J}}_{j})}}\delta_{Y^{*}_{j}}(\cdot),
\label{postpm}
\EndEq
where, 
$
R_{n({\bf p})}:=T/\[T+\sum_{j=1}^{n({\bf p})}h({\mbox
      {\~J}}_{j})\],
$
\begin{thm} Let $\{Y_1,\ldots, Y_n\}$ be iid $P$ where the prior law of $P$ is
 determined by the weighted Poisson measure $Q(dN|\rho, \eta)$. Then, the  posterior
 distribution of $P|{\bf Y}$ is equivalent to the distribution of
 the random measure $P^{*}_{n}$ defined in~\mref{postpm} whose law
 is now determined by the joint probability measure 
\Eq
{\mathcal Q}(dN, {\bf {\mbox
{\~J}}}\in d{\bf s}|{\bf Y},h)\propto
w(N+\sum_{j=1}^{{n(\bf p)}}\delta_{s_j,Y^{*}_j}){\mathcal
  P}(dN|\rho, \eta)
{{\(T
+\sum_{j=1}^{n(\p)}h(s_j)\)}^{-n}}
\prod_{j=1}^{n(\bf
      p)}h(s_{j})^{e_{j,n}}\rho(ds_{j}|Y^{*}_j)
\label{posteriorP2}
\EndEq
The marginal distribution of $\{Y_1,\ldots,Y_{n}\}$ corresponds to the
un-normalised term on the right hand side of ~\mref{posteriorP2}, integrated over ${\mathcal
  P}(dN|\rho, \eta)$,
and multiplied by $\prod_{i=1}^{n({\bf p})}\eta(dY^{*}_{i})$. 
\end{thm}

\Remark
Notice that the expression above remains complicated even if $W(N)$ is
replaced by $1$. The result is more in line with Pitman, Perman and Yor
(1992) and Pitman (1995b). In comparison with section 5 the law of the
random measure $P^{*}_{n}$ is a bit more complex as it is decribed via
the biased jumps ${\js}$ which have a much more complex (non-independent) joint distribution. Under moment
conditions other forms of the posterior are easily obtained.  
\EndRemark

Now letting 
\Eq
T_{n(\bf p)}:=T-\sum_{j=1}^{n(\bf p)}h(s_{j})
\EndEq
it is quite clear that one can gain further interpretation by adapting
the descriptions given in Perman, Pitman and Yor (1992, section 4) and
Pitman (1995b). Note in particular when $n=1$, and $w(N):=g(T)$ and
$\rho$ is homogeneous an evaluation of the expectation of $P$ reveals
the structural distributions, 
\Eq
\PR_{g}({{\mbox {\~J}}}_{1}\in d{\bf s},T_1\in dt_1)=
{\frac{g(t_{1}+s)f_{T}(t_{1})dt_{1}{h(s)}\Omega(s)ds}
{{\(t_{1}+s\)E[g(T)]}}}={\frac{g(t_{1}+s)}{E[g(T)]}}\PR ({{\mbox
  {\~J}}}_{1}\in d{\bf s},T_1\in dt_1)
\EndEq
A suitable change of variable yields 
\Eq
\PR_{g}(T\in dt,T_1\in dt_1)=
{\frac{g(t)}{E[g(T)]}}
\PR (T\in dt,
T_1\in dt)={\frac{g(t)f_T(t)}{E[g(T)]}}\PR (T_1 \in dt_1|T=t)
\EndEq
and more generally for each $n$,
$$
\PR_{g}\({\bf p},{\bf J}\in d{\bf s},T_{n({\bf
    p})}\in dt\)=
{\frac{g(t +\sum_{j=1}^{n({\bf
      p})}h(s_j))}{E[g(T)]}}
\PR\({\bf p}, {\bf J}\in d{\bf s},T_{n({\bf
    p})}\in dt\).
$$
The joint distributions $\PR ({{\mbox
  {\~J}}}_{1}\in d{\bf s},T_1\in dt_1)$ and $\PR\({\bf p}, {\bf J}\in d{\bf s},T_{n({\bf
    p})}\in dt\)$  are given in
Pitman (1995b) and Perman, Pitman and Yor (1992).  
One can for instance obtain formulae for
the joint distribtuion of $(T,T_1,\ldots, T_n)$
via $\(g(t)/E[g(T)]\)\PR(T\in dt, T_1\in dt_{1},\ldots,T_{n}\in dt_{n})$ using 
Perman, Pitman and Yor (1992, Theoerem 2.1). All of these facts correspond with the
$PK(\rho,\gamma)$ concept. See additionally Pitman and Yor (1997,
Proposition 47). 

Combining the description in Pitman (1995b, Lemma 5
and equation(30)) with Theorem 8.1 the next result follows.

\begin{cor}
Suppose that $w(N+\sum_{j=1}^{n({\bf p})}\delta_{s_j,Y^*_j})
:=w(N;\sum_{j=1}^{n({\bf p})}\delta_{s_j})$ does not depend on
$Y^*$ for all $n$. Then a joint law of 
of ${\bf p},{\bf {\mbox {\~J}}},{\bf Y}^{*}$ is defined as;
\Eq
{\frac{1}{E[w(N)]}}{\[\int_{\mathcal M}{\frac{w(N;\sum_{j=1}^{n({\bf p})}\delta_{s_j}){\mathcal
    P}(dN|\rho,\eta)}{{(T+\sum_{j=1}^{n({\bf p})}h(s_j))}^{n}}}\]}
{\prod_{j=1}^{n(\bf
      p)}h({s_j})^{e_{j,n}}\rho(ds_{j}|Y^{*}_{j})}\eta(dY^{*}_{j}). 
\EndEq
As a consequence, the distribution of ${\bf Y}|{\bf {\mbox
    {\~J}}},T_{n({\bf p})},{\bf p}$ is conditionally independent of
$T_{n({\bf p})}$ such that the sequence 
$\{Y_1,\ldots, Y_{n}\}$ consists of $n({\bf p})$ unique values 
$\{Y^{*}_1,\ldots,Y^{*}_{n({\bf p})}\}$ which are independent with
distribution
\Eq
\PR (dY^{*}_{j}|{\js}_{j},{\bf p}))\propto
\rho({\mbox {\~J}}_{j}|Y^{*}_j)\eta(dY^{*}_{j}).
\label{djy}
\EndEq
for $j=1,\ldots, n({\bf p})$. Additionally the joint distribution of
${\bf {\mbox {\~J}}},{\bf p}$ is
\Eq
\PR({\bf {\mbox {\~J}}}\in d{\bf s},{\bf
  p})={\frac{1}{E[w(N)]}}{\[\int_{\mathcal
    M}{\frac{w(N;\sum_{j=1}^{n({\bf p})}\delta_{s_j}){\mathcal
    P}(dN|\rho,\eta)}{{(T+\sum_{j=1}^{n({\bf p})}h(s_j))}^{n}}}\]}
{\prod_{j=1}^{n(\bf
      p)}h({s_j})^{e_{j,n}}\int_{\mathcal Y}\rho(ds_{j}|y)\eta(dy)}. 
\EndEq
If additionally $\rho$ is homogeneous then
\Eq
\PR (d{\bf Y})=p_{\rho}(e_1,\ldots, e_{n({\bf p})})
\prod_{j=1}^{n(\p)}H(dY^{*}_{j}),
\EndEq
where, the EPPF is 
\Eq
p_{\rho}(e_1,\ldots, e_{n({\bf p})})={\frac{1}{E[w(N)]}}\int_{{\mathcal M}\times{{\mathcal
    S}^{n({\bf p})}}}
{\frac{w(N;\sum_{j=1}^{n({\bf p})}\delta_{s_j}){\mathcal P}(dN|\rho,\eta)\prod_{j=1}^{n(\bf
      p)}{h(s_{j})}^{e_{j,n}}\rho(ds_{j})}
{{\(T+\sum_{j=1}^{n(\p)}h({s_j})\)}^{n}}}. 
\label{EPPFw}
\EndEq
Moreover for every $n$, the joint distribution of the unique values of
$\{Y_1,\ldots,Y_n\}|{\bf p}$, that is  $\{Y^{*}_{1},\ldots, Y^{*}_{n({\bf p})}\}$, are
now iid $H$. If $W(N):=g(T)$ then the EPPF is;
\Eq
\int_{{\rr}^{+}\times{{\mathcal S}^{n({\bf p})}}}
\[
{\frac{g(t+\sum_{j=1}^{n({\bf p})}s_j)f_{T}(t)dt\prod_{j=1}^{n(\bf
      p)}{h(s_{j})}^{e_{j,n}}\rho(ds_{j})}
{{\(t+\sum_{j=1}^{n(\p)}h({s_j})\)}^{n}}}\].
\label{EPPFpitman}
\EndEq 
When $h(s)=s\in (0,\infty)$ and $g(T)=1$, the expression in~\mref{EPPFpitman}
corresponds exactly with equation~(31) of Pitman~(1995b).
\end{cor} 

\Remark
One could simply apply a Fubini argument to Theorem 8.1
to deduce the existence of the relevant joint distributions. However,
without any interpretation gained via Pitman, Perman and Yor~(1992) and
Pitman (1995b) the result is somewhat vacuous. Again some of those
interpretations may also be of interests to practicing Bayesians
statisticians as it
certainly goes beyond the usual mean/variance assessment to compare
different random probability measures.
\EndRemark

\Remark
The results above serve also to add to the explicit formulae given
in Pitman (1995b) for the length biased case. 
Note that one could apply~\mref{gammaid} to obtain alternate
representations for the EPPF as in Pitman (1995b, Corollary 6). An
additional interesting feature of Corollary 8.1 is the conditional independence
result in~\mref{djy}. 
\EndRemark

\subsection{Mixture models}
Theorem 8.1 and its corollary  can certainly handle structures of
the form 
\Eq
{\mathcal{Q}}(dN|\rho,\eta)\prod_{i=1}^{n}\int_{\mathcal {Y}}K(X_i|Y_i)P(dY_i)
\label{ssm}
\EndEq
A description of the relevant posterior laws is 
presented below,

\begin{thm}
The posterior
distribution of ${\bf Y},P|{\bf X}$ based on the
model~\mref{likelihood} is representable as, 
$$
\PR(d{\bf Y},dP|{\bf X})\propto 
{\mathcal{P}}(dP|{\bf Y})
\[\prod_{i=1}^{n(\p)}\PR (dY^{*}_{j}|{\bf {\mbox {\~J}}},{\bf p},{\bf
  X})\]
\PR({\bf {\mbox {\~J}}}\in d{\bf s},{\bf
  p}|{\bf X}),
$$
where a conditional distribution of ${\bf Y}
|{\bf {\mbox {\~J}}},{\bf p},{\bf X}$ is 
given such that the sequence
$\{Y_1,\ldots,Y_n\}$ consists of 
$n({\bf p} )$ unique values ${\bf Y}^{*}=\{Y_{1}^{*},\ldots ,Y_{n(\p
)}^{*}\}$ which are independent with respective distributions,  
\Eq
\PR (dY^{*}_{j}|{\bf {\mbox {\~J}}},{\bf p},{\bf X})
={\frac{\[\prod_{i\in C_{j}}K_i(X_i|Y^{*}_{j})\]
\PR (dY^{*}_{j}|{\js}_{j},{\bf p})}{
\int_{\mathcal Y}\[\prod_{i\in C_{j}}K_i(X_i|Y^{*}_{j})\]
\PR (dY^{*}_{j}|{\js}_{j},{\bf p})}
}.
\label{postY}
\EndEq 
for $j=1,\ldots,n({\bf p})$.
In addition 
$
\PR({\bf {\mbox {\~J}}}\in d{\bf s},{\bf
  p}|{\bf X})
\propto
\PR({\bf {\mbox {\~J}}}\in d{\bf s},{\bf p}){\prod_{j=1}^{n(\p)}}\int_{\mathcal{Y}}
\[\prod_{i\in
    C_j}K(X_{i}|Y^{*}_j)\]\PR(dY^{*}_{j}|{\js}_{j},{\bf p})
$
\end{thm}

It follows that when $\rho$ does not depend on $y$~\mref{postY} does
not depend on {\textbf {\~ J}},${\bf T}$. Hence in that case posterior
characterizations of $P,{\bf Y}$ do not differ in form from the results in Lo
(1984) and Ishwaran and James (2001a). The general case however is a
different mattter entirely. 

\Section{Acknowledgements}

I wish to express my deep gratitude to Professor Jim Pitman for
encouraging me to look at this problem in the Poisson generality in which it
appears. Without his gentle persistent nudging and advice this
manuscript may not have been written. 
I thank him for his mentorship despite the fact
that at the commencement of this project we had never met. 
The essence of the technique and style which I use, I have learned from my teacher  
Albert Y. Lo. I thank him  for his steady influence and encouragement
and always reminding me to think simply.
I also mention that the general style of Fubini
calculus presented here was developed over conversations between Albert
Lo and  Lucien Le Cam, of which I have benefitted from. I wish to
thank Dr. C.K. Lim for inspiring me. I wish to thank also
Professors Kallenberg, Paulauskas and Talagrand  for fielding my
somewhat vague questions on the existence of a L{\'e}vy-Khinchine result in
quite abstract spaces. I hope they understand now that I was curious
if there was such a general analogue of Lemma 2.1. and what one might
do with it.

\vskip0.2in

\centerline{\Heading References}
\vskip0.4in
\tenrm
\def\smc{\tensmc}
\def\sl{\tensl}
\def\bf{\tenbold}
\baselineskip0.15in

\Ref
\by    Aalen, O. O.
\yr    1975
\book  Statistical inference for a family of counting processes. Ph.D. Thesis
\publ  University of California
\publaddr Berkeley
\EndRef

\Ref
\by    Aalen, O. O.      
\yr    1978
\paper Nonparametric inference for a family of counting processes
\jour  \AnnStat
\vol   6
\pages 535-545
\EndRef

\Ref
\by    Aalen, O. O.      
\yr    1992
\paper Modelling hetoregeneity in survival analysis by the compund
       Possion distribution
\jour  Annals of Applied Probability
\vol   2
\pages 951-972
\EndRef

\Ref
\by    Aldous, D. J.
\yr    1985
\paper Exchangeability and related topics.  In {\it \'Ecole d'\'Et\'e de
       Probabilit\'es de Saint-Flour XII} (P. L. Hennequin, editor).
       Springer Lecture Notes in Mathematics, Vol. 1117
\EndRef

\Ref
\by    Andersen, P. K., Borgan, O. , Gill, R. D. and  Keiding, N.
\yr    1993
\book  Statistical Models Based On Counting Processes
\publ  Springer-Verlag
\publaddr New York
\EndRef

\Ref   
\by    Antoniak, C. E. 
\yr    1974
\paper Mixtures of Dirichlet processes with applications to Bayesian
nonparametric problems
\jour  \AnnStat
\vol   2
\pages 1152-1174
\EndRef

\Ref
\by     Banjevic, D., Ishwaran, H. and Zarepour, M.  
\yr     2002
\paper  A recursive method for functionals of Poisson processes. To
        appear in {\it Bernoulli}
\EndRef

\Ref
\by    Bar-Lev, S.K. and Enis, P. 
\yr    1986    
\paper Reproducibility and natural exponential families with power
       variance functions
\jour  \AnnStat
\vol   14
\pages 1507-1522
\EndRef

\Ref   
\by    Barndorff-Nielsen, O.E. and Shephard, N.
\yr    2001
\paper Normal modified stable processes. Preprint
\EndRef

\Ref
\by    Berk, R. and Savage , I. R.
\yr    1979
\paper Dirichlet processes
  produce discrete measures: an elementary proof. \textit{Contributions 
    to Statistics, Jaroslav Hajek Memorial Volume.} Edited by Jana
  Jure\v ckova.  Reidel, Dordrecht-Boston, Mass.-London. pp. 25-31 
\EndRef

\Ref
\by      Bertoin, J.
\yr       2001
\paper    Homogeneous fragmentation processes 
\jour     Probab. Theory Related Fields 
\vol      121
\pages   301--318
\EndRef

\Ref
\by    Blackwell, D.
\yr    1973
\paper Discreteness of Ferguson selections
\jour  \AnnStat
\vol   1
\pages 356-358
\EndRef

\Ref
\by    Blackwell, D. and MacQueen, J. B.
\yr    1973
\paper Ferguson distributions via P\'olya urn schemes
\jour  \AnnStat
\vol   1
\pages 353-355
\EndRef

\Ref
\by    Blackwell, D. and Maitra A.
\yr    1984
\paper Factorization of probability measures and absolutely measurable
sets
\jour  Proc. Amer. Math. Soc.
\vol   92
\pages 251-254
\EndRef

\Ref   
\by    Blum, J. and Susarla, V. 
\yr    1977
\paper On the posterior distribution
of a Dirichlet process given random right censored data 
\jour  Stochastic Process. Appl
\vol   5
\pages  207-211
\EndRef

\Ref
\by    Brix, A.
\yr    1999
\paper Generalized Gamma measures and shot-noise Cox processes
\jour  Adv. in Appl. Probab. 
\vol   31 
\pages 929-953
\EndRef

\Ref
\by    Brunner, L. J., Chan, A. T., James, L. F. and Lo, A. Y. 
\yr    2001
\paper Weighted Chinese restaurant processes and Bayesian mixture models.
       preprint
\EndRef

\Ref
\by    Brunner, L. J. and Lo, A. Y.
\yr    1989
\paper  Bayes methods for a symmetric unimodal density and its mode 
\jour  \AnnStat
\vol   17
\pages 1550-1566
\EndRef   

\Ref
\by     Carlton, M.A. 
\yr     1999
\paper   Applications of the Two-Parameter Poisson Dirichlet
distribution. Ph.D. Thesis. Univeristy of California, Los
Angeles. Dept. of Statistics
\EndRef

\Ref   
\by    Cifarelli, D.M. and Regazzini, E. 
\yr    1990
\paper Some remarks on the distribution of the means of a Dirichlet
       process
\jour  \AnnStat
\vol   18
\pages 429-442
\EndRef

\Ref    
\by      Constantine, G. M. 
\yr      1999
\paper   Identities over set partitions
\jour    Discrete Mathematics
\vol     204
\pages   155-162
\EndRef

\Ref
\by      Constantine, G. M. and Savits, T. H. 
\yr      1994
\paper   A stochastic process interpretation of partition identities
\jour    Siam J. Discrete Math.
\vol     7
\pages   194-202
\EndRef

\Ref
\by        Daley, D. J. and Vere-Jones, D.
\yr        1988
\book      An Introduction to the Theory of Point Processes
\publ      Springer-Verlag
\publaddr  New York
\EndRef

\Ref
\by      Dellacherie, C and Meyer, P. A.
\yr      1978
\book    Probabilities and Potential
\publ    North Holland
\publ    Amsterdam
\EndRef

\Ref
\by      Dey, J. 
\yr      1999
\paper   Some Properties and Characterizations of Neutral-to-the-Right
Priors and Beta Processes. Ph.D. Thesis.  Michigan State University
\EndRef

\Ref
\by       Dey, J, Erickson, R.V. and Ramamoorthi, R.V.
\yr       2000
\paper    Neutral to right priors- A review. Preprint
\EndRef 

\Ref
\by        Di Nardo, E. and Senato, D.
\yr        2001
\paper     Umbral nature of the Poisson random variables. In {\it
  Algebraic Combinatorics Computer Science: A Tribute to Gian-Carlo
  Rota} Springer 245-267
\EndRef

\Ref
\by    Diaconis, P. and Kemperman, J. 
\yr    1996
\paper Some new tools for Dirichlet priors. Bayesian Statistics 5 
       (J.M. Bernardo, J.O. Berger, A.P. Dawid and A.F.M. Smith eds.), 
       Oxford University Press, pp. 97-106
\EndRef 

\Ref
\by    Doksum, K. A.
\yr    1974
\paper Tailfree and neutral random probabilities
and their posterior distributions
\jour  Ann. Probab
\vol   2
\pages 183-201
\EndRef

\Ref
\by    Donnelly, P. and Tavar\'e, S.
\yr    1987
\paper The population genealogy of the infinitely-many neutral 
       alleles model
\jour  J. Math. Biol. 
\vol   25
\pages 381-391
\EndRef

\Ref 
\by        Doss, H. 
\yr        1994
\paper     Bayesian nonparametric estimation for
  incomplete data via successive substitution sampling
\jour      \AnnStat
\vol       22
\pages     1763-1786
\EndRef

\Ref
\by    Dykstra, R. L. and Laud, P. W.
\yr    1981
\paper A Bayesian nonparametric approach to reliability
\jour  \AnnStat
\vol   9
\pages 356-367
\EndRef

\Ref
\by    Engen, S.
\yr    1978
\book  Stochastic Abundance Models with Emphasis on Biological
Communities and Species Diversity
\publ  Chapman and Hall
\EndRef

\Ref
\by    Ewens, W. J.
\yr    1972
\paper The sampling theory of selectively neutral alleles
\jour  Theor. Popul. Biol.
\vol   3
\pages 87-112
\EndRef

\Ref
\by    Ewens, W. and Tavar\'e, S.
\yr    1997
\paper Multivariate Ewens distribution.  In {\it Discrete Multivariate
       Distributions} (S. Kotz and N. Balakrishnan, eds.).  Wiley, New
       York
\EndRef

\Ref
\by    Ferguson, T. S.
\yr    1973
\paper A Bayesian analysis of some nonparametric problems
\jour  \AnnStat
\vol   1
\pages 209-230
\EndRef

\Ref
\by    Ferguson, T. S.
\yr    1974
\paper Prior distributions on spaces of probability measures
\jour  \AnnStat
\vol   2
\pages 615-629
\EndRef

\Ref
\by     Ferguson, T. S. and Klass, M. J. 
\yr     1972
\paper  A representation of independent increment processes without 
Gaussian components 
\jour   Ann. Math. Statist.
\vol    43
\pages  1634-1643
\EndRef

\Ref    
\by     Ferguson, T. S. and Phadia, E.
\yr     1979
\paper  Bayesian nonparametric estimation based on censored data
\jour   \AnnStat
\vol    7
\pages  163-186
\EndRef

\Ref   
\by     Freedman, D. A.
\yr     1963
\paper  On the asymptotic behaviour of
  Bayes estimates in the discrete case
\jour   Ann. Math. Statist
\vol    34
\pages  1386-1403
\EndRef

\Ref
\by      Fitzsimmons, P., Pitman, J. and Yor, M. 
\yr      1992
\paper   Markovian bridges: construction, Palm interpretation and
splicing, In {\it Seminar on stochastic process 1992} Editors Cinlar,
Chung, Sharpe. Birkhauser, Boston. 101-135. 
\EndRef 

\Ref
\by      Gill, R. D. and Johansen, S. 
\yr      1990 
\paper   Survey of product-integration with a view towards
applications in survival analysis 
\jour    \AnnStat
\vol     18
\pages   1501-1555
\EndRef

\Ref
\by     Groeneboom, P.
\yr     1996
\paper  Lectures on inverse problems. Lectures on probability theory and statistics (Saint-Flour, 1994), 67--164, Lecture Notes in Math., 1648, Springer, Berlin
\EndRef

\Ref
\by     Groeneboom, P. and  Wellner, J. A. 
\yr     1992
\paper  Information bounds and nonparametric maximum likelihood estimation. DMV Seminar, 19. Birkhauser Verlag, Basel
\EndRef

\Ref   
\by    Gyllenberg, M, and Koski, T. 
\yr    2001
\paper Probabilistic models for bacterial taxonomy
\jour  International Statistical Review
\vol   69
\pages 249-276
\EndRef

\Ref
\by    Hansen, B. and Pitman, J.
\yr    2000
\paper Prediction rules for exchangeable sequences related to
       species sampling
\vol   46
\jour  Statist. Prob. Letters
\pages 251-256
\EndRef

\Ref      
\by      Hjort, N. L.
\yr      1990
\paper   Nonparametric Bayes estimators based on Beta processes in
models for life history data 
\jour    \AnnStat
\vol         18
\pages    1259-1294
\EndRef

\Ref   
\by    Hougaard, P. 
\yr    1986
\paper Survival models for heterogeneous populations derived from
       stable distributions
\jour  Biometrika
\vol   73
\pages 387-396
\EndRef

\Ref   
\by    Hougaard, P., Lee, M. L. and Whitmore, G. 
\yr    1997
\paper Analysis of overdispersed count data by mixtures of Poisson
       variables and Poisson processes 
\jour  Biometrics
\vol   53
\pages 1225-1238
\EndRef

\Ref
\by      Huang, Y. and Louis, T. A. 
\yr      1998
\paper   Nonparametric estimation of
the joint distribution of survival time and mark variables 
\jour    Biometrika
\vol     85
\pages   785-798
\EndRef

\Ref
\by    Ishwaran, H. and James, L. F. 
\yr    2001a
\paper Generalized weighted Chinese restaurant processes for species
       sampling models. Manuscript
\EndRef

\Ref     
\by      Ishwaran, H. and James, L. F. 
\yr      2001b
\paper   Gibbs sampling methods for stick-breaking priors 
\jour    J. Amer. Stat. Assoc
\pages   161-173 
\EndRef

\Ref
\by     Ishwaran, H., James, L. F. and Sun, J. 
\yr     2001   
\paper  Bayesian model selection in finite mixtures by marginal
density decompositions  
\jour   J. Amer. Stat. Assoc 
\vol    96
\pages  1316-1332
\EndRef

\Ref   
\by    Jacod, J. 
\yr    1975
\paper      Multivariate point processes: predictable projection,
Radon-Nikodym derivatives, representation of martingales
\jour  \Zfw
\vol   35
\pages 1-37
\EndRef

\Ref   
\by    James, L. F. and Kwon, S.
\yr    2000
\paper A Bayesian nonparametric  approach for the joint distribution
       of survival time and mark variables under univariate censoring.
       Unpublished manuscript
\EndRef

\Ref
\by    James, L. F. 
\yr    2001a
\paper Bayesian calculus for Gamma processes with applications to
       semiparametric models. To appear in {\it Sankhy\=a Ser. A}
\EndRef 

\Ref
\by    James, L. F.
\yr    2001b
\paper An analysis of weighted generalised gamma process mixture
models. Unpublished notes
\EndRef

\Ref
\by    Jorgensen, B.
\yr    1997
\book  The Theory of Dispersion Models. 
       Monographs on Statistics and Applied Probability, 76. 
\publ  Chapman and Hall
\publaddr London 
\EndRef

\Ref
\by    Kallenberg, O. 
\yr    1986
\book  Random Measures, 4th Ed
\publ  Akademie-Verlag and Academic Press
\publaddr Berlin and London
\EndRef

\Ref
\by    Kallenberg, O.
\yr    1997
\book   Foundations of modern probability. Probability and its
Applications 
\publ  Springer-Verlag
\publaddr New York
\EndRef

\Ref
\by    Kaplan, E. L. and Meier, P. 
\yr    1958
\paper  Nonparametric estimation from incomplete observations
\jour   J. Amer. Statist. Assoc. 
\vol    53 
\pages  457-481
\EndRef
 
\Ref
\by    Kerov, S. 
\yr    1998
\paper Interlacing measures
\jour  Amer. Math. Soc. Transl.
\vol   181
\pages 35-83
\EndRef

\Ref
\by    Kerov, S. and Tsilevich, N. V.
\yr    1998
\paper The Markov-Krein correspondence in several dimensions.
POMI preprint No. 283, Steklov Institute of Mathematics, St. Petersburg
\EndRef

\Ref
\by    Khintchine, A. Ya. 
\yr    1937
\paper Sur theorie der unbeschrankt teilbaren Verteilungsgesetze
\jour  Mat. Sb.
\vol 44
\pages 79-119
\EndRef

\Ref
\by    Kim, Y.
\yr    1999
\paper  Nonparametric Bayesian estimators for counting processes 
\jour   \AnnStat
\vol    27 
\pages  562-588
\EndRef

\Ref
\by    Kingman, J. F. C.
\yr    1993
\book  Poisson Processes
\publ  Oxford University Press
\publaddr Oxford
\EndRef

\Ref
\by    Kingman, J. F. C.
\yr    1975
\paper Random discrete distributions
\jour  \JRSSB
\vol   37
\pages 1-22
\EndRef

\Ref
\by    Kingman, J. F. C. 
\yr    1967
\paper Completely random measures
\jour  Pacific J. Math. 
\vol   21
\pages 59-78
\EndRef

\Ref
\by    K{\"{u}}chler, U. and Sorenson, M.
\yr    1997
\book  Exponential Families of Stochastic Processes
\publ  Springer-Verlag
\publaddr New York
\EndRef

\Ref
\by    Le Cam, L.
\yr    1986
\book  Asymptotic Methods in Statistical Decision Theory 
\publ  Springer-Verlag
\publaddr New York
\EndRef

\Ref
\by    Le Cam, L. 
\yr    1961 
\paper  A stochastic description of precipitation. In {\it 1961
  Proc. 4th Berkeley Sympos. Math. Statist. and Prob.}, Vol. III
Univ. California Press, Berkeley Calif.  165-186
\EndRef

\Ref
\by    Lee, M.L.T. and Whitmore, G. 
\yr    1993
\paper Stochastic processes directed by randomized time
\jour  Journal of Applied Probability
\vol   30
\pages 302-314
\EndRef
  
\Ref
\by    Last, G. and Brandt, A.
\yr    1995
\book  Marked Point Proceses on the Real Line: The Dynamic Approach
\publ  Springer
\publaddr New York
\EndRef

\Ref
\by    Lindsay, B.
\yr    1995
\book  Mixture models: theory, geometry, and applications. NSF-CBMS
       Regional Conference Series in Probability and Statistics, Volume 5 
\publ  Institute for Mathematical Statistics: Hayward, CA
\publaddr Hayward, CA
\EndRef

\Ref
\by    Lo, A. Y.
\yr    1982
\paper Bayesian nonparametric statistical inference for Poisson point
       processes 
\jour  \Zfw
\vol   59
\pages 55-66
\EndRef

\Ref
\by    Lo, A. Y.
\yr    1984
\paper On a class of Bayesian nonparametric estimates: I. Density estimates
\jour  \AnnStat
\vol   12
\pages 351-357
\EndRef

\Ref
\by    Lo, A. Y.
\yr    1993
\paper A Bayesian bootstrap for censored data 
\jour  \AnnStat
\vol    21 
\pages  100--123
\EndRef

\Ref
\by    Lo, A. Y. and Weng, C. S.
\yr    1989
\paper On a class of Bayesian nonparametric estimates: II. Hazard
       rates estimates 
\jour  Ann. Inst. Stat. Math
\vol   41
\pages 227-245
\EndRef

\Ref
\by    Lo, A.Y., Brunner, L.J. and Chan, A.T.
\yr    1996
\paper Weighted Chinese restaurant processes and Bayesian mixture
       model. Research Report Hong Kong University of
       Science and Technology
\EndRef

\Ref
\by    MacEachern, S. N., Clyde, M. and Liu, J. S.
\yr    1999
\paper Sequential importance sampling for nonparametric Bayes models:
       the next generation
\jour  Canadian J. Statist.
\vol   27
\pages 251-267
\EndRef

\Ref
\by    Matthes, K. Kerstan, J., and Mecke, J. 
\yr    1978
\book  Infinitely Divisible Point Processes. English Edition
\publ  Wiley
\publaddr Chichester
\EndRef

\Ref
\by     McCloskey, J. W. 
\yr     1965
\paper  A Model for the Distribution of Individuals by Species in an 
        Environment. Ph.D. Thesis. Michigan State University
\EndRef

\Ref    
\by    Pachl, J. K.
\yr    1978
\paper Disintegration and compact measures
\jour  Math. Scand.
\vol   43
\pages 157-168
\EndRef

\Ref
\by    Perman, M., Pitman, J. and Yor, M. 
\yr    1992
\paper Size-biased sampling of Poisson point processes and
       excursions
\jour  Probab. Theory Related Fields 
\vol   92
\pages  21-39
\EndRef
 
\Ref  
\by    Pitman, J.
\yr    1995a
\paper Exchangeable and partially exchangeable random partitions
\jour  \newZfw
\vol   102
\pages 145-158
\EndRef

\Ref
\by    Pitman, J. 
\yr    1995b
\paper Poisson-Kingman partitions
\jour  Available at www.stat.berkeley.edu/users/pitman
\EndRef

\Ref  
\by    Pitman, J.
\yr    1996
\paper Some developments of the Blackwell-MacQueen urn scheme. In {\it
       Statistics, Probability and Game Theory} (T.S. Ferguson,
       L.S. Shapley and J.B. MacQueen, eds.) 245-267.  IMS Lecture
       Notes-Monograph series, Vol 30
\EndRef

\Ref
\by    Pitman, J. 
\yr    1997a
\paper Partition structures derived from Brownian motion and stable 
       subordinators 
\jour  Bernoulli
\vol   3 
\pages 79-96 
\EndRef

\Ref
\by    Pitman, J. 
\yr    1997b
\paper Some probabilistic aspects of set partitions 
\jour   Amer. Math. Monthly, 
\pages  201-209
\EndRef

\Ref
\by    Pitman, J.
\yr    1999
\paper Coalescents with multiple collisions
\jour  Ann. Probab 
\vol   27 
\pages 1870-1902
\EndRef

\Ref
\by    Pitman, J. and Yor, M. 
\yr    1992
\paper Arcsine laws and interval partitions derived from a stable 
       subordinator 
\jour  Proc. London Math. Soc. 
\vol   65
\pages 326--356 
\EndRef

\Ref
\by    Pitman, J. and Yor, M.
\yr    1997
\paper The two-parameter Poisson-Dirichlet distribution derived from
       a stable subordinator
\jour  \AnnProb
\vol   25
\pages 855-900
\EndRef

\Ref
\by    Pitman, J. and Yor, M. 
\yr    2001
\paper On the distribution of ranked heights of excursions of a
       Brownian  bridge 
\jour  \AnnProb 
\vol   29
\pages 361-384
\EndRef

\Ref   
\by    Pollard, D.
\yr    2001
\book  User's Guide to Measure Thoeretic Probability
\publ  Cambridge University Press
\EndRef

\Ref
\by    Rosinski, J. 
\yr    2001
\paper Series representations of L{\'e}vy processes from
the perspective of point processes. In {\it L{\'e}vy Processes: Theory and
   Applications.} Eds. Brandorff-Nielsen, Mikosch, and Resnick, pp. 401-415.
 Birkhauser, Boston.
\EndRef

\Ref
\by    Rota, G-C.
\yr    1964
\paper The number of partitions of a set
\jour  American Mathematical Monthly
\vol   71
\pages 498-504
\EndRef

\Ref
\by     Sato, K.
\yr     1999
\book   L{\'e}vy Processes and Infinitely Divisible proceses
\publ   Cambridge University Press
\publaddr Cambridge, UK
\EndRef

\Ref
\by     Susarla, V. and Van Ryzin, J. 
\yr     1976
\paper   Nonparametric Bayesian estimation of survival curves from
incomplete observations
\jour    J. Amer. Statist. Assoc
\vol     71
\pages   897-902
\EndRef
 
\Ref
\by    Tsilevich, N. V.
\yr    1997
\paper Distribution of the mean value for certain random measures.
       POMI preprint No. 240, Steklov Institute of Mathematics,
       St. Petersburg. English translation in Journal of Mathematical
       Sciences, vol. 96 (1999), No. 5, pp. 3616-3623.
\EndRef

\Ref
\by    Tsilevich, N. V., Vershik, A. M, and Yor, M.
\yr    2001
\paper An infinite-dimensional analogue of the Lebesque measure and
       distinguished properties of the gamma process 
\jour   J. Funct. Anal 
\vol    185 
\pages   274-296
\EndRef

\Ref
\by    Tsilevich, N. V., Vershik, A. M, and Yor, M.
\yr    2000
\paper Distingusihed properties of the gamma process and related
       topics, Pr\'epublication du Laboratoire de Probabilit\'es et Mod\`eles 
       Al\'eatoires no. 575, Mars 2000 
\EndRef

\Ref
\by    Turnbull, B. W.
\yr    1976
\paper The empirical distribution function with arbitrarily grouped,
censored and truncated data 
\jour   J. Roy. Statist. Soc. Ser. B 
\vol    38
\pages  290-295
\EndRef
  
\Ref
\by    Tweedie, M.C.K.  
\yr    1984
\paper An index which distinguishes between some important exponential
       families. In {\it Statistics: Applications and New
       Directions} Eds. J.K. Ghosh and J. Roy, pp. 579-604. Indian
       Statistical Institute, Calcutta
\EndRef 

\Ref
\by    Walker, S. and Muliere, P.
\yr    1997
\paper Beta-Stacy processes and a generalization of the P{\'o}lya-urn
scheme
\jour   \AnnStat
\vol    25 
\pages  1762-1780
\EndRef

\Ref   
\by    Wolpert, R. L. and Ickstadt, K.
\yr    1988a
\paper Poisson/Gamma random field models for spatial statistics
\jour  Biometrika
\vol   85
\pages 251-267
\EndRef

\Ref
\by    Wolpert,  R. L. and  Ickstadt, K.
\yr    1998b
\paper Simulation of L{\'e}vy random fields. In {\it Practical
       Nonparametric and Semiparametric Bayesian Statistics} (D. Dey,
       P. Mueller and D. Sinha, eds.) 227-241. Springer Lecture Notes
\EndRef

\vskip0.75in

\smc

\Tabular{ll}

Lancelot F. James\\
The Hong Kong University of Science and Technology\\
Department of Information Systems and Management\\
Clear Water Bay, Kowloon\\
Hong Kong\\
\rm lancelot\at ust.hk\\

\EndTabular

\end{document}